\documentclass[10pt,reqno]{amsart}
\usepackage{amssymb,mathrsfs,graphicx}
\usepackage{placeins,cite}
\usepackage[a4paper,textwidth=170mm,textheight=235mm,centering]{geometry}
\numberwithin{equation}{section}

\usepackage[colorlinks=true,linkcolor=blue,citecolor=blue,urlcolor=blue,hyperfootnotes=false]{hyperref}

\newtheorem{theorem}{Theorem}[section]
\newtheorem{proposition}[theorem]{Proposition}%

\newtheorem{example}{Example}%
\newtheorem{assumption}{Assumption}

\newtheorem{lemma}{Lemma}[section]
\newcommand{\E}{\mathbb{E}}
\newcommand{\PP}{\mathbb{P}}
\newcommand{\RR}{\mathbb{R}}

\begin{document}
	\title[Domain-preserving projected Milstein scheme]
{An explicit domain-preserving projected Milstein scheme with strong order one for stochastic systems}

\author[X. Hu]{Xingwei Hu\textsuperscript{1,2}}
\author[Y. Hu]{Yaozhong Hu\textsuperscript{2}}
\author[A. Xiao]{Aiguo Xiao\textsuperscript{1,3,*}}

\address{\textsuperscript{1} School of Mathematics and Computational Science,
Xiangtan University, Hunan 411105, China}
\address{\textsuperscript{2} Department of Mathematical and Statistical Sciences,
University of Alberta, Edmonton, Alberta T6G 2G1, Canada}
\address{\textsuperscript{3} Hunan Research Center of the Basic Discipline Fundamental
Algorithmic Theory and Novel Computational Methods \&
China National Center for Applied Mathematics in Hunan,
Xiangtan University, Hunan 411105, China}
\email[Aiguo Xiao]{xag@xtu.edu.cn}
\thanks{Corresponding author: Aiguo Xiao.}

\begin{abstract} 
		A variety of domain-preserving numerical schemes have been developed for scalar stochastic differential equations and for several classes of multidimensional systems, including those defined on hyperrectangles and positive orthants. To accommodate more general domain constraints, we consider stochastic systems on bounded or unbounded open convex polyhedral domains, where coefficients may exhibit boundary singularities and polynomial growth. Under local regularity, growth and boundary assumptions, we establish well-posedness and integrability estimates for the exact solution. Then we construct an explicit projected Milstein scheme that preserves the domain by applying a Euclidean projection onto a step-size-dependent interior set after each Milstein update. For systems driven by multidimensional commutative noise, we prove strong order-one convergence for the maximum grid error in $L^p$, for suitable $p > 1$ and sufficiently small step sizes under additional weighted regularity and stability conditions. The proof combines these integrability estimates with local Milstein consistency, localization and weighted stability, without requiring a priori step-size-uniform moment or boundary inverse-moment bounds for the numerical approximation. Applications to constrained stochastic models and numerical experiments illustrate the observed first-order convergence behavior.
	\end{abstract}


	\maketitle

	\section{Introduction}
	\label{sec.1}
	Stochastic differential equations (SDEs) are widely used to model dynamical systems subject to random perturbations. In many applications, their solutions must remain in a prescribed domain, such as a positive orthant, a bounded interval, a simplex, or a more general convex domain. Closed-form solutions are rarely available, so numerical approximations are needed. These approximations should preserve the prescribed domain to retain the model’s interpretation and keep coefficient evaluations well defined. In general, standard explicit schemes fail to preserve the domain, as they may exit it with positive probability. 
	
	For SDEs with non-globally Lipschitz coefficients, taming
	\cite{HJK2012,WG2013,S2013,S2016} and truncation \cite{M15,M16} have been used to control coefficient growth and construct strongly convergent explicit schemes. Convergence analyses under such conditions have also been developed for modified It\^o--Taylor schemes \cite{JKN2009}, fully implicit schemes \cite{TZ2013} and implicit Milstein schemes \cite{W2023}. Projected Euler--Maruyama and Milstein schemes have been analyzed in \cite{BIK2016,BIK2017} using stochastic C-stability and B-consistency under global monotonicity and polynomial regularity assumptions. In these schemes, the current approximation is projected onto a ball whose radius increases as the step size decreases, before applying the time-stepping update. Notably, these results concern equations with coefficients defined on the entire space. However, for equations posed on a prescribed domain, additional analysis is needed to establish domain preservation.
	
	Domain-preserving approximations have been extensively studied
	for scalar SDEs \cite{AA2005,AA2013,SNS2012,YH2021,CGW2021,LW2023,C2023,T2025,WG2025,DFWM2023,MS2007,brc26,HHX20262}. Transformation-based schemes form an important part of this literature. Lamperti transformations can reduce scalar SDEs to equations with additive noise \cite{CJI2016,NS2014}, while logarithmic and related transformations can map the original interval onto the real line \cite{PPlogTM,LG2023,TX2024,HXW,HHX2027}. \cite{NS2014} combined a Lamperti transformation with an implicit discretization and established strong order-one convergence for a domain-preserving scheme. \cite{LNX2025} combined this transformation with semi-implicit discretization and taming to obtain a boundary-preserving scheme with strong order $1.5$ under additional assumptions. 
	
	Extending scalar approaches to multidimensional systems requires accounting for interactions between components, the geometry of the domain and the structure of the noise. Domain-preserving schemes have nevertheless been developed for several classes of stochastic systems
	\cite{MF21,APPL2023,Niko2015,EL2026,CWF2024,hxw2026}, with related work on epidemic systems with delays \cite{THX2026}. In particular, \cite{DW2026} proposed a positivity-preserving explicit Milstein-type scheme for stochastic Lotka--Volterra competition systems driven
	by multidimensional noise and established strong order-one convergence. More recently, \cite{br26} developed explicit domain-preserving schemes for stochastic systems on hypercubes. Their general scheme attains strong order $\frac 12$ and weak order one, while a modified version achieves strong order one for systems driven by a scalar Brownian motion. These results raise the following question: Can domain-preserving numerical schemes be constructed for stochastic systems on more general domains?
	
	In this paper, we consider stochastic systems on bounded or
	unbounded convex polyhedral domains. The coefficients may exhibit boundary singularities and polynomial growth. Specifically, the system takes the form
	\begin{flalign}\label{ssde}
		\mathrm{d}X_t = f(X_t) \mathrm{d}t + \sum_{l=1}^m g_l(X_t)\,\mathrm dB_t^l,\quad t\in(0,T],\quad X_0=x_0\in \mathcal{D},
	\end{flalign}
	where $x_0$ is deterministic, $B_t$ is a standard  $m$-dimensional Brownian motion and $\mathcal{D}\subset \RR^d$ is a non‑empty open convex polyhedral domain of the form
	\begin{flalign*}
		\mathcal{D} = \{x \in \mathbb{R}^d : \phi_j(x) > 0, \, j = 1,2, \dots, J\}, \quad \phi_j(x) = c_j + \nu_j^\top x, \quad |\nu_j| = 1.
	\end{flalign*}
	Here, for each $j=1,2,\ldots,J$, $c_j\in\mathbb R$ and $\nu_j\in\mathbb R^d$ are fixed. The vector $\nu_j$ is normal to the hyperplane $\{x:\phi_j(x)=0\}$ and points into the half-space $\{x:\phi_j(x)>0\}$. In particular, $\mathcal{D}$ may be an open hyperrectangle, a positive orthant, or the interior of a simplex. This setting also allows domains such as $(0,1) \times \mathbb{R}$, where one coordinate is bounded and the other is unbounded. For a sufficiently small step size \(h=\frac TN\), we introduce a closed convex set \(\mathcal D_h\subset\mathcal D\) containing \(x_0\). This step-size-dependent set keeps the numerical approximation a positive distance from the boundary and bounds its size when $\mathcal{D}$ is unbounded. For systems driven by multidimensional commutative noise,
	the proposed projected Milstein scheme takes the form
	\begin{flalign*}
		Y_{n+1} =\Pi_{\mathcal D_h}\bigl(\mathcal M_{h,n}(Y_n)\bigr),\quad n=0,1,\ldots,N-1, \quad Y_0=x_0,
	\end{flalign*}
	where \(\mathcal M_{h,n}\) is the standard commutative Milstein update and \(\Pi_{\mathcal D_h}\) denotes the Euclidean projection onto $\mathcal{D}_h$ (see Section \ref{sec3.1}). The projection ensures that the numerical approximation
	remains in the prescribed domain at every time step. 
	
	The exact solution may enter the boundary layers excluded
	by $\mathcal D_h$ or, on unbounded domains, exceed the imposed
	bounds. The contribution of these events to the approximation
	error must therefore be controlled. Moreover, the lower bounds on the
	boundary distances vanish as $h\to0$, so domain preservation
	alone does not yield step-size-uniform boundary inverse-moment
	estimates for the numerical approximation.
	
	To address these difficulties, we first establish well-posedness
	and moment, boundary inverse-moment and exponential integrability
	estimates for the exact solution under suitable local regularity,
	boundary and growth conditions (see Section \ref{sec2.2}). We then combine these estimates with geometric localization, local Milstein consistency and weighted stability estimates (see Sections \ref{sec3.2} and \ref{sec3.3}). Under suitable weighted regularity and stability assumptions, we establish the following error estimate for appropriately chosen projection parameters:
	\begin{flalign*}
		\mathbb E\left[ \max_{0\le n\le N}|X(t_n)-Y_n|^p \right]\le Ch^p, \quad t_n=nh,
	\end{flalign*}
	for suitable \(p>1\) and sufficiently small \(h\). The analysis does not require a priori step-size-uniform moment or inverse-moment estimates for the numerical approximation.
	
	We apply the theory to representative constrained models
	in population dynamics and epidemiology and specify
	the conditions under which the convergence theorem applies. We also construct an SIR model with population-dependent noise as a test problem for the proposed method. Numerical experiments exhibit convergence rates close to one over the tested step sizes and include a comparison with the scheme of \cite{br26} on a common
	hypercube problem driven by a scalar Brownian motion.
	
	The paper is organized as follows. Section~\ref{sec.2} introduces the notation and preliminary inequalities and establishes well-posedness and integrability of the exact solution. Section~\ref{sec3} presents the projected Milstein scheme and proves its strong order-one convergence. Section~\ref{sec4} discusses applications to representative constrained stochastic models. Section~\ref{sec5} reports numerical experiments, and Section~\ref{sec6} concludes the paper.

	\section{Preliminaries and properties of the exact solution}
	\label{sec.2}
	\subsection{Notation and preliminary inequalities}
	Throughout this paper, let
	$(\Omega,\mathcal F,\{\mathcal F_t\}_{t\ge 0},\mathbb P)$
	be a filtered probability space satisfying the usual conditions.
	Let $B=(B^1,B^2,\ldots,B^m)^\top$ be an $m$-dimensional standard Brownian motion with respect to the filtration $(\mathcal F_t)_{t\geq0}$. We denote expectation by $\mathbb E$. For an event or a set $A$, $\mathbb I_A$ denotes its indicator, and $A^c$ denotes its complement. For $x,y\in\mathbb R^d$, $|x|$ denotes the Euclidean norm and
	$\langle x,y\rangle$ the Euclidean inner product. For a matrix $A$, $\|A\|$ and $\|A\|_F$ denote the induced Euclidean operator norm and the Frobenius norm, respectively, and $A^\top$ denotes its transpose. For a random variable $X$ and $p\ge1$, we write $\|X\|_{L^p}:=(\mathbb E|X|^p)^{\frac 1p}$. For sufficiently smooth scalar or vector functions,
	$D$ and $D^2$ denote the first and second derivatives, respectively. The notation $C_{\rm loc}^k(\mathcal{D};\mathbb R^m)$ denotes the space of functions that are $k$ times continuously differentiable locally on $\mathcal{D}$. For a set $A$, $\overline A$ denotes its closure, and $A\Subset \mathcal{D}$ means that $\overline A$ is a compact subset of $\mathcal{D}$. For $x,y\in\mathbb R^d$, we write $[x,y]:=\{(1-\theta)x+\theta y:\theta\in[0,1]\}$. If $K\subset\mathbb R^d$ is a nonempty closed convex set, $\Pi_K$ denotes the Euclidean projection onto $K$, i.e., $\Pi_K(x):=\arg\min_{y\in K}|x-y|$. In particular, $\Pi_K$ is nonexpansive: $|\Pi_K(x)-\Pi_K(y)|\le |x-y|$. For real numbers $m,n$, we write $m\wedge n:=\min\{m,n\}$ and
	$m\vee n:=\max\{m,n\}$. Let $N\in\mathbb N$, $h=T/N$, and $t_n:=nh$, $n=0,1,\ldots,N$. We set $\mathbb E_n[\cdot] :=\mathbb E[\cdot\mid\mathcal F_{t_n}]$. For a scalar function $\varphi$, we write $D^2\varphi(x)[u,v]=u^\top D^2\varphi(x)v$ for $u,v\in\mathbb R^d$. For a vector-valued function $\varphi=(\varphi_1,\varphi_2,\ldots,\varphi_s)^\top:
	\mathcal D\to\mathbb R^s$, we define
	\[
	D^2\varphi(x)[u,v]
	=
	\begin{pmatrix}
		u^\top D^2\varphi_1(x)v\\
		\vdots\\
		u^\top D^2\varphi_s(x)v
	\end{pmatrix}
	\in\mathbb R^s,
	\quad u,v\in\mathbb R^d.
	\]
	Throughout the paper, $C$ denotes a generic positive constant
	independent of $h$ and $n$, whose value may change from line to line.
	
	\begin{proposition}[Ville's inequality]\label{ville}
		Suppose that $Z_t$ is a nonnegative right-continuous supermartingale. Then it follows that for every $L>0$,
		\begin{flalign*}
			\PP \Big( \sup_{0 \le t \le T} Z_t \ge L \Big) \le\frac{ \E[Z_0]}{L},
		\end{flalign*}
	\end{proposition}
	
	We omit the proof here, as it can be found in \cite{JOb2015,HR2020}.
	
	\begin{proposition}[Discrete Lenglart domination inequality]
		Let $\{Z_n\}_{n=0}^N$ be a nonnegative adapted process and
		$\{A_n\}_{n=0}^N$ a nondecreasing predictable process with $A_0=0$.
		Suppose that $\mathbb E[Z_\tau]\le \mathbb E[A_\tau]$ for every stopping time $\tau\le N$. Then, for every $0<\theta<1$,
		\begin{flalign*}
			\mathbb E\Big[ \Big(\max_{0\le n\le N}Z_n\Big)^\theta \Big] \le C_\theta\mathbb E[A_N^\theta].
		\end{flalign*}
	\end{proposition}
	
	We omit the proof here, as it can be found in \cite{1977,Mehri2021}.
	
	\subsection{Well-posedness and integrability of the exact solution}
	\label{sec2.2}
	Recall the stochastic system \eqref{ssde}. Fix $x^{\circ}\in\mathcal{D}$. We impose the following assumptions.
	\begin{assumption}\label{ass0}
		For $l=1,2,\ldots,m$, $f,g_l\in C^2_{loc}(\mathcal{D},\RR^d)$.
	\end{assumption}
	\begin{assumption}\label{ass1}
		If $\mathcal{D}$ is unbounded, then there exist $p_0 > 2, K_1 > 0, K_2 \ge 0$ and $\gamma_1 \ge 0$ such that for all $x\in\mathcal D$,
		\begin{flalign*}
			\langle x - x^{\circ}, f(x) \rangle + \frac{p_0 - 1}{2} \sum_{l=1}^m \vert{}g_l(x)\vert{}^2 \le K_1 (1 + \vert{}x - x^{\circ}\vert{}^{2}) - K_2 (1 + \vert{}x - x^{\circ}\vert{}^{2}) \vert{}x - x^{\circ}\vert{}^{\gamma_1}.
		\end{flalign*}
		If $\gamma_1 > 0$, then $K_2 > 0$. If $\mathcal{D}$ is bounded, then this assumption is omitted.
	\end{assumption}
	
	\begin{assumption}\label{ass2}
		For each $j$, there exist $q_j > 0, K_{3,j} > 0, K_{4,j} \ge 0, \gamma_{2,j} \ge 0$ and $\gamma_{3,j} \ge 0$ such that for all $x\in\mathcal D$,
		\begin{flalign*}
			-\phi_j(x) \nu_j^\top f(x) + \frac{q_j + 1}{2} \sum_{l=1}^m |\nu_j^\top g_l(x)|^2 \le K_{3,j} \phi_j(x)^2 (1 + |x - x^{\circ}|^{\gamma_{3,j}}) - K_{4,j} \phi_j(x)^{2 - \gamma_{2,j}}.
		\end{flalign*}
		If $\mathcal{D}$ is unbounded, then $\max_{1\leq j\leq J}\gamma_{3,j} \leq \gamma_1$.
	\end{assumption}
	
	For each $j=1,2,\ldots,J$, define
	\begin{equation*}
		q_j^\dagger :=
		\begin{cases}
			q_j, & \mathcal D \text{ is bounded}, \\
			q_j, & \mathcal D \text{ is unbounded and } \gamma_{3,j} < \gamma_1, \\
			q_j, & \mathcal D \text{ is unbounded and } \gamma_{3,j} = \gamma_1 = 0, \\
			\dfrac{p_0 K_2 q_j}{p_0 K_2 + K_{3,j} q_j}, & \mathcal D \text{ is unbounded and } 0 < \gamma_{3,j} = \gamma_1.
		\end{cases}
	\end{equation*}
	For $\varphi\in C^2(\mathcal{D};\mathbb R)$, define the generator 
	\begin{flalign*}
		&\mathbb L\varphi(x)
		=D\varphi(x)f(x)
		+\frac12\sum_{l=1}^m
		D^2\varphi(x)[g_l(x),g_l(x)]\\
		=&\sum_{i=1}^d f_i(x)
		\frac{\partial\varphi}{\partial x_i}(x)
		+\frac12\sum_{l=1}^m\sum_{i,j=1}^d
		g_{l,i}(x)g_{l,j}(x)
		\frac{\partial^2\varphi}{\partial x_i\partial x_j}(x)
	\end{flalign*}
	for $x\in\mathcal{D}$. Here, $f_i(x)$ and $g_{l,i}(x)$ denote the $i$-th components of $f(x)$ and $g_l(x)$, respectively.
	
	The following lemma establishes global well-posedness of \eqref{ssde} and provides moment and boundary inverse-moment bounds for the exact solution in both bounded and unbounded domains. If $\mathcal D$ is unbounded and $\gamma_1>0$, then it also yields an exponential integrability estimate.
	
	\begin{lemma}\label{Lm4.1}
		Suppose that Assumptions \ref{ass0} and \ref{ass2} hold. When $\mathcal D$ is unbounded, assume in addition that
		Assumption~\ref{ass1} holds. Then the SDE \eqref{ssde} admits a unique strong solution $\{X_t\}_{t\in[0,T]}$ satisfying
		\begin{flalign*}
			\PP(X_t\in \mathcal{D}, \forall t\in[0,T])=1,
		\end{flalign*}
		and the following estimates hold.
		\begin{enumerate}
			\item[(1)] If $\mathcal D$ is bounded, then for every
			$p^*>0$,
			\begin{flalign*}
				\mathbb{E}\Big[ \sup_{0 \le t \le T} \left(1 + |X_t - x^{\circ}|^2\right)^{\frac{p^*}{2}} \Big] \leq C.
			\end{flalign*}
			Furthermore, for each $j=1,2,\ldots,J$ and every
			$0<q_j^*<q_j$,
			\begin{flalign*}
				\mathbb{E}\Big[ \sup_{0 \le t \le T} \phi_j(X_t)^{-q_j^*} \Big] \leq C.
			\end{flalign*}
			\item[(2)] If $\mathcal{D}$ is unbounded, then for every $0<p^* < p_0$,
			\begin{flalign*}
				\mathbb{E}\Big[ \sup_{0 \le t \le T} \left(1 + |X_t - x^{\circ}|^2\right)^{\frac{p^*}{2}} \Big] \leq C.
			\end{flalign*}
			Furthermore, for each $j=1,2,\ldots,J$ and every
			$0<q_j^*<q_j^\dagger$,
			\begin{flalign*}
				\mathbb{E}\Big[ \sup_{0 \le t \le T} \phi_j(X_t)^{-q_j^*} \Big] \leq C.
			\end{flalign*}
			If $\gamma_1>0$, then for every $0<K<p_0K_2$,
			\begin{equation*}
				\mathbb E\Big[\exp\Big( K\int_0^T|X_s-x^{\circ}|^{\gamma_1}\mathrm ds\Big)
				\Big]\leq C.
			\end{equation*}
		\end{enumerate}
		\begin{proof}
			By Assumption \ref{ass0}, the coefficients are locally Lipschitz on $\mathcal D$, and hence there exists a unique maximal local strong solution $\{X_t\}_{0 \le t < \tau_e}$ to \eqref{ssde}, where $\tau_e$ denotes the explosion time. Define 
			\begin{flalign*}
				\mathcal{D}_M=\{x\in\mathcal{D}: |x-x^{\circ}|<M,\quad \phi_j(x)>\frac{1}{M}, j=1,2,\ldots,J\},
			\end{flalign*}
			where $M\geq M_0$ and $M_0$ is sufficiently large such that $X_0\in\mathcal{D}_{M_0}$. Clearly, its closure $\overline{\mathcal{D}}_M$ is compact and satisfies $\overline{\mathcal{D}}_M\subset \mathcal{D}$. Define the stopping time
			$ \tau_M = \inf \{ t \in [0,\tau_e) : X_t\notin\mathcal{D}_M \} \wedge \tau_e$. Since $\overline{\mathcal{D}}_M \Subset \mathcal{D}$, we have $\tau_M < \tau_e$ on $\{\tau_e < \infty\}$, and hence $X_{t \wedge \tau_M}$ is well defined. Since $ \mathcal{D}_M\subset \mathcal{D}_{M+1}$, we have $\tau_M\leq \tau_{M+1}$. Set $\tau_\infty = \lim_{M\rightarrow\infty}\tau_M$. Since $\tau_M\leq \tau_e$, we have $\tau_\infty\leq \tau_e$ a.s. We next show that $\tau_\infty=\tau_e$ a.s. Suppose, to the contrary, that $\PP(\tau_\infty<\tau_e)>0$. On the event $\{\tau_\infty<\tau_e\}$, the continuity of $X_t$ implies that
			$K_\tau:=\{X_t:0\leq t\leq \tau_\infty\} \Subset \mathcal{D}$. Hence
			\begin{flalign*}
				\sup_{0 \le t \le \tau_\infty}|X_t-x^{\circ}|<\infty\quad\text{and}\quad \min_{\substack{0\le t\le \tau_\infty \\ 1\le j\le J}} \phi_j(X_t) > 0.
			\end{flalign*}
			Moreover, there exists a sufficiently large $M_* \ge M_0$ such that $K_\tau \Subset \mathcal{D}_{M_*}$. Since $\tau_\infty < \tau_e$ and $X_t$ is continuous, there exists $\varepsilon > 0$ such that $\tau_\infty + \varepsilon < \tau_e$ and $X_t \in \mathcal{D}_{M_*}$, $0 \le t \le \tau_\infty + \varepsilon$.
			By the definition of $\tau_{M_*}$, this implies $\tau_{M_*} \ge \tau_\infty + \varepsilon > \tau_\infty$, which contradicts $\tau_{M_*} \le \tau_\infty = \lim_{M \to \infty} \tau_M$. Hence, $\tau_M \uparrow \tau_e $ a.s.
			
			(A) We first establish moment bounds for the stopped solution. 
			
			If $\mathcal{D}$ is bounded, then $C_{\mathcal{D}} := \sup_{x \in \mathcal{D}} |x - x^{\circ}| < \infty$, and hence for every $p^* > 0$,
			\begin{flalign*}
				\sup_{M \ge M_0} \mathbb{E} \Big[ \sup_{0 \le t \le T} (1 + |X_{t \wedge \tau_M} - x^{\circ}|^2)^{\frac{p^*}{2}} \Big] \le (1 + C_{\mathcal{D}}^2)^{\frac{p^*}{2}}.
			\end{flalign*}
			Therefore, it remains to consider the case where $\mathcal{D}$ is unbounded.
			Define $ Q(x) = (1 + |x - x^{\circ}|^2)^{\frac{p_0}{2}} $. Then by Assumption \ref{ass1}, we get
			\begin{flalign}\label{4.1}
				\mathbb L Q(x) & \leq p_0 (1 + |x - x^{\circ}|^2)^{\frac{p_0}{2} - 1} \Big( \langle x - x^{\circ}, f(x) \rangle + \frac{p_0 - 1}{2} \sum_{i=1}^m |g_i(x)|^2 \Big) \nonumber\\
				& \le p_0 (1 + |x - x^{\circ}|^2)^{\frac{p_0}{2}-1} \left( K_1 (1 + |x - x^{\circ}|^2) - K_2 (1 + |x - x^{\circ}|^2) |x - x^{\circ}|^{\gamma_1} \right) \nonumber\\
				& = K_1 p_0 Q(x) - K_2 p_0 Q(x) |x - x^{\circ}|^{\gamma_1}\,, 
			\end{flalign} 
			which implies  $\mathbb L Q(x)\leq K_1 p_0 Q(x)$. Define $\widehat{Z}_{Q,M}(t) := e^{-K_1 p_0 (t \wedge \tau_M)} Q(X_{t \wedge \tau_M}).$ By It\^o's formula,
			\begin{flalign*}
				\mathrm{d}\widehat{Z}_{Q,M}(t) = e^{-K_1 p_0 t} \mathbb I_{\{t < \tau_M\}} \left( \mathbb LQ(X_t) - K_1 p_0 Q(X_t) \right) \mathrm{d}t + e^{-K_1 p_0 t} \mathbb I_{\{t < \tau_M\}} \sum_{l=1}^m DQ(X_t) g_l(X_t) \, \mathrm{d}B_t^l.
			\end{flalign*}
			Since $\overline{\mathcal{D}}_M \Subset \mathcal{D}$, the stopped stochastic integrand is bounded. Hence the stochastic integral is a square-integrable martingale, and therefore $\widehat{Z}_{Q,M}$ is a nonnegative supermartingale. Taking $0 < p^* < p_0$, $\theta^* = \frac{p^*}{p_0} \in (0, 1)$ and using $Q(X_{t \wedge \tau_M}) = e^{K_1 p_0 (t \wedge \tau_M)} \widehat{Z}_{Q,M}(t) \le e^{K_1 p_0 T} \widehat{Z}_{Q,M}(t)$, as well as Proposition \ref{ville}, we have
			\begin{flalign*}
				\mathbb{E} \Big[ \Big( \sup_{0 \le t \le T} \widehat{Z}_{Q,M}(t) \Big)^{\theta^*} \Big] \le \frac{Q(X_0)^{\theta^*}}{1 - \theta^*}
			\end{flalign*}
			and
			\begin{flalign*}
				\E\Big[\sup_{0 \le t \le T} (1 + |X_{t \wedge \tau_M} - x^{\circ}|^2)^{\frac{p^*}{2}} \Big]= \E\Big[\sup_{0 \le t \le T} Q(X_{t \wedge \tau_M})^{\theta^*}\Big]\leq \frac{e^{K_1 p_0 T \theta^*}Q(X_0)^{\theta^*}}{1-\theta^*}.
			\end{flalign*}
			Hence,
			\begin{flalign}\label{4.3}
				\sup_M \E\Big[\sup_{0 \le t \le T} (1 + |X_{t \wedge \tau_M} - x^{\circ}|^2)^{\frac{p^*}{2}} \Big]\leq C. 
			\end{flalign}
			\indent If $\mathcal{D}$ is unbounded and $\gamma_1>0$, then $K_2>0$. Fix any constant $0<K < K_2 p_0$ and define $ I_{Q, M}(t) = \int_0^{t \wedge \tau_M} |X_s - x^{\circ}|^{\gamma_1} \mathrm{d}s $, as well as
			$ E(t) = \exp \left( - K_1 p_0 (t \wedge \tau_M) + K I_{Q, M}(t) \right) $
			and $ Z_{Q, M}(t) = E(t) Q(X_{t \wedge \tau_M}) $. Applying the chain rule to $E(t)$ yields
			\begin{flalign*}
				\mathrm{d}E(t) = E(t) \mathbb{I}_{\{t < \tau_M\}} \left( -K_1 p_0 + K |X_t - x^{\circ}|^{\gamma_1} \right) \mathrm{d}t
			\end{flalign*}
			and
			\begin{flalign*}
				\mathrm{d}Q(X_{t \wedge \tau_M}) = \mathbb{I}_{\{t < \tau_M\}} \mathbb{L} Q(X_t) \mathrm{d}t + \mathbb{I}_{\{t < \tau_M\}} \sum_{l=1}^m D Q(X_t) g_l(X_t) \mathrm{d}B_t^l.
			\end{flalign*}
			Combining these identities with the inequality \eqref{4.1}, we obtain 
			\begin{flalign*}
				\mathrm{d} Z_{Q, M}(t) & = E(t) \mathbb{I}_{\{t < \tau_M\}} \left( \mathbb{L} Q(X_t) - K_1 p_0 Q(X_t) + K Q(X_t) |X_t - x^{\circ}|^{\gamma_1} \right) \mathrm{d}t \\
				& \quad + E(t) \mathbb{I}_{\{t < \tau_M\}} \sum_{l=1}^m D Q(X_t) g_l(X_t) \mathrm{d}B_t^l \\
				& \le E(t) \mathbb{I}_{\{t < \tau_M\}} (K_1 p_0 Q(X_t) - K_2 p_0 Q(X_t) |X_t - x^{\circ}|^{\gamma_1} - K_1 p_0 Q(X_t)  + K Q(X_t) |X_t - x^{\circ}|^{\gamma_1}) \mathrm{d}t\\
				& \quad + E(t) \mathbb{I}_{\{t < \tau_M\}} \sum_{l=1}^m D Q(X_t) g_l(X_t) \mathrm{d}B_t^l \\
				& = -(K_2 p_0 - K) E(t) \mathbb{I}_{\{t < \tau_M\}} Q(X_t) |X_t - x^{\circ}|^{\gamma_1} \mathrm{d}t + E(t) \mathbb{I}_{\{t < \tau_M\}} \sum_{l=1}^m D Q(X_t) g_l(X_t) \mathrm{d}B_t^l.
			\end{flalign*} 
			Integrating both sides of the above  inequality  from $0$ to $t$ yields
			\begin{align*}
				Z_{Q, M}(t) & \le Q(X_0) - (K_2 p_0 - K) \int_0^{t \wedge \tau_M} e^{-K_1 p_0 s + K I_{Q, M}(s)} |X_s - x^{\circ}|^{\gamma_1} Q(X_s) \mathrm{d}s + N_{Q, M}(t),
			\end{align*} 
			where
			\begin{align*}
				N_{Q, M}(t) = \sum_{l=1}^m \int_0^{t \wedge \tau_M} e^{-K_1 p_0 s + K I_{Q, M}(s)} D Q(X_s) g_l(X_s) \mathrm{d}B_s^l.
			\end{align*} 
			Since the stopped stochastic integrand is bounded, $N_{Q,M}$ is a square-integrable martingale, and hence $Z_{Q,M}$ is a nonnegative supermartingale. Therefore, $\E [Z_{Q, M}(T)]\le Z_{Q,M}(0)=Q(X_0)$.
			On the other hand, $Z_{Q,M}(T)\geq e^{-K_1p_0T}e^{KI_{Q,M}(T)}$. Therefore,
			\begin{flalign}\label{4.4}
				\sup_{M\geq M_0}\E\Big[\exp\Big(K\int_{0}^{T\wedge \tau_M}|X_s-x^{\circ}|^{\gamma_{1}}\mathrm{d}s\Big)\Big]\leq e^{p_0K_1T}Q(X_0)=:C.
			\end{flalign}
			\indent (B) We next establish the inverse-moment estimates of the stopped solution for both bounded and unbounded domains. 
			
			Define $ V_j(x) = \phi_j(x)^{-q_j} $. Then by Assumption \ref{ass2},
			\begin{flalign}\label{4.2}
				\mathbb{L} V_j(x) 
				& = q_j \phi_j(x)^{-q_j - 2} \Big( -\phi_j(x) \nu_j^T f(x) + \frac{q_j + 1}{2} \sum_{i=1}^m |\nu_j^T g_i(x)|^2 \Big) \nonumber \\
				& \le q_j \phi_j(x)^{-q_j - 2} \Big( K_{3,j} \phi_j(x)^2 (1 + |x - x^{\circ}|^{\gamma_{3,j}}) - K_{4,j}\phi_j(x)^{2 - \gamma_{2,j}} \Big) \nonumber\\
				& = K_{3,j} q_j V_j(x) (1 + |x - x^{\circ}|^{\gamma_{3,j}}) - K_{4,j} q_j \phi_j(x)^{-\gamma_{2,j}} V_j(x).
			\end{flalign} 
			Besides, we define
			$ A_{j, M}(t) = K_{3,j} q_j \int_0^{t \wedge \tau_M} (1 + |X_s - x^{\circ}|^{\gamma_{3,j}}) \mathrm{d}s $
			and
			$ Z_{j, M}(t) = e^{-A_{j, M}(t)} V_j(X_{t \wedge \tau_M}) $.
			By calculation and the inequality \eqref{4.2}, we obtain
			\begin{flalign*}
				\mathrm{d}A_{j, M}(t) &= K_{3,j} q_j \left( 1 + |X_t - x^{\circ}|^{\gamma_{3,j}} \right) \mathbb I_{\{t < \tau_M\}} \mathrm{d}t
			\end{flalign*}
			and
			\begin{flalign*}
				\mathrm{d}Z_{j, M}(t) = &e^{-A_{j, M}(t)} \mathbb I_{\{t < \tau_M\}} ( \mathbb{L} V_j(X_t) - K_{3,j} q_j \left( 1 + |X_t - x^{\circ}|^{\gamma_{3,j}} \right) V_j(X_t) ) \mathrm{d}t \\
				&+ e^{-A_{j, M}(t)} \mathbb I_{\{t < \tau_M\}} \sum_{l=1}^m D V_j(X_t) g_l(X_t) \mathrm{d}B_t^l\\
				\leq&-K_{4,j}q_je^{-A_{j, M}(t)} \mathbb I_{\{t < \tau_M\}}\phi_j(X_t)^{-\gamma_{2,j}} V_j(X_t)\mathrm{d}t+e^{-A_{j, M}(t)} \mathbb I_{\{t < \tau_M\}} \sum_{l=1}^m D V_j(X_t) g_l(X_t) \mathrm{d}B_t^l.
			\end{flalign*}
			Integrating this from $0$ to $t$ yields
			\begin{flalign*}
				Z_{j, M}(t) \le V_j(X_0) - K_{4,j} q_j \int_0^{t \wedge \tau_M} e^{-A_{j, M}(s)} \phi_j(X_s)^{-\gamma_{2,j}} V_j(X_s) \mathrm{d}s + N_{j, M}(t),
			\end{flalign*}
			where
			\begin{align*}
				N_{j, M}(t) = \sum_{l=1}^m \int_0^{t \wedge \tau_M} e^{-A_{j, M}(s)} D V_j(X_s) g_l(X_s) \mathrm{d}B_s^l.
			\end{align*} 
			Similarly, $N_{j,M}$ is a square-integrable martingale, and hence $Z_{j,M}$ is a nonnegative supermartingale. Therefore, by Proposition \ref{ville}, for every $0<u<1$ and every $j$,
			\begin{flalign*}
				\E \Big[ \Big(\sup_{0 \le t \le T}  Z_{j, M}(t)\Big)^u \Big]\leq \frac{V_j(X_0)^u}{1-u} < \infty.
			\end{flalign*}
			Fix $j\in\{1,2,\ldots,J\}$ and $0<q_j^*<q_j^\dagger$, $\theta_j = \frac{q_j^*}{q_j} \in (0, 1)$. Since $V_j(X_{t \wedge \tau_M}) = e^{A_{j, M}(t)} Z_{j, M}(t)$, we have
			\begin{flalign}\label{sup0T}
				\sup_{0 \le t \le T} \phi_j(X_{t \wedge \tau_M})^{-q_j^*} \le e^{\theta_j A_{j, M}(T)} \big(\sup_{0 \le t \le T}Z_{j, M}(t)\big)^{\theta_j}.
			\end{flalign}
			For any $u\in(\theta_j,1)$, define $\Lambda_{j,u}:=\frac{K_{3,j}q_jq_j^*u}{q_ju-q_j^*}$. By H\"older's inequality and \eqref{sup0T},
			\begin{flalign*}
				\mathbb E \Big[ \sup_{0\leq t\leq T} \phi_j(X_{t\wedge\tau_M})^{-q_j^*} \Big] 
				\leq \Big( \mathbb E  \Big( \sup_{0\leq t\leq T} Z_{j,M}(t) \Big)^u  \Big)^{\frac{\theta_j}{u}} 
				\Big( \mathbb E \exp \Big( \Lambda_{j,u} \int_0^{T\wedge\tau_M} \bigl( 1+|X_s-x^{\circ}|^{\gamma_{3,j}} \bigr) \mathrm ds \Big) \Big)^{\frac{u-\theta_j}{u}}.
			\end{flalign*}
			
			If $\mathcal{D}$ is bounded, then $A_{j,M}(T)\leq K_{3,j}q_jT(1+C_\mathcal{D}^{\gamma_{3,j}})=:C_{j,T}$. Substituting this into \eqref{sup0T} yields 
			\begin{flalign*}
				\E\Big[\sup_{0 \le t \le T} \phi_j(X_{t \wedge \tau_M})^{-q_j^*} \Big]\le e^{\theta_j C_{j,T}} \E\Big[\Big(\sup_{0 \le t \le T}Z_{j, M}(t)\Big)^{\theta_j}\Big]\leq e^{\theta_j C_{j,T}} \frac{V_j(X_0)^{\theta_j}}{1-\theta_j}.
			\end{flalign*}
			Therefore, if $\mathcal{D}$ is bounded, then it holds that
			\begin{flalign}\label{bounded}
				\sup_{M \ge M_0}\E\Big[\sup_{0 \le t \le T} \phi_j(X_{t \wedge \tau_M})^{-q_j^*} \Big]\leq C.
			\end{flalign}
			
			On the other hand, if $\mathcal{D}$ is unbounded and $\gamma_{1}=0$, then Assumption \ref{ass2} implies \(\gamma_{3,j}=0\). Hence, one has $A_{j,M}(T)\leq 2K_{3,j}q_jT$. Further, it holds that for $0<q_j^*<q_j^\dagger$,
			\begin{flalign}\label{unbounded1}
				\sup_{M \ge M_0}\E\Big[\sup_{0 \le t \le T} \phi_j(X_{t \wedge \tau_M})^{-q_j^*} \Big]\leq C.
			\end{flalign}
			
			If $\mathcal D$ is unbounded and $ 0\leq\gamma_{3,j}<\gamma_1$, choose any $u_j\in(\theta_j,1)$. Then $\Lambda_{j,u_j}\leq C$. For any $\varepsilon > 0$, there exists $C_{\varepsilon,j}>0$ such that
			$ |X_s - x^{\circ}|^{\gamma_{3,j}} \le \varepsilon |X_s - x^{\circ}|^{\gamma_1} + C_{\varepsilon,j}, $
			which yields
			\begin{flalign*}
				\Lambda_{j, u_j} \int_0^{T \wedge \tau_M} (1 + |X_s - x^{\circ}|^{\gamma_{3,j}}) \mathrm{d}s \le C_{\varepsilon, j}T + \Lambda_{j, u_j} \varepsilon \int_0^{T \wedge \tau_M}  |X_s - x^{\circ}|^{\gamma_1}\mathrm{d}s.
			\end{flalign*}
			Choose \(\varepsilon>0\) and \(K\) such that $ \Lambda_{j, u_j} \varepsilon<K < p_0K_2. $
			By the inequality \eqref{4.4}, we show that for $0<q_j^*<q_j^\dagger$,
			\begin{flalign}\label{unbounded2}
				\sup_{M\geq M_0}\E \Big[ \sup_{0 \le t \le T} \phi_j(X_{t \wedge \tau_M})^{-q_j^*} \Big] \le C.
			\end{flalign}
			
			If $\mathcal{D}$ is unbounded and $\gamma_{1}=\gamma_{3,j}>0$, then we let $0<q_j^*<q_j^\dagger=\frac{p_0K_2q_j}{p_0K_2+K_{3,j}q_j}$. Since $\lim_{u \uparrow 1} \Lambda_{j,u} = \frac{K_{3,j} q_j q_j^*}{q_j - q_j^*}<p_0K_2$, there exists $u_j\in(\theta_j,1)$ sufficiently close to $1$ such that $ \Lambda_{j, u_j}<p_0K_2$. By using H\"older's inequality, we derive
			\begin{flalign*}
				\E \Big[ \sup_{0 \le t \le T} \phi_j(X_{t \wedge \tau_M})^{-q_j^*} \Big] 
				& \le \Big( 
				\E\Big[\Big(\sup_{0 \le t \le T}Z_{j, M}(t)\Big)^{u_j}\Big] \Big)^{\frac{\theta_j}{u_j} } \Big( \E \Big[ \exp \Big( \frac{\theta_j u_j}{u_j - \theta_j} A_{j, M}(T) \Big) \Big] \Big)^{\frac{u_j - \theta_j}{u_j}} \\
				& \le C \Big( \E \Big[ \exp \Big( \Lambda_{j, u_j} \int_0^{T \wedge \tau_M} (1 + |X_s - x^{\circ}|^{\gamma_{3,j}}) \mathrm{d}s \Big) \Big] \Big)^{\frac{u_j - \theta_j}{u_j}}.
			\end{flalign*}
			Choose $K$ such that $ \Lambda_{j, u_j}< K <p_0K_2$. Then
			\begin{flalign*}
				\Lambda_{j,u_j} \int_0^{T \wedge \tau_M} \left(1 + |X_s - x^{\circ}|^{\gamma_1}\right) \mathrm{d}s   \leq \Lambda_{j,u_j} T + K \int_0^{T \wedge \tau_M} |X_s - x^{\circ}|^{\gamma_1} \mathrm{d}s .
			\end{flalign*}
			Hence, by using the inequality \eqref{4.4}, we derive
			\begin{flalign*}
				\sup_{M\geq M_0}\E \Big[ \sup_{0 \le t \le T} \phi_j(X_{t \wedge \tau_M})^{-q_j^*} \Big] \le C.
			\end{flalign*}
			Combining the above cases, for every $j=1,2,\ldots,J$ and every $0<q_j^*<q_j^\dagger$,
			\begin{flalign*}
				\sup_{M\geq M_0}\E \Big[ \sup_{0 \le t \le T} \phi_j(X_{t \wedge \tau_M})^{-q_j^*} \Big] \le C.
			\end{flalign*}
			
			(C) Now we remove the stopping times to get the conclusions.
			
			If $\mathcal{D}$ is unbounded, then by Markov's inequality, the inequalities \eqref{4.3}, and the bound above, we obtain
			\begin{equation*}
				\begin{aligned}
					\PP(\tau_M \le T) &\le \frac{\mathbb{E}\Big[ \sup_{0 \le t \le T} \left(1 + |X_{t \wedge \tau_M} - x^{\circ}|^2\right)^{p^*/2} \Big]}{(1 + M^2)^{p^*/2}} + \sum_{j=1}^J \frac{\mathbb{E}\left[ \sup_{0 \le t \le T} \phi_j(X_{t \wedge \tau_M})^{-q_j^*} \right]}{M^{q_j^*}} \\
					&\le C M^{-p^*} + C \sum_{j=1}^J M^{-q_j^*}.
				\end{aligned}
			\end{equation*}
			
			If $\mathcal{D}$ is bounded, then taking $M>C_\mathcal{D}$ gives
			\begin{equation*}
				\PP(\tau_M \le T) \le  C \sum_{j=1}^J M^{-q_j^*}.
			\end{equation*}
			In either case, $\lim_{M \to \infty} \mathbb{P}(\tau_M \le T) = 0$. Since $\tau_M \uparrow \tau_e$ a.s., it follows that 
			\begin{flalign*}
				\mathbb{P}(\tau_e \le T) = \lim_{M \to \infty} \mathbb{P}(\tau_M \le T) = 0.
			\end{flalign*}
			As $T > 0$ is arbitrary, we conclude that $\tau_e = \infty$ a.s.
			Since $\tau_M \uparrow \tau_e = \infty$ a.s.,
			for every fixed $T > 0$, $X_{t \wedge \tau_M} = X_t$, $0 \le t \le T$ for all sufficiently large $M$, a.s. Hence, by Fatou's lemma, we derive
			\begin{flalign*}
				\mathbb{E} \Big[ \sup_{0 \le t \le T} (1 + |X_t - x^{\circ}|^2)^{\frac{p^*}{2}} \Big] \le C,
			\end{flalign*}
			where $0 < p^* < p_0$ if $\mathcal{D}$ is unbounded, while $p^* > 0$ if $\mathcal{D}$ is bounded. Similarly, for every $j = 1,2, \dots, J$ and $0 < q_j^* < q_j^\dagger$,
			\begin{flalign*}
				\mathbb{E} \Big[ \sup_{0 \le t \le T} \phi_j(X_t)^{-q_j^*} \Big] \le C.
			\end{flalign*}
			Moreover, if $\mathcal D$ is unbounded and $\gamma_1 > 0$, then for every $0 < K < p_0 K_2$,
			\begin{flalign*}
				\mathbb E \Big[ \exp \Big( K \int_0^T |X_s - x^{\circ}|^{\gamma_1} \mathrm{d}s \Big) \Big] 
				\leq \liminf_{M \to \infty} \mathbb E \Big[ \exp \Big( K \int_0^{T \wedge \tau_M} |X_s - x^{\circ}|^{\gamma_1} \mathrm{d}s \Big) \Big] \leq C,
			\end{flalign*}
			where Fatou's lemma and \eqref{4.4} are used. Therefore, all the assertions of this lemma hold. This completes the proof. \hfill
		\end{proof}
	\end{lemma}
	
	Based on the preceding results, we establish exponential
	integrability of time integrals of inverse powers of
	$\phi_j(X_t)$.
	\begin{lemma}\label{Lm4.2}
		Suppose the assumptions in Lemma \ref{Lm4.1} hold. Fix $j\in\{1,2,\ldots,J\}$ satisfying $\gamma_{2,j}>0$ and $K_{4,j} > 0$. Then, for every $0 < \bar{\mathscr{X}}_j < \mathscr{X}_j^\dagger$, there exists $C > 0$ such that 
		\begin{flalign*}
			\mathbb{E} \Big[ \exp \Big( \bar{\mathscr{X}}_j \int_0^T \phi_j(X_s)^{-\gamma_{2,j}} \, \mathrm{d}s \Big) \Big] \le C,
		\end{flalign*}
		where
		\begin{flalign*}
			\mathscr{X}_j^\dagger = \begin{cases}
				\dfrac{p_0 K_2 b_j}{a_j + p_0 K_2}, & \mathcal{D} \text{ is unbounded and } 0 < \gamma_{3,j} = \gamma_1,\\
				b_j, & otherwise, 
			\end{cases}
		\end{flalign*}
		with  $
		b_j = \frac{K_{4,j} q_j}{2}$ and $ a_j = K_{3,j} q_j + b_j$.
		\begin{proof}
			Define $W_j(x) = 1 + V_j(x) = 1 + \phi_j(x)^{-q_j}$ and set $b_j = \frac{K_{4,j} q_j}{2}$, $a_j = K_{3,j} q_j + b_j$. Clearly, it holds that
			$ \phi_j(x)^{-\gamma_{2,j}} \le 1 + \phi_j(x)^{-q_j - \gamma_{2,j}}$.
			Then we have
			\begin{equation*}
				\begin{aligned}
					-2b_j \phi_j(x)^{-q_j - \gamma_{2,j}} \le b_j - b_j \phi_j(x)^{-\gamma_{2,j}} - b_j \phi_j(x)^{-q_j - \gamma_{2,j}} = b_j - b_j \phi_j(x)^{-\gamma_{2,j}}W_j(x).
				\end{aligned}
			\end{equation*}	
			Substituting this into $\mathbb{L} V_j$, we obtain
			\begin{flalign}\label{LWj}
				\mathbb{L} W_j(x) =\mathbb{L} V_j(x)\leq &K_{3,j} q_j V_j(x) (1 + |x - x^{\circ}|^{\gamma_{3,j}})+ b_j - b_j \phi_j(x)^{-\gamma_{2,j}} W_j(x) \nonumber\\
				\le &a_j (1 + |x - x^{\circ}|^{\gamma_{3,j}}) W_j(x) - b_j \phi_j(x)^{-\gamma_{2,j}} W_j(x).
			\end{flalign}
			Define $J_{j, M}(t) = \int_0^{t \wedge \tau_M} (1 + |X_s - x^{\circ}|^{\gamma_{3,j}}) \mathrm{d}s$ and $I_{j, M}(t) = \int_0^{t \wedge \tau_M} \phi_j(X_s)^{-\gamma_{2,j}} \mathrm{d}s$. For any $0 < \mathscr{X}_j < b_j$, define $ U_{j, M}(t) = \exp ( -a_j J_{j, M}(t) + \mathscr{X}_j I_{j, M}(t) ) W_j(X_{t \wedge \tau_M})$.
			By \eqref{LWj}, we obtain
			\begin{flalign*}
				\mathrm{d}U_{j, M}(t)=&e^{-a_j J_{j, M}(t) + \mathscr{X}_j I_{j, M}(t)}\mathbb{I}_{\{t<\tau_M\}}\Big[\mathbb{L} W_j(X_t)-a_j(1 + |X_t - x^{\circ}|^{\gamma_{3,j}})W_j(X_t)+\mathscr{X}_j\phi_j(X_t)^{-\gamma_{2,j}}W_j(X_t)\Big]\mathrm{d}t\\
				&+e^{-a_j J_{j, M}(t) + \mathscr{X}_j I_{j, M}(t)}\mathbb{I}_{\{t<\tau_M\}}\sum_{l=1}^{m}DW_j(X_t)g_l(X_t)\mathrm{d}B_t^l\\
				\leq &- (b_j - \mathscr{X}_j)e^{-a_j J_{j, M}(t) + \mathscr{X}_j I_{j, M}(t)}\mathbb{I}_{\{t<\tau_M\}}\phi_j(X_t)^{-\gamma_{2,j}}W_j(X_t)\mathrm{d}t\\
				&+e^{-a_j J_{j, M}(t) + \mathscr{X}_j I_{j, M}(t)}\mathbb{I}_{\{t<\tau_M\}}\sum_{l=1}^{m}DW_j(X_t)g_l(X_t)\mathrm{d}B_t^l.
			\end{flalign*}
			Integrating this from $0$ to $t$ gives
			\begin{flalign*}
				U_{j, M}(t) &\le W_j(X_0) - (b_j - \mathscr{X}_j) \int_0^{t \wedge \tau_M} e^{-a_j J_{j, M}(s) + \mathscr{X}_j I_{j, M}(s)} \cdot \phi_j(X_s)^{-\gamma_{2,j}} W_j(X_s) \mathrm{d}s + \tilde{N}_{j, M}(t),
			\end{flalign*}
			where
			\begin{flalign*}
				\tilde{N}_{j, M}(t)=\sum_{l=1}^{m}\int_{0}^{t\wedge\tau_M}e^{-a_j J_{j, M}(s) + \mathscr{X}_j I_{j, M}(s)}DW_j(X_s)g_l(X_s)\mathrm{d}B_s^l.
			\end{flalign*}
			Since the stopped stochastic integrand is bounded, the stochastic integral is a martingale. Hence $U_{j,M}$ is a nonnegative supermartingale. Hence, $\mathbb{E}[U_{j, M}(t)] \le W_j(X_0)$. Since $
			U_{j, M}(T) \ge \exp ( -a_j J_{j, M}(T) + \mathscr{X}_j I_{j, M}(T) )$, we get 
			\begin{flalign*}
				\mathbb{E} [ \exp ( \mathscr{X}_j I_{j, M}(T) - a_j J_{j, M}(T)) ] \le W_j(X_0).
			\end{flalign*}
			Let $0 < \bar{\mathscr{X}}_j < \mathscr{X}_j^\dagger$. We choose $\mathscr{X}_j \in (\bar{\mathscr{X}}_j, b_j)$ appropriately below and set $\bar{\theta}_j := \frac{\bar{\mathscr{X}}_j}{\mathscr{X}_j} \in (0, 1)$. Then $\bar{\theta}_j \mathscr{X}_j = \bar{\mathscr{X}}_j$, and hence
			\begin{flalign*}
				e^{\bar{ \mathscr{X}}_j I_{j, M}(T)} = ( e^{\mathscr{X}_j I_{j, M}(T) - a_j J_{j, M}(T)} )^{\bar{\theta}_j} e^{\bar{\theta}_j a_j J_{j, M}(T)}.
			\end{flalign*}
			Then we further have
			\begin{equation*}
				\mathbb{E} \left[ e^{\bar{ \mathscr{X}}_j I_{j, M}(T)} \right] \le \Big( \mathbb{E} \Big[ e^{\mathscr{X}_j I_{j, M}(T) - a_j J_{j, M}(T)} \Big] \Big)^{\bar{\theta}_j} \Big( \mathbb{E} \Big[ e^{\frac{\bar{\theta}_j a_j}{1 - \bar{\theta}_j} J_{j, M}(T)} \Big] \Big)^{1 - \bar{\theta}_j} \le W_j(X_0)^{\bar{\theta}_j} \Big( \mathbb{E} \Big[ e^{\frac{\bar{\theta}_j a_j}{1 - \bar{\theta}_j} J_{j, M}(T)} \Big] \Big)^{1 - \bar{\theta}_j}.
			\end{equation*}
			
			If $\mathcal{D}$ is bounded, then $J_{j,M}(T)\leq T(1+C_{\mathcal{D}}^{\gamma_{3,j}})$. Thus, by the inequality above, we obtain $\mathbb{E} [  e^{\bar{ \mathscr{X}}_j I_{j, M}(T)}] \le C$. If $\mathcal{D}$ is unbounded and $\gamma_{3,j}=\gamma_{1}=0$, then $J_{j,M}(T)\leq 2T$. Hence, it follows that $\mathbb{E} [e^{\bar{ \mathscr{X}}_j I_{j, M}(T)} ] \le C$. If $\mathcal{D}$ is unbounded and $0 \leq \gamma_{3,j} < \gamma_1$, then by the inequality $ |X_s - x^{\circ}|^{\gamma_{3,j}} \le \varepsilon |X_s - x^{\circ}|^{\gamma_1} + C_{\varepsilon,j}$, we have 
			\begin{flalign*}
				\frac{\bar{\theta}_j a_j}{1 - \bar{\theta}_j} J_{j, M}(T)\leq \frac{\bar{\theta}_j a_j}{1 - \bar{\theta}_j}(1+C_{\varepsilon,j})T+\frac{\bar{\theta}_j a_j}{1 - \bar{\theta}_j}\varepsilon \int_0^{T \wedge \tau_M}|X_s - x^{\circ}|^{\gamma_{1}} \mathrm{d}s.
			\end{flalign*}
			Taking a sufficiently small $ \varepsilon $ such that $\frac{\bar{\theta}_j a_j}{1 - \bar{\theta}_j}\varepsilon<K<p_0K_2$ and using \eqref{4.4}, we have $\mathbb{E} [ e^{\bar{ \mathscr{X}}_j I_{j, M}(T)}] \le C$. If $\mathcal{D}$ is unbounded and $0 < \gamma_{3,j} = \gamma_1$, since $ \bar{\mathscr{X}}_j < \mathscr{X}_j^\dagger = \frac{p_0 K_2 b_j}{a_j + p_0 K_2}$, we may choose $\mathscr{X}_j \in (\frac{a_j + p_0 K_2}{p_0 K_2}\bar{\mathscr{X}}_j, b_j)$. Then $\frac{\bar{\theta}_j a_j}{1 - \bar{\theta}_j} < p_0 K_2$. Hence, we may choose $K$ such that $\frac{\bar{\theta}_j a_j}{1 - \bar{\theta}_j} < K < p_0 K_2$. Since $\gamma_{3,j} = \gamma_1$, we have
			\begin{flalign*}
				\frac{\bar{\theta}_j a_j}{1 - \bar{\theta}_j} J_{j, M}(T)\leq \frac{\bar{\theta}_j a_j}{1 - \bar{\theta}_j}T+K\int_0^{T \wedge \tau_M}|X_s - x^{\circ}|^{\gamma_{1}} \mathrm{d}s.
			\end{flalign*}
			Hence, by \eqref{4.4}, we obtain $\mathbb{E} \left[ e^{\bar{ \mathscr{X}}_j I_{j, M}(T)} \right] \le C$. Therefore,
			\begin{flalign*}
				\sup_{M \ge M_0}\E\Big[\exp\Big(\bar{\mathscr{X}}_j \int_0^{T \wedge \tau_M} \phi_j(X_s)^{-\gamma_{2,j}} \mathrm{d}s\Big)\Big] \le C, \quad 0< \bar{\mathscr{X}}_j < \mathscr{X}_j^\dagger.
			\end{flalign*}
			Since $\tau_M \uparrow \infty$ a.s., by Lemma \ref{Lm4.1}, Fatou's lemma yields
			\begin{flalign*}
				\mathbb{E} \Big[\exp \Big( \bar{\mathscr{X}}_j \int_0^T \phi_j(X_s)^{-\gamma_{2,j}} \mathrm{d}s \Big) \Big]
				&\le \liminf_{M \to \infty} \mathbb{E}\Big[ \exp \Big( \bar{\mathscr{X}}_j \int_0^{T \wedge \tau_M} \phi_j(X_s)^{-\gamma_{2,j}} \mathrm{d}s \Big) \Big]\le C,
			\end{flalign*}
			which validates the assertion. \hfill
		\end{proof}
	\end{lemma}
	
	To account for coefficient growth and boundary singularities
	in the one-step stability analysis, we introduce the following
	weight function
	\begin{flalign*}
		\rho(x) = 1 + d_0 |x-x^{\circ}|^{\kappa_0} + \sum_{j=1}^J d_j \phi_j(x)^{-\kappa_j},
	\end{flalign*}
	where $d_0,d_j\geq 0$, $\kappa_0,\kappa_j\geq 0$. Define $\mathcal{J}_\rho := \{ j \in \{1,2,\dots, J\} : d_j > 0, \, \kappa_j > 0 \}$, $\mathcal{J}_{\rho,c} := \{ j \in\mathcal{J}_\rho : \kappa_j=\gamma_{2,j} > 0 \}$,  $\delta_{0,c}:=\mathbb{I}_{\{\mathcal{D} \, \text{is unbounded}, \, d_0>0, \, \kappa_0=\gamma_{1}>0\}}$, and \begin{flalign*}
		\frac{1}{\mathscr{X}_\rho^\dagger}
		:=
		\begin{cases}
			\dfrac{d_0}{p_0K_2}
			+\displaystyle\sum_{j\in\mathcal J_{\rho,c}}
			\dfrac{d_j}{\mathscr X_j^\dagger},
			& \delta_{0,c}=1, \\[2ex]
			\displaystyle\sum_{j\in\mathcal J_{\rho,c}}
			\dfrac{d_j}{\mathscr X_j^\dagger},
			& \delta_{0,c}=0,
		\end{cases}
	\end{flalign*}
	with the convention that $\mathscr X_\rho^\dagger=\infty$ if the right-hand side is zero.
	
	\begin{lemma}\label{Lm4.3}
		Suppose the assumptions in Lemma \ref{Lm4.1} hold. Assume further that
		$K_{4,j}>0$, $0<\kappa_j\le\gamma_{2,j}$, $j\in\mathcal J_\rho$.
		If $\mathcal D$ is unbounded and $d_0\kappa_0>0$, then assume additionally that $0<\kappa_0\le\gamma_1$. Then for every $0<\mathscr{X}<\mathscr{X}_\rho^\dagger$, there exists $C>0$ such that
		\begin{flalign*}
			\mathbb E\Big[\exp\Big(\mathscr{X}\int_0^T\rho(X_s) \mathrm{d}s\Big)\Big]\le C.
		\end{flalign*}
		\begin{proof}
			Fix $0<\mathscr X<\mathscr X_\rho^\dagger$. If $\mathcal D$ is unbounded and $0<\kappa_0<\gamma_1$, then, for every $\varepsilon_0>0$, there exists $C_{\varepsilon_0}>0$ such that $|x-x^{\circ}|^{\kappa_0} \le \varepsilon_0|x-x^{\circ}|^{\gamma_1}
			+ C_{\varepsilon_0}$. If $\mathcal D$ is bounded, then $|x-x^{\circ}|$ is uniformly bounded. Similarly, for every $j\in\mathcal J_\rho$ with $0<\kappa_j<\gamma_{2,j}$ and every $\varepsilon_j>0$, there exists $C_{\varepsilon_j}>0$ such that $\phi_j(x)^{-\kappa_j} \le \varepsilon_j\phi_j(x)^{-\gamma_{2,j}} +
			C_{\varepsilon_j}$. Consequently,
			\begin{flalign*}
				\int_0^T\rho(X_s)\,\mathrm{d}s\le C_T+\tilde{d}_0\int_0^T|X_s-x^{\circ}|^{\gamma_1}\mathrm ds+\sum_{j\in\mathcal J_\rho}\tilde{d}_j\int_0^T	\phi_j(X_s)^{-\gamma_{2,j}}\,\mathrm ds,
			\end{flalign*}
			where
			\begin{flalign*}
				\tilde{d}_0=
				\begin{cases}
					d_0,
					& \mathcal D\text{ is unbounded and }
					\kappa_0=\gamma_1>0,\\
					\varepsilon_0d_0,
					& \mathcal D\text{ is unbounded and }
					0<\kappa_0<\gamma_1,\\
					0,
					& \text{otherwise},
				\end{cases}
				\quad\text{and}\quad
				\tilde{d}_j=
				\begin{cases}
					d_j,
					& \kappa_j=\gamma_{2,j}>0,\\
					\varepsilon_jd_j,
					& 0<\kappa_j<\gamma_{2,j},\\
					0,
					& \kappa_j=0.
				\end{cases}
			\end{flalign*}
			In the preceding bound, the term involving $\gamma_1$
			is omitted when $\tilde d_0=0$. By the definition of $\mathscr X_\rho^\dagger$ and $\mathscr X<\mathscr X_\rho^\dagger$, if $\delta_{0,c}=1$, then $\frac{\mathscr X d_0}{p_0K_2} +
			\sum_{j\in\mathcal J_{\rho,c}} \frac{\mathscr Xd_j}{\mathscr X_j^\dagger} <1$, whereas, if $\delta_{0,c}=0$, then $\sum_{j\in\mathcal J_{\rho,c}} \frac{\mathscr Xd_j}{\mathscr X_j^\dagger} <1$. Hence, $\varepsilon_0$ and $\varepsilon_j$ can be chosen sufficiently small. If $\tilde{d}_0>0$, then $K_2>0$, and the above choices can be made so that $\frac{\mathscr X\tilde{d}_0}{p_0K_2} + \sum_{j\in\mathcal J_\rho:\,\tilde{d}_j>0}
			\frac{\mathscr X\tilde{d}_j}{\mathscr X_j^\dagger} <1$.
			Hence there exist positive numbers $\theta_0$ and
			$\theta_j$, $j\in\mathcal J_\rho$ with $\tilde{d}_j>0$, such that
			$\theta_0+\sum_{j:\,\tilde{d}_j>0}\theta_j=1$, and $\frac{\mathscr X\tilde{d}_0}{\theta_0}<p_0K_2$, 
			$\frac{\mathscr X\tilde{d}_j}{\theta_j}<\mathscr X_j^\dagger$.
			Therefore, by H\"older's inequality and Lemmas \ref{Lm4.1}
			and \ref{Lm4.2},
			\begin{align*}
				\mathbb E\left[ \exp\left(
				\mathscr X\int_0^T\rho(X_s)\,\mathrm ds \right)
				\right] \le& C \left[ \mathbb E \exp\left(
				\frac{\mathscr X\tilde{d}_0}{\theta_0} \int_0^T|X_s-x^{\circ}|^{\gamma_1}\,\mathrm ds \right)
				\right]^{\theta_0}\\
				&\cdot \prod_{\substack{j\in\mathcal J_\rho\\ \tilde{d}_j>0}}
				\left[ \mathbb E
				\exp\left( \frac{\mathscr X\tilde{d}_j}{\theta_j}
				\int_0^T \phi_j(X_s)^{-\gamma_{2,j}}\,\mathrm ds
				\right) \right]^{\theta_j}\le C.
			\end{align*}
			
			If $\tilde{d}_0=0$ and at least one $\tilde d_j>0$, the parameters chosen above satisfy $\sum_{\substack{j\in\mathcal J_\rho\\\tilde{d}_j>0}} \frac{\mathscr X\tilde{d}_j}{\mathscr X_j^\dagger}<1$. Hence there exist positive numbers $\theta_j$,
			$j\in\mathcal J_\rho$ with $\tilde{d}_j>0$, such that $\sum_{\substack{j\in\mathcal J_\rho\\\tilde{d}_j>0}}\theta_j=1$ and $\frac{\mathscr X\tilde{d}_j}{\theta_j}<\mathscr X_j^\dagger$.
			Then H\"older's inequality and Lemma \ref{Lm4.2} yield
			\begin{align*}
				\mathbb E\left[ \exp\left(
				\mathscr X\int_0^T\rho(X_s)\,\mathrm ds \right)
				\right]\le C \prod_{\substack{j\in\mathcal J_\rho\\ \tilde{d}_j>0}} \left[ \mathbb E\exp\left(
				\frac{\mathscr X\tilde{d}_j}{\theta_j} \int_0^T
				\phi_j(X_s)^{-\gamma_{2,j}}\,\mathrm ds
				\right) \right]^{\theta_j} \le C.
			\end{align*}
			
			Finally, if $\tilde d_0=0$ and $\tilde d_j=0$ for every
			$j\in\mathcal J_\rho$, then
			$\int_0^T\rho(X_s)\,\mathrm ds\le C_T$. Hence,
			\[
			\mathbb E\left[
			\exp\left(\mathscr X\int_0^T\rho(X_s)\,\mathrm ds\right)
			\right]
			\le e^{\mathscr X C_T}\leq C.
			\]
			This completes the proof. \hfill
		\end{proof}	
	\end{lemma}

	\section{Domain-preserving projected Milstein scheme for stochastic systems}
	\label{sec3}
	\subsection{The projected Milstein scheme and weighted assumptions}
	\label{sec3.1}
	We now construct the domain-preserving scheme. For $x\in\mathcal D$, define $H_{lk}(x) := Dg_k(x)g_l(x)$, $ l,k=1,2,\ldots,m$. For $n=0,1,\ldots,N-1$, set $\Delta B_n^l:=B_{t_{n+1}}^l-B_{t_n}^l$ and $\Xi_{lk,n} := \Delta B_n^l \Delta B_n^k - \delta_{lk}h$. Here $\delta_{lk}$ is the Kronecker delta. Throughout this section, we assume that the commutative noise condition holds: $H_{lk}(x) = H_{kl}(x), \, l, k = 1, 2,\ldots, m$.
	
	Define the one-step Milstein map by
	\begin{flalign*}
		\mathcal M_{h,n}(x) = x + h f(x) + \sum_{k=1}^m g_k(x) \Delta B_n^k + \frac{1}{2} \sum_{l,k=1}^m H_{l,k}(x) \Xi_{lk,n}.
	\end{flalign*}
	Take $0 < \alpha < \frac{1}{2}$. If $\mathcal{D}$ is unbounded, then take $\beta > 0$ and define
	\begin{flalign*}
		\mathcal{D}_h = \{ x \in \bar{\mathcal{D}} : \phi_j(x) \ge h^\alpha, j=1,2,\ldots,J,\quad |x-x^{\circ}| \le h^{-\beta} \}.
	\end{flalign*}
	If $\mathcal{D}$ is bounded, then 
	\begin{flalign*}
		\mathcal{D}_h = \{ x \in \bar{\mathcal{D}} : \phi_j(x) \ge h^\alpha, j=1,2,\ldots,J \}.
	\end{flalign*}
	Since $X_0\in\mathcal D$, there exists $h_0>0$ such that
	$X_0\in\mathcal D_h$ for every $0<h\le h_0$. For such $h$, $\mathcal D_h$ is a nonempty closed convex subset of $\mathcal D$, so its Euclidean projection is well defined, single-valued and satisfies
	\begin{flalign*}
		|\Pi_{\mathcal D_h}(x)-\Pi_{\mathcal D_h}(y)| \le |x-y|, \quad x,y\in\mathbb R^d.
	\end{flalign*}  
	The projected Milstein scheme is defined by
	\begin{flalign*}
		Y_{n+1} = \Pi_{\mathcal{D}_h} (\mathcal M_{h,n}(Y_n)), \quad Y_0 = X_0.
	\end{flalign*}
	By construction, $Y_n\in\mathcal D_h\subset\mathcal D$ for $n=0,1,\ldots,N$ pathwise. For $z=\mathcal M_{h,n}(Y_n)$, computing $\Pi_{\mathcal D_h}(z)$ amounts to minimizing $\frac 12|y-z|^2$ subject to $\phi_j(y)\ge h^\alpha$, $j=1,2,\ldots,J$, with the additional constraint $|y-x^\circ|\le h^{-\beta}$ when $\mathcal D$ is unbounded. The bounded-domain problem is a strictly convex quadratic program with linear constraints; the additional ball constraint in the unbounded case gives a convex quadratically constrained problem. All coefficients in the Milstein step are evaluated at the known state $Y_n$, while the projection may require numerical optimization for a general domain. Projection procedures for the numerical examples are described in Section \ref{sec5}.         
	
	For $\mu\geq 0$ and $\eta=(\eta_1,\eta_2,\ldots,\eta_J)\in[0,\infty)^{J}$, define
	\begin{flalign*}
		\mathcal{A}_{\mu, \eta}(x) = 1 + |x-x^{\circ}|^\mu + \sum_{j: \eta_j > 0} \bar{\phi}_j(x)^{-\eta_j},
	\end{flalign*}
	where $ \bar{\phi}_j(x) = 1 \wedge \phi_j(x)$. 
	
	Let the labels $I$, $R$ and $S$ refer to the increment estimate,
	the local Milstein remainder and the stability analysis,
	respectively. For each $L\in\{I,R,S\}$, choose $\mu_L\ge0$ and
	$\eta_L=(\eta_{L,1},\eta_{L,2},\ldots,\eta_{L,J})\in[0,\infty)^J$.
	Set
	\begin{flalign*}
		\mathcal A_L(x):=\mathcal A_{\mu_L,\eta_L}(x),\quad x\in\mathcal D.
	\end{flalign*}
	
	We write $g(x)=(g_1(x),g_2(x),\ldots,g_m(x))\in\mathbb R^{d\times m}$. For a sufficiently smooth function $u$ on $\mathcal{D}$, define 
	\begin{flalign*}
		L^0 u(x) = Du(x) f(x) + \frac{1}{2} \sum_{l=1}^m D^2 u(x)[g_l(x), g_l(x)]\quad \text{and}\quad
		L^l u(x) = Du(x) g_l(x), \quad l = 1,2,\ldots, m.
	\end{flalign*}
	\begin{assumption}\label{ass4}
		Suppose that there exists $C>0$ such that for all $x\in\mathcal{D}$,
		\begin{flalign*}
			\begin{cases}
				|f(x)| + \|g(x)\|_F \le C \mathcal{A}_I(x), \\
				|L^0 f(x)| + \left( \sum_{l} |L^l f(x)|^2 \right)^{\frac 12} + \left( \sum_k |L^0 g_k(x)|^2 \right)^{\frac 12} \le C \mathcal{A}_R(x), \\
				\|D f(x)\| + \left( \sum_{l} \|D g_l(x)\|^2 \right)^{\frac 12} + \left( \sum_{l,k} \|D H_{l,k}(x)\|^2 \right)^{\frac 12} \le C \mathcal{A}_S(x).
			\end{cases}
		\end{flalign*}
	\end{assumption}
	\begin{assumption}\label{ass5}
		Suppose that there exist $C>0$ and $ \bar{q} > 4$ such that for any $x \in \mathcal{D}$ and $ v \in \mathbb{R}^d$,
		\begin{flalign*}
			\langle v, D f(x) v \rangle + \frac{\bar{q}-1}{2} \sum_{l=1}^m | D g_l(x) v |^2 \le C\rho(x) |v|^2.
		\end{flalign*}
	\end{assumption}
	
	The following lemma provides the moment bounds for the three weights  in Assumption~\ref{ass4}.
	\begin{lemma}\label{Lm4.4}
		Suppose the assumptions in Lemma \ref{Lm4.1} hold. Fix $q\in[4,\bar q)$ and set $m_I:=\max\{4q,\bar q\}$. Assume the following conditions hold.
		\begin{enumerate}
			\item[(1)] If $\mathcal D$ is bounded, then assume that $\max\{m_I\eta_{I,j},\,2q\eta_{R,j},\,4q\eta_{S,j}\}
			<q_j$ for $j=1,2,\ldots,J$.
			\item[(2)] If $\mathcal{D}$ is unbounded, then assume that
			\[
			\max\{m_I\eta_{I,j},\,2q\eta_{R,j},\,4q\eta_{S,j}\}
			<q_j^\dagger,\qquad j=1,2,\ldots,J,
			\]
			and that $\max\{m_I\mu_I,2q\mu_R,4q\mu_S\}<p_0$.
		\end{enumerate}
		Then there exists $C>0$ such that
		\begin{flalign}\label{IRS}
			\mathbb E\Big[\sup_{0\le t\le T}\mathcal A_I(X_t)^{m_I}\Big]+\mathbb E\Big[\sup_{0\le t\le T}\mathcal A_R(X_t)^{2q}\Big]+\mathbb E\Big[\sup_{0\le t\le T}\mathcal A_S(X_t)^{4q}\Big]\le C.
		\end{flalign}
	\end{lemma}
	\begin{proof}
		Let $(\ell,a)\in\{(I,m_I),(R,2q),(S,4q)\}$. Since $a\ge1$ and $\bar\phi_j(x)^{-a\eta_{\ell,j}}\le1+\phi_j(x)^{-a\eta_{\ell,j}}$, we have
		\begin{flalign*}
			\mathcal A_\ell(x)^a \le C\Big[ 1+\vert{}x-x^{\circ}\vert{}^{a\mu_\ell} + \sum_{j:\eta_{\ell,j}>0}\phi_j(x)^{-a\eta_{\ell,j}} \Big].
		\end{flalign*}
		The assumptions give $a\eta_{\ell,j}<q_j^\dagger$ whenever $\eta_{\ell,j}>0$. Hence the required boundary inverse-moment bounds follow from Lemma~\ref{Lm4.1}.
		
		The term $|x-x^\circ|^{a\mu_\ell}$ is bounded when $\mathcal D$ is bounded or $\mu_\ell=0$; otherwise, $0<a\mu_\ell<p_0$,
		and Lemma~\ref{Lm4.1} provides its moment bound. Consequently,
		\begin{align*}
			\mathbb E\Big[
			\sup_{0\le t\le T}\mathcal A_\ell(X_t)^a
			\Big]\le C\Big(1+\mathbb E\Big[
			\sup_{0\le t\le T}|X_t-x^\circ|^{a\mu_\ell}
			\Big]+\sum_{j:\eta_{\ell,j}>0}\mathbb E\Big[\sup_{0\le t\le T}\phi_j(X_t)^{-a\eta_{\ell,j}}\Big]\Big)\le C.
		\end{align*}
		Summing the estimates for the three choices of $(\ell,a)$
		proves~\eqref{IRS}.\hfill
	\end{proof}

	\subsection{Local consistency estimates}
	\label{sec3.2}
	Lemma~\ref{Lm4.4} and the coefficient growth bound in
	Assumption~\ref{ass4} yield the following short-time
	increment estimate for the exact solution.
	\begin{lemma}\label{Lm4.5}
		Suppose Assumption \ref{ass4} and the assumptions of Lemma \ref{Lm4.4} hold. Then, for every $2\le a\le m_I$, there exists a constant $C>0$ such that
		\begin{equation}
			\mathbb E\Big[ \sup_{s\le u\le t}|X_u-X_s|^a \Big] \le C|t-s|^{\frac{a}{2}}, \qquad 0\le s\le t\le T.
		\end{equation}
		\begin{proof}
			By Lemma \ref{Lm4.4},
			\begin{equation*}
				\mathbb E\Big[ \sup_{0\le t\le T}\mathcal A_I(X_t)^{m_I} \Big]\leq C.
			\end{equation*}
			Since $2\le a\le m_I$, Assumption \ref{ass4}, together with Lemma \ref{Lm4.4}, yields
			\begin{equation*}
				\sup_{0\le v\le T} \mathbb E\left[ |f(X_v)|^a+\|g(X_v)\|_F^a \right] \le C.
			\end{equation*}
			Indeed, we have $|f(x)|+\|g(x)\|_F \le C\mathcal A_I(x)$
			and 
			\begin{equation*}
				|f(X_v)|^a+\|g(X_v)\|_F^a \le C\mathcal A_I(X_v)^a \le C\left(1+\mathcal A_I(X_v)^{m_I}\right).
			\end{equation*}
			For $s\le u\le t$, we obtain
			\begin{equation*}
				X_u-X_s = \int_s^u f(X_v)\,\mathrm dv + \sum_{l=1}^m \int_s^u g_l(X_v)\,\mathrm dB_v^l.
			\end{equation*}
			Therefore, by H\"older's inequality and the Burkholder--Davis--Gundy inequality,
			\begin{flalign*}
				\mathbb E\Big[ \sup_{s\le u\le t}|X_u-X_s|^a \Big] \le C\mathbb E \Big( \int_s^t|f(X_v)|\,\mathrm dv \Big)^a+ C\mathbb E \Big( \int_s^t\|g(X_v)\|_F^2\,\mathrm dv \Big)^{\frac{a}{2}}.
			\end{flalign*}
			H\"older's inequality gives
			\begin{equation*}
				\mathbb E \Big( \int_s^t|f(X_v)|\,\mathrm dv \Big)^a \le |t-s|^{a-1} \int_s^t \mathbb E|f(X_v)|^a\,\mathrm dv \le C|t-s|^a,
			\end{equation*}
			and
			\begin{flalign*}
				\mathbb E \Big( \int_s^t\|g(X_v)\|_F^2\,\mathrm dv \Big)^{\frac{a}{2}} \le |t-s|^{\frac{a}{2}-1} \int_s^t \mathbb E\|g(X_v)\|_F^a\,\mathrm dv\le C|t-s|^{\frac{a}{2}}.
			\end{flalign*}
			Consequently,
			\begin{equation*}
				\mathbb E\Big[ \sup_{s\le u\le t}|X_u-X_s|^a \Big] \le C|t-s|^a + C|t-s|^{\frac{a}{2}} \le C|t-s|^{\frac{a}{2}}.
			\end{equation*}
			This completes the proof. \hfill
		\end{proof}
	\end{lemma}

	To establish local consistency, set $X_n=X(t_n)$ and define
	the one-step Milstein remainder by
	$\mathcal R_{n+1}:=X_{n+1}-\mathcal M_{h,n}(X_n)$.
	Using the weight moment bounds in Lemma~\ref{Lm4.4},
	the increment estimate in Lemma~\ref{Lm4.5}, and the
	coefficient and derivative bounds in Assumption~\ref{ass4},
	we obtain the following remainder estimates.
	\begin{lemma}\label{Lm4.6}
		Suppose that Assumption \ref{ass4} and the assumptions of Lemma \ref{Lm4.4} hold. 
		Then there exists a constant $C>0$, independent of $h$ and $n$, such that
		\begin{equation}\label{4.15}
			\|\mathcal  R_{n+1}\|_{L^{2q}}
			\leq Ch^{\frac 32}\quad\text{and}\quad \|\mathbb E_n \mathcal  R_{n+1}\|_{L^{2q}}
			\leq Ch^2.
		\end{equation}
		Moreover, if $k_n:=1+h^{-2}|\mathbb E_n \mathcal R_{n+1}|+h^{-\frac 32}
		\left(\mathbb E_n|\mathcal R_{n+1}|^q\right)^{1/q}$,
		then
		\begin{equation}\label{4.16}
			\sup_{h,n}\mathbb E[k_n^{2q}]\leq C.
		\end{equation}
		\begin{proof}
			We decompose $\mathcal R_{n+1} = \mathcal R_{n+1}^f+\mathcal R_{n+1}^g$, where
			\begin{flalign*}
				\mathcal R_{n+1}^f = \int_{t_n}^{t_{n+1}} \bigl(f(X_t)-f(X_n)\bigr)\,\mathrm dt\,
			\end{flalign*}
			and
			\begin{flalign*}
				\mathcal R_{n+1}^g = \sum_{l=1}^m \int_{t_n}^{t_{n+1}} \bigl(g_l(X_t)-g_l(X_n)\bigr)\,\mathrm dB_t^l-\frac12\sum_{l,k=1}^mH_{lk}(X_n)\Xi_{lk,n}.
			\end{flalign*}
			By It\^o's formula, we have
			\begin{flalign*}
				f(X_t)-f(X_n) = \int_{t_n}^t L^0 f(X_u)\,\mathrm du + \sum_{l=1}^m \int_{t_n}^t L^l f(X_u)\,\mathrm dB_u^l.
			\end{flalign*}
			Hence,
			\begin{equation*}
				\mathcal R_{n+1}^f=\int_{t_n}^{t_{n+1}}\int_{t_n}^tL^0f(X_u)\,\mathrm du\,\mathrm dt +\sum_{l=1}^m\int_{t_n}^{t_{n+1}}
				\int_{t_n}^tL^lf(X_u)\,\mathrm dB_u^l\,\mathrm dt.
			\end{equation*}
			By Assumption \ref{ass4} and Lemma \ref{Lm4.4},
			\begin{flalign*}
				\sup_{0\leq u\leq T}\left\|L^0f(X_u)\right\|_{L^{2q}}
				\leq C \quad \text{and} \quad \sup_{0\leq u\leq T}\Big\|\Big(\sum_{l=1}^m|L^lf(X_u)|^2\Big)^{\frac 12}\Big\|_{L^{2q}} \leq C.
			\end{flalign*}
			Therefore, by Minkowski's inequality,
			\begin{equation}\label{L0f}
				\Big\|\int_{t_n}^{t_{n+1}}\int_{t_n}^tL^0f(X_u)\,\mathrm du\,\mathrm dt\Big\|_{L^{2q}}\leq\int_{t_n}^{t_{n+1}}\int_{t_n}^t
				\|L^0f(X_u)\|_{L^{2q}}\,\mathrm du\,\mathrm dt\leq Ch^2.
			\end{equation}
			By the stochastic Fubini theorem,
			\begin{flalign}\label{Llf}
				\sum_{l=1}^m\int_{t_n}^{t_{n+1}}\int_{t_n}^tL^lf(X_u)\,\mathrm dB_u^l\,\mathrm dt=\sum_{l=1}^m\int_{t_n}^{t_{n+1}}(t_{n+1}-u)L^lf(X_u)\,\mathrm dB_u^l.
			\end{flalign}
			The Burkholder--Davis--Gundy inequality gives
			\begin{equation*}
				\Big\|\sum_{l=1}^m\int_{t_n}^{t_{n+1}}(t_{n+1}-u)L^lf(X_u)\,\mathrm dB_u^l\Big\|_{L^{2q}}\leq C\Big[\int_{t_n}^{t_{n+1}}(t_{n+1}-u)^2\Big\|\Big(\sum_{l=1}^m|L^lf(X_u)|^2\Big)^{\frac 12}\Big\|_{L^{2q}}^2\,\mathrm du \Big]^{\frac 12}\leq Ch^{\frac 32}.
			\end{equation*}
			Combining this and \eqref{L0f}, we obtain
			\begin{equation}\label{Rnf}
				\|\mathcal R_{n+1}^f\|_{L^{2q}}\leq Ch^{\frac 32}.
			\end{equation}
			Moreover, the stochastic integral in \eqref{Llf} has zero conditional expectation with respect to $\mathcal F_{t_n}$. Hence, 
			\begin{equation}\label{ERnf}
				\|\mathbb E_n \mathcal R_{n+1}^f\|_{L^{2q}}\leq\Big\|\int_{t_n}^{t_{n+1}}
				\int_{t_n}^tL^0f(X_u)\,\mathrm du\,\mathrm dt\Big\|_{L^{2q}} \leq
				Ch^2.
			\end{equation}
			
			Next, applying It\^o's formula to $g_k(X_t)$ gives
			\begin{equation*}
				g_k(X_t) = g_k(X_n) + \int_{t_n}^t L^0 g_k(X_u)\,\mathrm du + \sum_{l=1}^m \int_{t_n}^t H_{lk}(X_u)\,\mathrm dB_u^l.
			\end{equation*}
			Using $	H_{lk}(X_u) = H_{lk}(X_n) + \bigl(H_{lk}(X_u)-H_{lk}(X_n)\bigr)$
			and the commutativity condition $H_{lk}=H_{kl}$, we obtain
			\begin{flalign*}
				\sum_{l,k=1}^m H_{lk}(X_n) \int_{t_n}^{t_{n+1}}\int_{t_n}^t \,\mathrm dB_u^l\,\mathrm dB_t^k = \frac12 \sum_{l,k=1}^m H_{lk}(X_n)\Xi_{lk,n}.
			\end{flalign*}
			Substituting this identity into the definition of $\mathcal R^g_{n+1}$ yields
			\begin{equation}\label{Rngg}
				\mathcal R_{n+1}^g=\sum_{k=1}^m\int_{t_n}^{t_{n+1}}\Big[\int_{t_n}^tL^0g_k(X_u)\,\mathrm du +\sum_{l=1}^m\int_{t_n}^t\bigl(H_{lk}(X_u)-H_{lk}(X_n)
				\bigr)\,\mathrm dB_u^l\Big]\,\mathrm dB_t^k.
			\end{equation}
			For the first term in \eqref{Rngg}, the BDG inequality, Assumption \ref{ass4}, and Lemma \ref{Lm4.4} yield
			\begin{equation*}
				\Big\|\sum_{k=1}^m\int_{t_n}^{t_{n+1}}\int_{t_n}^t
				L^0g_k(X_u)\,\mathrm du\,\mathrm dB_t^k\Big\|_{L^{2q}} \leq C\Big[\int_{t_n}^{t_{n+1}}\Big(\int_{t_n}^t\Big\|\Big( \sum_{k=1}^m|L^0g_k(X_u)|^2\Big)^{\frac 12}\Big\|_{L^{2q}}\,\mathrm du\Big)^2\mathrm dt\Big]^{\frac 12}\leq Ch^{\frac 32}.
			\end{equation*}
			It remains to estimate $H_{lk}(X_u)-H_{lk}(X_n)$.
			For $\theta\in[0,1]$, define $
			Z_{\theta,u}:=(1-\theta)X_n+\theta X_u$.
			Since $\mathcal D$ is convex, $Z_{\theta,u}\in\mathcal D$. The affine form of $\phi_j$ gives $\phi_j(Z_{\theta,u})=(1-\theta)\phi_j(X_n)+\theta\phi_j(X_u)$, and hence $\bar{\phi_j}(Z_{\theta,u})^{-\eta_j}\leq (1-\theta)\bar{\phi}_j(X_n)^{-\eta_j}+\theta\bar{\phi}_j(X_u)^{-\eta_j}$ for $\eta_j\geq 0$. By the integral mean-value formula,
			\begin{equation*}
				H_{lk}(X_u)-H_{lk}(X_n) = \int_0^1 DH_{lk}(Z_{\theta,u})(X_u-X_n)\,\mathrm d\theta.
			\end{equation*}
			On the other hand, for any $\mu\geq 0$,
			\begin{flalign*}
				|Z_{\theta,u}-x^{\circ}|^\mu\leq C(|X_n-x^{\circ}|^\mu+|X_u-x^{\circ}|^\mu).
			\end{flalign*}
			Then we obtain
			\begin{equation*}
				\mathcal A_S(Z_{\theta,u})\leq C\bigl(
				\mathcal A_S(X_n)+\mathcal A_S(X_u)\bigr).
			\end{equation*}
			Hence Assumption \ref{ass4} and Lemma \ref{Lm4.4} imply
			\begin{equation*}
				\sup_{\theta\in[0,1]}\|DH_{lk}(Z_{\theta,u})\|_{L^{4q}}\leq C.
			\end{equation*}
			Since $4q\leq m_I$, Lemma \ref{Lm4.5} gives
			\begin{equation*}
				\|X_u-X_n\|_{L^{4q}}\leq C(u-t_n)^{\frac 12}.
			\end{equation*}
			Therefore, by H\"older's inequality,
			\begin{flalign*}
				\|H_{lk}(X_u)-H_{lk}(X_n)\|_{L^{2q}}
				&\leq\int_0^1\|DH_{lk}(Z_{\theta,u})(X_u-X_n)\|_{L^{2q}}
				\,\mathrm d\theta \\&\leq\int_0^1
				\|DH_{lk}(Z_{\theta,u})\|_{L^{4q}}\|X_u-X_n\|_{L^{4q}}\,\mathrm d\theta\leq C(u-t_n)^{\frac 12}.
			\end{flalign*}
			Applying the BDG inequality twice and using the inequality above, we obtain
			\begin{flalign*}
				&\Big\|\sum_{l,k=1}^m\int_{t_n}^{t_{n+1}}
				\int_{t_n}^t\bigl(H_{lk}(X_u)-H_{lk}(X_n)
				\bigr)\,\mathrm dB_u^l\,\mathrm dB_t^k\Big\|_{L^{2q}}^2 \\
				\leq & C\int_{t_n}^{t_{n+1}}\int_{t_n}^t\sum_{l,k=1}^m
				\|H_{lk}(X_u)-H_{lk}(X_n)\|_{L^{2q}}^2\,\mathrm du\,\mathrm dt \leq C \int_{t_n}^{t_{n+1}}\int_{t_n}^t(u-t_n)\,\mathrm du\,\mathrm dt\leq Ch^3.
			\end{flalign*}
			Thus,
			\begin{equation}\label{Rng}
				\|\mathcal R_{n+1}^g\|_{L^{2q}}\leq Ch^{\frac 32}.
			\end{equation}
			Since the integrands in~\eqref{Rngg}  are adapted and square integrable, each term is an It\^o integral over $[t_n,t_{n+1}]$. Hence
			\begin{equation}\label{ERng}
				\mathbb E_n [\mathcal R_{n+1}^g]=0.
			\end{equation}
			Combining \eqref{Rnf}, \eqref{ERnf}, \eqref{Rng}, and \eqref{ERng}, we obtain
			\begin{flalign*}
				\|\mathcal R_{n+1}\|_{L^{2q}}\leq Ch^{\frac 32}\quad \text{and}\quad
				\|\mathbb E_n \mathcal R_{n+1}\|_{L^{2q}}\leq Ch^2,
			\end{flalign*}
			which proves \eqref{4.15}. Finally, by this, the elementary inequality
			$(a+b+c)^{2q}\leq C(a^{2q}+b^{2q}+c^{2q})$, and the conditional Jensen inequality, we have
			\begin{flalign*}
				\mathbb E[k_n^{2q}]\leq& C+C h^{-4q}\mathbb E|\mathbb E_n \mathcal R_{n+1}|^{2q} +C h^{-3q}\mathbb E[\left(\mathbb E_n|\mathcal R_{n+1}|^q\right)^2]\\
				\leq& C+Ch^{-4q}\|\E_n\mathcal R_{n+1}\|^{2q}_{L^{2q}} +Ch^{-3q}\E[|\mathcal R_{n+1}|^{2q}]\leq C.
			\end{flalign*}
			We conclude that
			$\sup_{h,n}\mathbb E[k_n^{2q}]\leq C$, completing  the proof. \hfill
		\end{proof}
	\end{lemma}

	\subsection{Localization and stability estimates}
	\label{sec3.3}
	The coefficient derivatives in Assumption~\ref{ass4} need not be globally bounded. We therefore localize the exact solution and the numerical error to obtain step-size-dependent derivative bounds along $[X_n,Y_n]$. Relative neighborhoods are also introduced to compare
	the stability weight at grid points with its values along the exact trajectory. For sufficiently small $h>0$, if $\mathcal D$ is unbounded, then set
	\begin{flalign*}
		K_h := \{ x\in\mathcal D: \phi_j(x)\geq 8h^\alpha,\ j=1,2,\ldots,J, \, |x-x^{\circ}|\leq \frac18h^{-\beta}\}.
	\end{flalign*}
	If $\mathcal D$ is bounded, then set
	\begin{flalign*}
		K_h := \{ x\in\mathcal D: \phi_j(x)\geq 8h^\alpha,\ j=1,2,\ldots,J\}.
	\end{flalign*}
	For both bounded and unbounded domains, fix $0<\varepsilon<\frac{1}{16}$.
	For $x\in \mathcal D$, define the relative neighborhood
	\begin{flalign*}
		\mathcal T_\varepsilon(x) := \left\{ y\in\mathcal D: \vert{}y-x\vert{}\leq\varepsilon(1+|x-x^{\circ}|), \, (1-\varepsilon)\phi_j(x) \leq \phi_j(y) \leq (1+\varepsilon)\phi_j(x),\ j=1,2,\ldots,J \right\}.
	\end{flalign*}
	For later use, define
	\begin{flalign*}
		\delta_S:=
		\begin{cases}
			\max\big\{\beta\mu_S, \max_{1\leq j\leq J}\alpha\eta_{S,j}
			\big\},
			& \text{if }\mathcal D\text{ is unbounded},\\
			\max_{1\leq j\leq J}\alpha\eta_{S,j},
			& \text{if }\mathcal D\text{ is bounded}.
		\end{cases}
	\end{flalign*}
	With the convention $\min\varnothing=N+1$, define
	\[
	\sigma_X:=
	\min\left\{
	\begin{aligned}
		&\min\bigl\{k\in\{0,1,\ldots,N\}:X_k\notin K_h\bigr\},\\
		&\min\bigl\{k\in\{1,2,\ldots,N\}:
		\exists\,t\in[t_{k-1},t_k],\
		X_t\notin\mathcal T_\varepsilon(X_{k-1})\bigr\}
	\end{aligned}
	\right\}.
	\]
	Let $e_n:=X_n-Y_n$ and define
	\[
	\sigma_E:=\min\bigl\{k\in\{0,1,\ldots,N\}:
	k<\sigma_X\ \text{and}\ |e_k|>h^\alpha\bigr\}.
	\]
	Define $F_n:=f(X_n)-f(Y_n)$, $G_{l,n}:=g_l(X_n)-g_l(Y_n)$, and $M_{lk,n}:=H_{lk}(X_n)-H_{lk}(Y_n)$ with $|G_n|^2:=\sum_{l=1}^m|G_{l,n}|^2$, $|M_n|^2:=\sum_{l,k=1}^m|M_{lk,n}|^2$. The following lemma establishes the geometric implications of this localization, yielding explicit derivative bounds along the convex combination of the exact and numerical solutions.
	\begin{lemma}\label{Lm4.7}
		Suppose Assumption \ref{ass4} holds. Then for all sufficiently small $h>0$, the following statements hold.
		\begin{enumerate}
			\item[(1)]	If $n<\sigma_X$ and $n<\sigma_E$, then $X_n\in K_h$ and $|e_n|\le h^\alpha$. For every $z\in[X_n,Y_n]$, we have  $\phi_j(z)\ge4h^\alpha$, $j=1,2,\ldots,J$ and
			\[
			|z-x^\circ|\le\frac14h^{-\beta}
			\]
			when $\mathcal D$ is unbounded. Consequently,
			\begin{flalign*}
				\sup_{z\in[X_n,Y_n]} \Big[ \|Df(z)\| +\Big(\sum_{l=1}^m\|Dg_l(z)\|^2\Big)^{\frac 12} +\Big(\sum_{l,k=1}^m\|DH_{lk}(z)\|^2\Big)^{\frac 12} \Big] \le Ch^{-\delta_S}
			\end{flalign*}
			and
			\begin{flalign*}
				|F_n|+|G_n|+|M_n|\le Ch^{-\delta_S}|e_n|.
			\end{flalign*}
			\item[(2)] If $n+1<\sigma_X$, then
			$X_{n+1}\in K_h\subset\mathcal D_h$, and hence
			\begin{flalign*}
				|e_{n+1}| \le |X_{n+1}-\mathcal M_{h,n}(Y_n)|.
			\end{flalign*} 
		\end{enumerate}
		\begin{proof}
			On $\{n < \sigma_X, n < \sigma_E\}$, the definitions of $\sigma_X$ and $\sigma_E$ ensure that $X_n \in K_h$ and $|e_n| \le h^\alpha$. Let $z = Y_n + \theta e_n$, $\theta \in [0,1]$. Since $Y_n \in \mathcal{D}_h$ and $\phi_j(Y_n) \geq h^\alpha > 0$ for all $j$, we have $Y_n \in \mathcal{D}$. Moreover, $n < \sigma_X$ implies $X_n \in K_h \subset \mathcal{D}$. Since $\mathcal{D}$ is convex and
			\begin{flalign*}
				z = (1 - \theta)Y_n + \theta X_n,
			\end{flalign*}
			it follows that $z \in \mathcal{D}$. Since $e_n=X_n-Y_n$, we have $ z=(1-\theta)Y_n+\theta X_n$ and therefore $|z-X_n| = (1-\theta)|e_n|\le h^\alpha$. Since each $\phi_j$ is affine with $|\nabla\phi_j|=|\nu_j|=1$, we obtain $|\phi_j(z)-\phi_j(X_n)| \leq |z-X_n|$. Since $X_n\in K_h$, $\phi_j(X_n)\geq 8h^\alpha$. Thus $\phi_j(z) \geq \phi_j(X_n)-\vert{}z-X_n\vert{} \geq 8h^\alpha-h^\alpha = 7h^\alpha \geq 4h^\alpha $. When $\mathcal D$ is unbounded, the triangle inequality gives $|z-x^{\circ}| \leq |X_n-x^{\circ}|+|z-X_n| \leq \frac18h^{-\beta}+h^\alpha$. Since $\alpha,\beta>0$, for all sufficiently small $h$, $h^\alpha\leq\frac18h^{-\beta}$. Hence $|z-x^{\circ}|\leq\frac14h^{-\beta}$. 
			
			We next estimate the growth weight $\mathcal A_S$ along the segment $[X_n,Y_n]$. For sufficiently small $h$, $\bar\phi_j(z)^{-\eta_{S,j}} \le C h^{-\alpha\eta_{S,j}}$. When $\mathcal D$ is unbounded, the inequality $|z-x^{\circ}|\le \frac14h^{-\beta}$ gives $|z-x^{\circ}|^{\mu_S} \le C h^{-\beta\mu_S}$. Therefore, by the definition of $\delta_S$,
			\begin{flalign*}
				\sup_{z\in[X_n,Y_n]}\mathcal A_S(z) \leq Ch^{-\delta_S}.
			\end{flalign*} 
			In the bounded-domain case, the radial term is uniformly bounded and the same conclusion follows with $\delta_S = \max_j\alpha\eta_{S,j}$. Assumption \ref{ass4} yields
			\begin{flalign}\label{XnYn}
				\sup_{z\in[X_n,Y_n]} \Big[ \|Df(z)\| +\Big(\sum_{l=1}^m\|Dg_l(z)\|^2\Big)^{\frac 12} +\Big(\sum_{l,k=1}^m\|DH_{lk}(z)\|^2\Big)^{\frac 12} \Big] \le Ch^{-\delta_S}.
			\end{flalign}
			By the integral mean-value formula, $F_n = \int_0^1 Df(Y_n+\theta e_n)e_n\,\mathrm d\theta $, and for each $l$, $G_{l,n} = \int_0^1 Dg_l(Y_n+\theta e_n)e_n\,\mathrm d\theta $. Similarly, $M_{lk,n} = \int_0^1 DH_{lk}(Y_n+\theta e_n)e_n\,\mathrm d\theta$. Therefore, using the inequality \eqref{XnYn}, we have $|F_n| \leq Ch^{-\delta_S}|e_n|$, $|G_n|\leq Ch^{-\delta_S}|e_n|$, and $|M_n|\leq Ch^{-\delta_S}|e_n|$. Consequently, $|F_n|+|G_n|+|M_n| \leq Ch^{-\delta_S}|e_n|$. 
			
			Finally, suppose $n+1<\sigma_X$. Then $X_{n+1}\in K_h$. By the definition of $K_h$, $\phi_j(X_{n+1})\geq 8h^\alpha\ge h^\alpha, \, j=1,2,\ldots,J$, and, when $\mathcal D$ is unbounded, $|X_{n+1}-x^{\circ}| \leq \frac18h^{-\beta} \le h^{-\beta} $. Hence $X_{n+1}\in K_h\subset\mathcal D_h $, so that $\Pi_{\mathcal D_h}X_{n+1}=X_{n+1}$. Since $Y_{n+1} = \Pi_{\mathcal D_h}\bigl(\mathcal M_{h,n}(Y_n)\bigr)$ and the metric projection onto the closed convex set $\mathcal D_h$ is nonexpansive,
			\begin{flalign*}
				|e_{n+1}|= |X_{n+1} - \Pi_{\mathcal D_h}\bigl(\mathcal M_{h,n}(Y_n)\bigr) |= |\Pi_{\mathcal D_h}(X_{n+1}) - \Pi_{\mathcal D_h}\bigl(\mathcal M_{h,n}(Y_n)\bigr) |\leq | X_{n+1}-\mathcal M_{h,n}(Y_n) |.
			\end{flalign*}
			This completes the proof. \hfill
		\end{proof}
	\end{lemma}

	Combining the remainder bounds in Lemma~\ref{Lm4.6} with the localization estimates in Lemma~\ref{Lm4.7} and the stability condition in Assumption~\ref{ass5}, we obtain the following one-step estimate for the localized error.
	\begin{lemma}\label{Lm4.8}
		Suppose that  the assumptions of Lemmas \ref{Lm4.6} and \ref{Lm4.7} hold. Assume further that Assumption \ref{ass5} and $\delta_S<\frac12$ hold.
		Define $\zeta_n:=hF_n+\sum_{l=1}^mG_{l,n}\Delta B_n^l+ \frac12\sum_{l,k=1}^mM_{lk,n}\Xi_{lk,n}$.
		Then there exists a constant $\Lambda_0>0$, independent of $h$ and $n$, such that, for all sufficiently small $h>0$,
		\begin{flalign}\label{one-step-stability}
			\mathbb E_n\left[|e_{n+1}|^q
			\mathbb I_{\{n+1<\sigma_X,\ n+1\leq\sigma_E\}}\right]\leq
			\mathbb I_{\{n<\sigma_X,\ n<\sigma_E\}}\big[\bigl(1+\Lambda_0 h\rho(X_n)\bigr)|e_n|^q+Ch^{q+1}(1+k_n^q)\big].
		\end{flalign}
		\begin{proof}
			Set $A_n:=\{n<\sigma_X,\ n<\sigma_E\},\,
			B_{n+1}:=\{n+1<\sigma_X,\ n+1\leq\sigma_E\}$.
			Clearly, $B_{n+1}\subset A_n$. Moreover, on $\{n+1<\sigma_X\}$, Lemma \ref{Lm4.7} and the nonexpansivity of the projection imply
			\begin{flalign*}
				|e_{n+1}| \leq|X_{n+1}-\mathcal M_{h,n}(Y_n)|.
			\end{flalign*}
			Since $\mathcal R_{n+1} = X_{n+1}-\mathcal M_{h,n}(X_n)$, we have
			\begin{flalign*}
				X_{n+1}-\mathcal M_{h,n}(Y_n) = e_n+\zeta_n+\mathcal R_{n+1}.
			\end{flalign*}
			Consequently,
			\begin{equation}\label{step-start}
				\mathbb E_n\left[
				|e_{n+1}|^q\mathbb I_{B_{n+1}} \right]
				\leq \mathbb I_{A_n} \mathbb E_n
				|e_n+\zeta_n+\mathcal R_{n+1}|^q.
			\end{equation}
			For $q\geq4$, there exists $C_q>0$ such that
			\begin{align}\label{Taylor-q}
				|x+y|^q	\leq|x|^q+q|x|^{q-2}\langle x,y\rangle+
				\frac{q(q-1)}2|x|^{q-2}|y|^2+C_q(|x|^{q-3}|y|^3+|y|^q).
			\end{align}
			We apply \eqref{Taylor-q} with $x=e_n$ and $y=\zeta_n+\mathcal R_{n+1}$.
			Choose $\delta>0$ sufficiently small such that
			\begin{equation*}
				(1+\delta)(q-1)<\bar q-1.
			\end{equation*}
			Then $ |\zeta_n+\mathcal R_{n+1}|^2\leq(1+\delta)|\zeta_n|^2+(1+\delta^{-1})|\mathcal R_{n+1}|^2$. Moreover, we have $\mathbb E_n\zeta_n=hF_n$,
			\begin{equation}\label{zeta2}
				\mathbb E_n|\zeta_n|^2 \leq h|G_n|^2 + Ch^2(|F_n|^2+|M_n|^2),
			\end{equation}
			and, for $a\in\{3,q\}$,
			\begin{equation}\label{zetaa}
				\mathbb E_n|\zeta_n|^a \leq C( h^a|F_n|^a
				+ h^{a/2}|G_n|^a + h^a|M_n|^a ).
			\end{equation}
			For $z_\theta:=Y_n+\theta e_n$, $0\leq\theta\leq1$, the integral mean-value formula gives
			\begin{flalign*}
				F_n = \int_0^1 Df(z_\theta)e_n\,\mathrm d\theta \quad \text{and}\quad
				G_{l,n} = \int_0^1 Dg_l(z_\theta)e_n\,\mathrm d\theta.
			\end{flalign*}
			Hence, by Jensen's inequality,
			\begin{flalign*}
				|G_n|^2 = \sum_{l=1}^m|G_{l,n}|^2\leq\int_0^1\sum_{l=1}^m
				|Dg_l(z_\theta)e_n|^2\,\mathrm d\theta.
			\end{flalign*}
			On the event $A_n$, Lemma \ref{Lm4.7} yields $ \phi_j(z_\theta) \geq
			\phi_j(X_n)-h^\alpha \geq \frac78\phi_j(X_n)$.
			Moreover, $|z_\theta-x^{\circ}| \leq |X_n-x^{\circ}|+h^\alpha \leq
			|X_n-x^{\circ}|+1$ for sufficiently small $h$. Therefore, by the definition of $\rho$,
			\begin{equation*}
				\rho(z_\theta)\leq C\rho(X_n), \,0\leq\theta\leq1.
			\end{equation*}
			Using this and applying Assumption \ref{ass5}, we obtain
			\begin{align*}
				\langle e_n,F_n\rangle+\frac{\bar q-1}{2}|G_n|^2\leq
				\int_0^1\Big[\langle e_n,Df(z_\theta)e_n\rangle
				+\frac{\bar q-1}{2}\sum_{l=1}^m|Dg_l(z_\theta)e_n|^2\Big]
				\,\mathrm d\theta\leq C\rho(X_n)|e_n|^2.
			\end{align*}
			Using this, we obtain
			\begin{align}\label{leading-drift-diff}
				q h|e_n|^{q-2}\langle e_n,F_n\rangle+\frac{q(q-1)}2(1+\delta)
				h|e_n|^{q-2}|G_n|^2\leq Ch\rho(X_n)|e_n|^q-c_0h|e_n|^{q-2}|G_n|^2,
			\end{align}
			where $c_0 := \frac q2 \left[ \bar q-1-(1+\delta)(q-1) \right] >0$.
			On the event $A_n$, Lemma \ref{Lm4.7} gives
			\begin{equation*}
				|F_n|+|G_n|+|M_n| \leq Ch^{-\delta_S}|e_n|.
			\end{equation*}
			Therefore,
			\begin{equation*}
				h^2|e_n|^{q-2}|F_n|^2 + h^2|e_n|^{q-2}|M_n|^2 \leq Ch^{2-2\delta_S}|e_n|^q \leq Ch|e_n|^q,
			\end{equation*}
			where the last inequality follows from $\delta_S<\frac 12$. Moreover,
			for $a\in\{3,q\}$,
			\begin{flalign*}
				h^a|e_n|^{q-a}(|F_n|^a+|M_n|^a)\leq Ch|e_n|^q.
			\end{flalign*}
			If \(e_n=0\), then $F_n=G_n=M_n=0$. For \(e_n\ne0\), we have
			\begin{align*}
				h^{a/2}|e_n|^{q-a}|G_n|^a=h|e_n|^{q-2}|G_n|^2\Big(
				\sqrt h\frac{|G_n|}{|e_n|}\Big)^{a-2}\leq Ch|e_n|^{q-2}|G_n|^2
				h^{(a-2)(\frac{1}{2}-\delta_S)}.
			\end{align*}
			For sufficiently small $h$, since $\delta_S<\frac12$, the latter terms can be absorbed by the negative diffusion term in \eqref{leading-drift-diff}. By the definition of $k_n$,
			\begin{equation*}
				|\mathbb E_n \mathcal R_{n+1}|\leq h^2k_n,\quad
				\mathbb E_n|\mathcal R_{n+1}|^a\leq
				\left(\mathbb E_n|\mathcal R_{n+1}|^q\right)^{\frac aq}\leq h^{\frac{3a}{2}}k_n^a, \quad 2\le a\leq q.
			\end{equation*}
			Since $e_n$ is $\mathcal F_{t_n}$-measurable, we have
			\begin{flalign*}
				\mathbb E_n\langle e_n,\mathcal R_{n+1}\rangle=\langle e_n,\mathbb E_n\mathcal R_{n+1}\rangle.
			\end{flalign*}
			Hence, for any $\tilde{\eta}>0$, Young's inequality yields
			\begin{flalign*}
				|e_n|^{q-1}|\mathbb E_n\mathcal R_{n+1}| \leq h^2k_n|e_n|^{q-1}
				\leq \tilde{\eta} h|e_n|^q + C_{\tilde{\eta}} h^{q+1}k_n^q.
			\end{flalign*}
			Moreover,
			\begin{flalign*}
				|e_n|^{q-2}\mathbb E_n|\mathcal R_{n+1}|^2 + |e_n|^{q-3}\mathbb E_n|\mathcal R_{n+1}|^3 + \mathbb E_n|\mathcal R_{n+1}|^q \leq
				\tilde{\eta} h|e_n|^q + C_{\tilde{\eta}} h^{q+1}(1+k_n^q).
			\end{flalign*} 
			Combining these estimates with \eqref{Taylor-q}, \eqref{zeta2}, \eqref{zetaa}, and \eqref{leading-drift-diff}, we obtain, on the event $A_n$,
			\begin{align*}
				\mathbb E_n |e_n+\zeta_n+\mathcal R_{n+1}|^q \leq \bigl(1+\Lambda_0 h\rho(X_n)\bigr)|e_n|^q + Ch^{q+1}(1+k_n^q)
			\end{align*}
			for $\Lambda_0>0$. Substituting this estimate into \eqref{step-start} proves \eqref{one-step-stability}. This completes the proof. \hfill
		\end{proof}
	\end{lemma}
	
	Set $Z_n:=|e_n|^q\mathbb I_{\{n<\sigma_X,\ n\leq\sigma_E\}}$, $a_n
	:=\Lambda_0\rho(X_n)\mathbb I_{\{n<\sigma_X,\ n<\sigma_E\}}$, $\xi_n
	:=C(1+k_n^q)\mathbb I_{\{n<\sigma_X,\ n<\sigma_E\}}$, and define $W_0:=1$, $W_{n+1}:= \frac{W_n}{1+ha_n}$, and $ S_n:=W_nZ_n$. The weights $W_n$ remove the multiplicative factor $1+ha_n$
	from the error recurrence. The following lemma uses the
	discrete Lenglart inequality to obtain a maximal estimate for $S_n$.
	\begin{lemma}\label{Lm4.9}
		Suppose the assumptions of Lemma \ref{Lm4.8} hold. 
		Then
		\begin{equation}\label{Zn-recurrence}
			\mathbb E_n [Z_{n+1}]\leq(1+ha_n)Z_n+h^{q+1}\xi_n.
		\end{equation}
		Moreover, for every $0<s<q$,
		\begin{equation}\label{weighted-max}
			\mathbb E \Big[ \Big( \max_{0\leq n\leq N}S_n  \Big)^{\frac{s}{q}} \Big]
			\leq Ch^s.
		\end{equation}
		\begin{proof}
			By Lemma \ref{Lm4.8}, we have 
			\begin{flalign*}
				\mathbb E_n\left[|e_{n+1}|^q\mathbb I_{\{n+1<\sigma_X,\ n+1\leq\sigma_E\}}\right] \leq\mathbb I_{\{n<\sigma_X,\ n<\sigma_E\}}\left[\bigl(1+\Lambda_0h\rho(X_n)\bigr)|e_n|^q
				+Ch^{q+1}(1+k_n^q)\right].
			\end{flalign*}
			Since $Z_{n+1} = |e_{n+1}|^q \mathbb I_{\{n+1<\sigma_X,\ n+1\leq\sigma_E\}}$ and $\mathbb I_{\{n<\sigma_X,\ n<\sigma_E\}}|e_n|^q \leq Z_n$, we have 
			\begin{flalign*}
				\mathbb E_n Z_{n+1} \leq (1+ha_n)Z_n+h^{q+1}\xi_n.
			\end{flalign*}
			Since $a_n\geq 0$, $ 0<W_{n+1}\leq W_n\leq1$.
			Moreover, $a_n$, $\xi_n$, and $W_{n+1}$ are $\mathcal F_{t_n}$-measurable. Multiplying \eqref{Zn-recurrence} by $W_{n+1}$, we obtain
			\begin{flalign*}
				\mathbb E_n S_{n+1} = W_{n+1}\mathbb E_n Z_{n+1}  \leq
				W_{n+1}(1+ha_n)Z_n	+ h^{q+1}W_{n+1}\xi_n  =
				S_n+h^{q+1}W_{n+1}\xi_n.
			\end{flalign*}
			
			Define the increasing predictable process $
			\mathcal B_0:=0$, $\mathcal B_n := \sum_{k=0}^{n-1} h^{q+1}W_{k+1}\xi_k$, $n\geq1$.
			It follows that $\mathbb E_n(S_{n+1}-\mathcal B_{n+1}) \leq S_n-\mathcal B_n$
			and therefore $(S_n-\mathcal B_n)_{n=0}^N$ is a supermartingale.
			Since $e_0=0$, we have $S_0=0$. Thus the discrete Lenglart domination inequality implies that, for every $0<\vartheta<1$,
			\begin{equation}\label{Lenglart-appl}
				\mathbb E \Big[ \Big( \max_{0\leq n\leq N}S_n \Big)^\vartheta
				\Big] \leq C_\vartheta\mathbb E[\mathcal B_N^\vartheta].
			\end{equation}
			We take $\vartheta=\frac{s}{q}\in(0,1)$. Since $W_n\leq 1$,
			\begin{flalign*}
				\mathcal B_N \leq Ch^{q+1}\sum_{n=0}^{N-1}(1+k_n^q)=
				Ch^q\Big[ h\sum_{n=0}^{N-1}(1+k_n^q) \Big].
			\end{flalign*}
			Consequently,
			\begin{flalign*}
				\mathbb E[\mathcal B_N^{\frac{s}{q}}] \leq Ch^s \mathbb E\Big[
				\Big(h\sum_{n=0}^{N-1}(1+k_n^q)\Big)^{\frac sq}
				\Big].
			\end{flalign*}
			Since $\frac{s}{q}\in(0,1)$, the function $x\mapsto x^{\frac{s}{q}}$ is concave. Hence Jensen's inequality gives
			\begin{flalign*}
				\mathbb E[\mathcal B_N^{\frac{s}{q}}] \leq Ch^s \Big[
				h\sum_{n=0}^{N-1} \mathbb E(1+k_n^q) \Big]^{\frac{s}{q}}.
			\end{flalign*}
			By Lemma \ref{Lm4.6} and H\"older's inequality, we obtain $\sup_{h,n}\mathbb E[k_n^q] \leq \sup_{h,n}
			\bigl(\mathbb E[k_n^{2q}]\bigr)^{\frac 12} \leq C$. Since $Nh=T$, $h\sum_{n=0}^{N-1} \mathbb E(1+k_n^q) \leq C$.
			Thus $\mathbb E[\mathcal B_N^{\frac{s}{q}}] \leq Ch^s$.
			Substituting this estimate into \eqref{Lenglart-appl} yields
			\begin{flalign*}
				\mathbb E \Big[ \Big( \max_{0\leq n\leq N}S_n \Big)^{\frac{s}{q}}
				\Big] \leq Ch^s.    
			\end{flalign*}
			Thus the assertion \eqref{weighted-max} holds.  \hfill
		\end{proof}
	\end{lemma}

	To remove the weights $W_n$ from the estimate in Lemma~\ref{Lm4.9}, we use the relative neighborhoods $\mathcal T_\varepsilon$ to compare the stopped sum $h\sum_{n=0}^{N-1}a_n$ with $\int_0^T\rho(X_t)\,\mathrm dt$.
	The exponential integrability estimate in Lemma~\ref{Lm4.3} is then used to control the inverse weights. Set $R_\varepsilon=\frac{1+\varepsilon}{1-\varepsilon}$ and $C_{\rho,\varepsilon} =\max\{1,R_\varepsilon^{\kappa_0}, \max_{1\le j\le J}R_\varepsilon^{\kappa_j}\}$.
	\begin{lemma}\label{Lm4.10}
		Suppose that the assumptions of Lemmas \ref{Lm4.3} and \ref{Lm4.9} hold. Let $q>s>r>0$ and assume that there exists a constant $\tilde{C}>C_{\rho,\varepsilon}$ such that
		\begin{flalign}\label{chireq}
			\tilde\chi:= \frac{\tilde{C}\Lambda_0rs}{q(s-r)} < \mathscr X_\rho^\dagger.
		\end{flalign}
		Then, for all sufficiently small $h>0$, 
		\begin{flalign}\label{enr}
			\mathbb E\Big[ \max_{\substack{0\leq n\leq N \\ n<\sigma_X,\ n\leq\sigma_E}} \vert{}e_n\vert{}^r \Big] \leq Ch^r.
		\end{flalign}
		\begin{proof}
			Recall that $Z_n = \vert{}e_n\vert{}^q \mathbb I_{\{n<\sigma_X,\ n\leq\sigma_E\}}$ and $S_n=W_n Z_n$, where $W_0=1$, $W_{n+1} = \frac{W_n}{1+ha_n}$ and $a_n = \Lambda_0\rho(X_n) \mathbb I_{\{n<\sigma_X,\ n<\sigma_E\}}$. By Lemma \ref{Lm4.9}, we have
			\begin{flalign}\label{Snsq}
				\mathbb E\Big[ \Big( \max_{0\leq n\leq N}S_n \Big)^{\frac sq} \Big] \leq Ch^s.
			\end{flalign}
			Fix a constant
			$\tilde C$ satisfying \eqref{chireq}. If $1\le n<\sigma_X$, then $X_n,X_t\in \mathcal T_\varepsilon(X_{n-1})$ for every $t\in[t_{n-1},t_n]$. Hence, $R_\varepsilon^{-1}\phi_j(X_t) \le\phi_j(X_n)\le R_\varepsilon\phi_j(X_t)$, $j=1,2,\ldots,J$, and $|X_n-x^\circ| \le R_\varepsilon|X_t-x^\circ|+(R_\varepsilon-1)$. For $\kappa_0>0$, the inequality $\tilde C>R_\varepsilon^{\kappa_0}$ implies that
			\begin{flalign*}
				(R_\varepsilon u+R_\varepsilon-1)^{\kappa_0} \le\tilde C u^{\kappa_0}+C_{\varepsilon,\tilde C}
			\end{flalign*}
			for $u\ge0$. Combining these observations with the definition of $\rho$ gives $\rho(X_n)\le\tilde C\rho(X_t) +C_{\varepsilon,\tilde C}$, $t\in[t_{n-1},t_n]$.
			Consequently,
			\begin{flalign*}
				h\rho(X_n) \le\tilde C\int_{t_{n-1}}^{t_n}\rho(X_t)\,\mathrm{d}t +C_{\varepsilon,\tilde C}h.
			\end{flalign*}
			For $1\le k\le N$, positivity and the definition of $a_n$ yield 
			\begin{align*}
				h\sum_{n=0}^{k-1}a_n
				\le h\Lambda_0\rho(X_0)
				+\tilde C\Lambda_0
				\sum_{n=1}^{k-1}\mathbb I_{\{n<\sigma_X,\,n<\sigma_E\}}
				\int_{t_{n-1}}^{t_n}\rho(X_t)\,\mathrm{d}t
				+\Lambda_0 C_{\varepsilon,\tilde C}T
				\le C+\tilde C\Lambda_0\int_0^T\rho(X_t)\,\mathrm{d}t.
			\end{align*}
			Here $X_0=x_0$ is deterministic and $h\le h_0$. Since $W_0^{-1}=1$ and $W_n^{-1}=\prod_{k=0}^{n-1}(1+ha_k)
			\le\exp\left(h\sum_{k=0}^{n-1}a_k\right)$, we obtain
			\begin{equation*}
				\max_{0\le n\le N}W_n^{-1}
				\le C\exp\left(\tilde C\Lambda_0
				\int_0^T\rho(X_t)\,\mathrm{d}t\right).
			\end{equation*}
			On $\{n<\sigma_X,\,n\le\sigma_E\}$, $\vert{}e_n\vert{}^r = Z_n^{\frac rq} = S_n^{\frac rq}W_n^{-\frac rq}$. Therefore, by H\"older's inequality, 
			\begin{flalign*}
				\mathbb E\Big[ \max_{\substack{0\leq n\leq N \\ n<\sigma_X,\ n\leq\sigma_E}} \vert{}e_n\vert{}^r \Big] \leq C \Big(\mathbb E \Big( \max_{0\leq n\leq N}S_n \Big)^{\frac sq}  \Big)^{\frac rs} \Big( \mathbb E \exp\Big( \frac{\tilde{C}\Lambda_0rs}{q(s-r)} \int_0^T\rho(X_t)\,\mathrm dt \Big) \Big)^{\frac{s-r}{s}}.
			\end{flalign*}
			Since $\tilde\chi =\frac{\tilde{C}\Lambda_0rs}{q(s-r)} <\mathscr X_\rho^\dagger$, Lemma \ref{Lm4.3} implies that the second factor is bounded. Combining this with \eqref{Snsq} yields
			\begin{flalign*}
				\mathbb E\Big[ \max_{\substack{0\leq n\leq N \\ n<\sigma_X,\ n\leq\sigma_E}} \vert{}e_n\vert{}^r \Big] \leq Ch^r.
			\end{flalign*}
			Thus the assertion \eqref{enr} holds.\hfill
		\end{proof}
	\end{lemma}

	\subsection{Global strong convergence}
	To remove the localization, we estimate $\mathbb P(\sigma_X\le N)$ and
	$\mathbb P(\sigma_E<\sigma_X)$. These estimates follow from the moment and boundary inverse-moment bounds in Lemma~\ref{Lm4.1}, the increment estimate in Lemma~\ref{Lm4.5} and the localized error estimate in Lemma~\ref{Lm4.10}.
	\begin{lemma}\label{Lm4.11}
		Suppose that the assumptions of Lemmas \ref{Lm4.1}, \ref{Lm4.5} and \ref{Lm4.10} hold. For each $j=1,2,\ldots,J$, choose $0<Q_{j,b}<q_j^\dagger$
		and define $Q_b:=\min_{1\leq j\leq J}Q_{j,b}$. If $\mathcal D$ is unbounded, choose $P_b$ such that $0<P_b<p_0$. Define
		\begin{flalign*}
			\Gamma_X := 
			\begin{cases}
				\min\{ \alpha Q_b, \, \beta P_b, \, \bar{q}(\frac{1}{2}-\alpha)-1 \}, & \mathcal D \text{ unbounded}, \\
				\min\{ \alpha Q_b, \, \bar{q}(\frac{1}{2}-\alpha)-1 \}, & \mathcal D \text{ bounded}\,, 
			\end{cases}
		\end{flalign*}
		and set $ \Gamma := \min\left\{ \Gamma_X, r(1-\alpha) \right\}$. Assume further that $
		\bar q(\frac12-\alpha)>1$. Then, for all sufficiently small $h>0$, the following estimates hold. 
		\begin{enumerate}
			\item[(1)]  $\mathbb P(\sigma_X\le N)\le Ch^{\Gamma_X}$.
			\item[(2)] $\mathbb P(\sigma_E<\sigma_X)
			\le Ch^{r(1-\alpha)}$.
			\item[(3)] For the event $G_h:=\{\sigma_X=N+1,\ \sigma_E=N+1\}$, we have
			\begin{equation}\label{bad-event-prob}
				\mathbb P(G_h^c)\le Ch^\Gamma.
			\end{equation}
		\end{enumerate}
		\begin{proof}
			We first estimate the probability that the exact solution leaves $K_h$ or the relative tube used in the definition of $\sigma_X$. For each $j=1,2,\ldots,J$, Lemma \ref{Lm4.1} and Markov's inequality give
			\begin{flalign}\label{boundry}
				\mathbb P\Big( \inf_{0\leq t\leq T}\phi_j(X_t)<8h^\alpha \Big) \leq C h^{\alpha Q_{j,b}} \mathbb E\Big[ \sup_{0\leq t\leq T} \phi_j(X_t)^{-Q_{j,b}} \Big] \leq Ch^{\alpha Q_{j,b}}.
			\end{flalign}
			Hence,
			\begin{flalign}\label{Pmin}
				\mathbb P\Big( \min_{1\leq j\leq J} \inf_{0\leq t\leq T}\phi_j(X_t) <8h^\alpha \Big) \leq Ch^{\alpha Q_b}.
			\end{flalign}
			
			If $\mathcal D$ is unbounded, then since $P_b<p_0$, Lemma \ref{Lm4.1} implies $\mathbb E[ \sup_{0\le t\le T} (1+|X_t-x^\circ|^2)^{\frac{P_b}{2}} ]\le C$. Therefore, 
			\begin{flalign}\label{Psup}
				\mathbb P\Big( \sup_{0\leq t\leq T}\vert{}X_t-x^{\circ}\vert{} > \frac18h^{-\beta} \Big) \leq Ch^{\beta P_b}.
			\end{flalign}
			It remains to control the event that the exact trajectory leaves the relative tube during one time step. Suppose that $X_n\in K_h$ and $X_t\notin\mathcal T_\varepsilon(X_n)$ for some $t\in[t_n,t_{n+1}]$. Then either $|X_t-X_n|>\varepsilon(1+|X_n-x^{\circ}|)$ or, for some $j$, $|\phi_j(X_t)-\phi_j(X_n)|> \varepsilon\phi_j(X_n)$. Since $\phi_j(X_n)\geq8h^\alpha$ and $|\phi_j(X_t)-\phi_j(X_n)| \leq |X_t-X_n|$, there exists $c_\varepsilon>0$ such that 
			\begin{flalign*}
				\sup_{t_n\leq t\leq t_{n+1}}\vert{}X_t-X_n\vert{} \geq c_\varepsilon h^\alpha.
			\end{flalign*}
			Therefore, by Markov's inequality and Lemma \ref{Lm4.5} with $a=\bar q$,
			\begin{flalign*}
				\mathbb P\left( X_n\in K_h,\,\exists t\in[t_n,t_{n+1}]:\, X_t\notin\mathcal T_\varepsilon(X_n) \right) \leq Ch^{-\alpha\bar q} \mathbb E\Big[ \sup_{t_n\leq t\leq t_{n+1}} |X_t-X_n|^{\bar q} \Big] \leq Ch^{\bar q(\frac 12-\alpha)}.
			\end{flalign*}
			Summing over $n=0,1,\ldots,N-1$ and using $N=T/h$, we obtain
			\begin{flalign}\label{tube}
				\mathbb{P}\left(\exists n:\,X_n\in K_h,\,\exists t\in[t_n,t_{n+1}]:\, X_t\notin\mathcal T_\varepsilon(X_n)\right)\leq Ch^{\bar q(\frac 12-\alpha)-1}.
			\end{flalign}
			
			In the unbounded case, combining this, \eqref{Pmin}, and \eqref{Psup} yields $\mathbb P(\sigma_X\leq N) \leq Ch^{\Gamma_X}$. If $\mathcal D$ is bounded, combining \eqref{Pmin} and \eqref{tube} yields
			\begin{flalign*}
				\mathbb P(\sigma_X\le N) \le Ch^{\Gamma_X},
			\end{flalign*}
			where
			\begin{flalign*}
				\Gamma_X = \min\{ \alpha Q_b,\, \bar q(\frac12-\alpha)-1 \}.
			\end{flalign*}
			On the other hand, on the event $\{\sigma_E<\sigma_X\}$, $|e_{\sigma_E}|>h^\alpha$. Hence, by Lemma \ref{Lm4.10},
			\begin{flalign*}
				\mathbb P(\sigma_E<\sigma_X) \leq h^{-\alpha r} \mathbb E\Big[ \max_{\substack{0\leq n\leq N \\ n<\sigma_X,\ n\leq\sigma_E}} \vert{}e_n\vert{}^r \Big] \leq Ch^{r(1-\alpha)}.
			\end{flalign*}
			Finally, $G_h^c \subset \{\sigma_X\leq N\} \cup \{\sigma_E<\sigma_X\}$ and therefore $\mathbb P(G_h^c) \leq C( h^{\Gamma_X} + h^{r(1-\alpha)} ) \leq Ch^\Gamma$. This completes the proof.  \hfill
		\end{proof}
	\end{lemma}

	We combine the localized error estimate in Lemma~\ref{Lm4.10} with the probability bounds in Lemma~\ref{Lm4.11} to obtain a global error bound. For unbounded domains, we choose $p<P_b<p_0$ and apply H\"older's inequality together with the moment bound in Lemma~\ref{Lm4.1} to control the contribution of the exact solution on the complementary event. The resulting parameter conditions are stated in the following theorem.
	\begin{theorem}\label{main-thm}
		Suppose that Assumptions \ref{ass0}--\ref{ass2} and \ref{ass4}--\ref{ass5} hold and that the commutative noise condition is satisfied. Assume further that the conditions of Lemmas \ref{Lm4.3} and \ref{Lm4.4} hold. Let $p>1$, let $\bar{q}> 4$ be as in Assumption~\ref{ass5} and choose $q\in[4,\bar{q})$. Choose $s,r$ such that $q>s>r>\max\{p,2\}$. Assume further that $\delta_S<\frac12$ and $\frac{C_{\rho,\varepsilon}\Lambda_0rs}{q(s-r)}
		<\mathscr X_\rho^\dagger$. Let $\Gamma$ be defined as in Lemma \ref{Lm4.11}. Assume the following conditions hold.
		\begin{enumerate}
			\item[(1)] If $\mathcal D$ is bounded, then assume that $\Gamma\ge p$.
			\item[(2)] If $\mathcal D$ is unbounded, then choose $P_b$
			such that $p<P_b<p_0$ and assume that
			\[
			\Gamma\Big(1-\frac{p}{P_b}\Big)\ge p
			\quad\text{and}\quad
			\Gamma-\beta p\ge p.
			\]
		\end{enumerate}
		Then there exists a constant $h_*\in(0,h_0]$ such that for all step sizes $h$ satisfying $0<h\le h_*$,
		\begin{flalign*}
			\mathbb E\Big[ \max_{0\le n\le N}|X(t_n)-Y_n|^p
			\Big] \le Ch^p.
		\end{flalign*}
		\begin{proof}
			Recall that $e_n=X_n-Y_n$. We decompose
			\begin{flalign*}
				\mathbb E\Big[\max_{0\leq n\leq N}|e_n|^p\Big]=
				\mathbb E\Big[\max_{0\leq n\leq N}|e_n|^p \mathbb I_{G_h}
				\Big]+\mathbb E\Big[\max_{0\leq n\leq N}|e_n|^p\mathbb I_{G_h^c}\Big].
			\end{flalign*}
			On the event $G_h$, we have $n<\sigma_X,\, n\leq\sigma_E,\, 0\leq n\leq N$. Hence, since $p<r$, Lemma \ref{Lm4.10} and H\"older's inequality yield
			\begin{equation}\label{error-Gh}
				\mathbb E\Big[\max_{0\leq n\leq N}|e_n|^p\mathbb I_{G_h}\Big]
				\leq\mathbb E\Big[\max_{\substack{0\leq n\leq N \\ n<\sigma_X,\ n\leq\sigma_E}}|e_n|^p\Big] \leq\Big(\mathbb E\Big[\max_{\substack{0\leq n\leq N \\ n<\sigma_X,\ n\leq\sigma_E}}|e_n|^r\Big]\Big)^{\frac pr}\leq	Ch^p.
			\end{equation}
			
			We next estimate the contribution on $G_h^c$.
			Suppose first that $\mathcal D$ is unbounded. Since $Y_n\in\mathcal D_h$,
			$|Y_n-x^{\circ}|\leq h^{-\beta}$. Therefore, $ |e_n| \leq |X_n-x^{\circ}|+|Y_n-x^{\circ}|
			\leq |X_n-x^{\circ}|+h^{-\beta}$ and hence $|e_n|^p \leq
			C(1+|X_n-x^{\circ}|)^p + Ch^{-\beta p}$. Consequently,
			\begin{equation*}
				\mathbb E\Big[ \max_{0\leq n\leq N}|e_n|^p \mathbb I_{G_h^c} \Big] \leq C\mathbb E\Big[ \sup_{0\leq t\leq T} (1+|X_t-x^{\circ}|)^p \mathbb I_{G_h^c} 	\Big] + Ch^{-\beta p}\mathbb P(G_h^c).
			\end{equation*}
			By H\"older's inequality,
			\begin{flalign*}
				\mathbb E\Big[\sup_{0\leq t\leq T}(1+|X_t-x^{\circ}|)^p
				\mathbb I_{G_h^c}\Big] \leq\Big(\mathbb E\Big[
				\sup_{0\leq t\leq T}(1+|X_t-x^{\circ}|)^{P_b}\Big]\Big)^{\frac{p}{P_b}}
				\mathbb P(G_h^c)^{1-\frac{p}{P_b}}.
			\end{flalign*}
			Since $P_b<p_0$, Lemma \ref{Lm4.1} implies that the first factor is bounded. Using Lemma \ref{Lm4.11}, $\mathbb P(G_h^c)\leq Ch^\Gamma$,
			we obtain
			\begin{equation*}
				\mathbb E\Big[\sup_{0\leq t\leq T}(1+|X_t-x^{\circ}|)^p
				\mathbb I_{G_h^c}\Big]\leq Ch^{\Gamma(1-\frac{p}{P_b})}.
			\end{equation*}
			Moreover,
			$  h^{-\beta p}\mathbb P(G_h^c) \leq Ch^{\Gamma-\beta p}$.
			Thus, by $\Gamma\Big(1-\frac{p}{P_b}\Big)\geq p$ and $\Gamma-\beta p\geq p$,
			\begin{flalign*}
				\mathbb E\Big[	\max_{0\leq n\leq N}|e_n|^p\mathbb I_{G_h^c}
				\Big]\leq Ch^p.
			\end{flalign*}
			
			If $\mathcal D$ is bounded, then both $X_n$ and $Y_n$ remain in the bounded set $\overline{\mathcal D}$, and hence $\max_{0\leq n\leq N}|e_n|^p\leq C$.
			Therefore,
			\begin{flalign*}
				\mathbb E\Big[ \max_{0\leq n\leq N}|e_n|^p \mathbb I_{G_h^c}
				\Big] \leq C\mathbb P(G_h^c) 	\leq Ch^\Gamma \leq Ch^p.
			\end{flalign*}
			Combining this estimate with \eqref{error-Gh} proves
			\begin{flalign*}
				\mathbb E\Big[\max_{0\leq n\leq N}|X_n-Y_n|^p
				\Big]\leq Ch^p.
			\end{flalign*}
			This completes the proof. \hfill
		\end{proof}
	\end{theorem}

	\section{Application to representative models}
	\label{sec4}
	In this section, we apply the convergence conclusions of Section \ref{sec3} to representative domain-constrained SDEs. For each model, we proceed in the following three steps:
	\begingroup
	\setlength{\leftmargini}{2em}
	\begin{enumerate}
		\renewcommand{\theenumi}{\roman{enumi}}
		\renewcommand{\labelenumi}{(\theenumi)}
		\item Verify Assumptions~\ref{ass0} and \ref{ass2} together with Assumption~\ref{ass1} for unbounded domains, thereby establishing the existence, uniqueness, and domain invariance of the exact solution;
		\item Verify Assumptions \ref{ass4}  and \ref{ass5}, along with the commutativity condition, and specify the compatibility conditions required in Theorem \ref{main-thm};
		\item Apply Theorem~\ref{main-thm} to obtain the strong order-one
		convergence of the projected Milstein approximation.
	\end{enumerate}
	\endgroup
	Throughout this section, $e_i$ denotes the $i$-th standard
	basis vector in the corresponding Euclidean space.
	
	\subsection{Stochastic Susceptible–Infected–Susceptible (SIS) epidemic model}
	
	Consider the stochastic SIS model \cite{GGHM2011}
	\begin{equation}\label{eq:SIS}
		\mathrm d I_t = I_t\bigl(\beta_e N_{\rm p}-\mu_e-\gamma_e-\beta_e I_t\bigr)\,\mathrm dt + \delta_e I_t(N_{\rm p}-I_t)\,\mathrm dB_t,
		\quad I_0\in(0,N_{\rm p}),
	\end{equation}
	where $\beta_e,\mu_e,\gamma_e,\delta_e>0$. Set $f(x)=x(a_0-\beta_e x)$ and $g(x)=\delta_e x(N_{\rm p}-x)$ with $a_0=\beta_e N_{\rm p}-\mu_e-\gamma_e$.
	
	(\romannumeral1) Take $\mathcal D=(0,N_{\rm p})$ with faces $\phi_1(x)=x$ and $\phi_2(x)=N_{\rm p}-x$. Clearly, $f,\,g\in C^{\infty}(\mathcal D)$, hence Assumption \ref{ass0} holds. For $\phi_1(x)=x$, we have $\nu_1=1$. For every
	$q_1>0$,
	\begin{align*}
		&-\phi_1(x)\nu_1 f(x) +\frac{q_1+1}{2}|\nu_1g(x)|^2 =
		-x^2(a_0-\beta_e x) +\frac{q_1+1}{2}\delta_e^2x^2(N_{\rm p}-x)^2 \\
		=& x^2\Big[ -a_0+\beta_e x +\frac{q_1+1}{2}\delta_e^2(N_{\rm p}-x)^2
		\Big] \le Cx^2.
	\end{align*}
	Hence Assumption~\ref{ass2} holds at $\phi_1(x)=x$. For $\phi_2(x)=N_{\rm p}-x$, we have $\nu_2=-1$. Then
	\begin{align*}
		&-\phi_2(x)\nu_2 f(x) +\frac{q_2+1}{2}|\nu_2g(x)|^2 =
		(N_{\rm p}-x)x(a_0-\beta_e x) +\frac{q_2+1}{2}\delta_e^2x^2(N_{\rm p}-x)^2 \\
		=& (N_{\rm p}-x)^2 \Big[ \frac{x(a_0-\beta_e x)}{N_{\rm p}-x} +\frac{q_2+1}{2}\delta_e^2x^2 \Big].
	\end{align*}
	Since $a_0-\beta_e N_{\rm p}=-(\mu_e+\gamma_e)<0$, we have $\frac{x(a_0-\beta_e x)}{N_{\rm p}-x}\longrightarrow-\infty$ as $x\uparrow N_{\rm p}$. Hence $\frac{x(a_0-\beta_e x)}{N_{\rm p}-x} +\frac{q_2+1}{2}\delta_e^2x^2\leq C$ and we obtain
	\begin{flalign*}
		-\phi_2(x)\nu_2 f(x) +\frac{q_2+1}{2}|\nu_2g(x)|^2\leq C\phi^2_2(x).
	\end{flalign*}
	Thus Assumption~\ref{ass2} holds for every finite $q_1,q_2>0$. Since $\mathcal D$ is bounded, $q_j^\dagger=q_j$, and consequently the boundary inverse-moment orders can be chosen arbitrarily large. Lemma~\ref{Lm4.1} therefore gives
	\begin{flalign*}
		\mathbb P(I_t\in(0,N_{\rm p}),\ 0\le t\le T)=1. 
	\end{flalign*}
	
	(\romannumeral2) Since $\mathcal D$ is bounded and $f$ and $g$ are polynomials, all the coefficients and derivatives appearing in Assumption~\ref{ass4} are bounded on $[0,N_{\rm p}]$. Hence we may choose $\mu_I=\mu_R=\mu_S=0$ and $\eta_I=\eta_R=\eta_S=0$. In particular, $\delta_S=0$. Thus, for every $\bar q>4$ and every $v\in\mathbb R$,
	\begin{flalign*}
		v^2f'(x) +\frac{\bar q-1}{2}v^2|g'(x)|^2\le C|v|^2. 
	\end{flalign*}
	Therefore Assumption~\ref{ass5} holds with $\rho(x)\equiv1$ and $\mathscr X_\rho^\dagger=\infty$. Since the equation is driven by a single Brownian motion, the commutativity condition is automatically satisfied. For any fixed $p>1$, choose $0<\alpha<\frac 12$ and then choose $r>\max\left\{p,2,\frac{p}{1-\alpha}\right\}$, $q>s>r$, and $q\ge4$. Since the boundary inverse-moment orders are arbitrary, we may choose $Q_b>\frac{p}{\alpha}$, and then choose $\bar q>q$ sufficiently large such that $\bar q\left(\frac12-\alpha\right)-1>p$. Consequently,
	\begin{flalign*}
		\Gamma = \min\Big\{ \alpha Q_b,\, \bar q\Big(\frac12-\alpha\Big)-1,\, r(1-\alpha) \Big\} > p. 
	\end{flalign*}
	Thus all the conditions of Theorem~\ref{main-thm} are fulfilled.
	
	(\romannumeral3) For every $p\ge1$ and all sufficiently small step sizes, the projected Milstein approximation satisfies
	\begin{flalign*}
		\mathbb E\left[ \max_{0\le n\le N}|I(t_n)-Y_n|^p
		\right] \le Ch^p. 
	\end{flalign*}
	For every $p>1$, the assertion follows directly from
	Theorem~\ref{main-thm}. For $p=1$, we have
	\begin{flalign*}
		\mathbb E\left[ \max_{0\le n\le N}|I(t_n)-Y_n|
		\right] \le \left( \mathbb E\left[ \max_{0\le n\le N}|I(t_n)-Y_n|^{p_1} \right] \right)^{\frac{1}{p_1}} \le Ch,
	\end{flalign*}
	where $p_1>1$. Hence the projected Milstein scheme converges strongly with order one for the stochastic SIS model.
	
	\subsection{Stochastic Susceptible--Infected--Recovered--Susceptible (SIRS) epidemic model with vaccination}
	
	We consider the following stochastic SIRS epidemic model with vaccination (see \cite{ZFWQ2021} for a related model)
	\begin{equation}\label{eq:SIRS}
		\left\{
		\begin{aligned}
			\mathrm dS_t
			&=\bigl[\Lambda-\beta_s S_tI_t-(\mu_s+\nu)S_t
			+\omega R_t\bigr]\mathrm dt
			-\sigma_s S_tI_t\,\mathrm dB_t,\\
			\mathrm dI_t
			&=\bigl[\beta_s S_tI_t-(\mu_s+\gamma_s)I_t\bigr]\mathrm dt
			+\sigma_s S_tI_t\,\mathrm dB_t,\\
			\mathrm dR_t
			&=\bigl[\gamma_s I_t+\nu S_t-(\mu_s+\omega)R_t\bigr]
			\mathrm dt,
		\end{aligned}
		\right.
	\end{equation}
	where $S_0>0$, $I_0>0$, $R_0>0$, and  $\Lambda,\beta_s,\gamma_s,\nu,\omega,\sigma_s\ge0$, $\mu_s>0$. Writing $X=(S,I,R)^\top$, the coefficients are
	\begin{flalign*}
		f(X) = \begin{pmatrix} \Lambda-\beta_s SI-(\mu_s+\nu)S+\omega R \\ \beta_s SI-(\mu_s+\gamma_s)I \\ \gamma_s I+\nu S-(\mu_s+\omega)R 
		\end{pmatrix} \quad \text{and}\quad g(X) = \sigma_s SI \begin{pmatrix} -1\\ 1\\ 0 \end{pmatrix}.
	\end{flalign*}
	
	(\romannumeral1) Define $N_t=S_t+I_t+R_t$. Adding the three equations \eqref{eq:SIRS} gives
	\begin{flalign*}
		\mathrm dN_t=(\Lambda-\mu_s N_t)\,\mathrm dt. 
	\end{flalign*}
	Hence
	\begin{flalign*}
		N_t = \frac{\Lambda}{\mu_s} + \Big( N_0-\frac{\Lambda}{\mu_s} \Big)e^{-\mu_s t}\quad \text{and}\quad N_0=S_0+I_0+R_0. 
	\end{flalign*}
	Choose $M > \max\{ N_0,\frac{\Lambda}{\mu_s}\}$ and define $ \mathcal D_M = \{ (S,I,R)\in\mathbb R^3: S>0,\ I>0,\ R>0,\ S+I+R<M\}$. The four boundary functions are $\phi_1(X)=S$, $\phi_2(X)=I$, $\phi_3(X)=R$, and $\phi_4(X)=\frac{M-S-I-R}{\sqrt3}$. For $\phi_1(X)$ and every $q_1>0$,
	\begin{align*}
		&-S f_1(X) +\frac{q_1+1}{2}|\nu_1^\top g(X)|^2=-S \bigl[ \Lambda-\beta_s SI-(\mu_s+\nu)S+\omega R \bigr] + \frac{q_1+1}{2}\sigma_s^2S^2I^2 \\
		= &S^2 \Big[ \beta_s I+\mu_s+\nu +\frac{q_1+1}{2}\sigma_s^2I^2
		\Big] -S(\Lambda+\omega R) \le CS^2.
	\end{align*}
	For $\phi_2(X)$,
	\begin{align*}
		-I f_2(X) +\frac{q_2+1}{2}|\nu_2^\top g(X)|^2 =
		I^2 \Big[ \mu_s+\gamma_s-\beta_s S +\frac{q_2+1}{2}\sigma_s^2S^2
		\Big] \le CI^2.
	\end{align*}
	For $\phi_3(X)$, since $\nu_3^\top g(X)=0$,
	\begin{flalign*}
		-Rf_3(X) = (\mu_s+\omega)R^2 -R(\gamma_s I+\nu S) \le(\mu_s+\omega)R^2. 
	\end{flalign*}
	Set $\tilde{N}=S+I+R$. For $\phi_4(X)$, $\nu_4=-\frac1{\sqrt3}(1,1,1)^\top$. Since
	\begin{flalign*}
		\sum_{i=1}^3f_i(X)=\Lambda-\mu_s \tilde{N}\quad \text{and}\quad
		(1,1,1)g(X)=0,
	\end{flalign*}
	we have
	\begin{flalign*}
		\nu_4^\top g(X)=0\quad\text{and}\quad
		\nu_4^\top f(X) = -\frac{\Lambda-\mu_s \tilde{N}}{\sqrt3}. 
	\end{flalign*}
	Consequently,
	\begin{flalign*}
		-\phi_4(X)\nu_4^\top f(X)=\frac{(M-\tilde{N})(\Lambda-\mu_s \tilde{N})}{3}. 
	\end{flalign*}
	Moreover, $\Lambda-\mu_s \tilde{N} = \mu_s(M-\tilde{N})+(\Lambda-\mu_s M) \le\mu_s(M-\tilde{N})$ because $M>\frac\Lambda\mu_s $. Hence
	\begin{flalign*}
		-\phi_4(X)\nu_4^\top f(X)\le\frac{\mu_s}{3}(M-\tilde{N})^2
		=\mu_s\phi_4(X)^2. 
	\end{flalign*}
	Thus Assumption~\ref{ass2} holds for every finite boundary exponent. Since $\mathcal D_M$ is bounded, Lemma~\ref{Lm4.1} implies
	\begin{flalign*}
		\mathbb P\bigl(X_t\in\mathcal D_M,\ 0\le t\le T
		\bigr)=1. 
	\end{flalign*}
	
	(\romannumeral2) The domain $\mathcal D_M$ is bounded and the coefficients are polynomial. Therefore all the coefficients and all the differential terms required in Assumption~\ref{ass4} are bounded, and we may take $\mu_I=\mu_R=\mu_S=0$ and
	$\eta_I=\eta_R=\eta_S=0$. In particular, $\delta_S=0$. Moreover, $Df$ and $Dg$ are bounded on $\mathcal D_M$. Hence for every $\bar q>4$,
	\begin{flalign*}
		\langle v,Df(x)v\rangle + \frac{\bar q-1}{2}|Dg(x)v|^2
		\le C|v|^2, \quad x\in\mathcal D_M. 
	\end{flalign*}
	Thus Assumption~\ref{ass5} holds with $\rho\equiv1$ and $\mathscr X_\rho^\dagger=\infty$. Since the equation is driven by a single Brownian motion, the commutativity condition is automatic. For any $p>1$, fix $0<\alpha<\frac 12$ and choose $r>\max\left\{p,2,\frac{p}{1-\alpha}\right\}$, $q>s>r$, $q\ge4$. Since all inverse-moment orders can be taken arbitrarily large, the auxiliary parameters can be chosen such that
	\begin{flalign*}
		\Gamma  = \min\Big\{ \alpha Q_b,\, \bar q\Big(\frac12-\alpha\Big)-1,\, r(1-\alpha) \Big\}>p. 
	\end{flalign*}
	Therefore all the assumptions of Theorem~\ref{main-thm} are satisfied.
	
	(\romannumeral3) Hence, for every $p>1$ and all sufficiently small step sizes, the projected Milstein approximation $Y_n$ satisfies
	\begin{flalign*}
		\mathbb E\left[ \max_{0\le n\le N} \left|X(t_n)-Y_n\right|^p \right] \le Ch^p,
	\end{flalign*}
	where $X=(S,I,R)^\top$. The case $p=1$ follows from any higher-moment estimate by H\"older's inequality. 
	
	\subsection{Competitive stochastic Lotka--Volterra (LV) system}
	
	Consider the following competitive stochastic LV system
	\begin{equation*}
		\left\{
		\begin{aligned}
			\mathrm{d}X_{1,t}=X_{1,t}(\tilde{a}_1-b_{11}X_{1,t}-b_{12}X_{2,t})\mathrm{d}t+\sigma_1X_{1,t}\mathrm{d}B_t^1, \\
			\mathrm{d}X_{2,t}=X_{2,t}(\tilde{a}_2-b_{21}X_{1,t}-b_{22}X_{2,t})\mathrm{d}t+\sigma_2X_{2,t}\mathrm{d}B_t^2, 
		\end{aligned}
		\right.
	\end{equation*}
	where $\tilde a_1,\tilde a_2\in\mathbb R$, $\sigma_1,\sigma_2\ge0$, $b_{11},b_{22}>0$ and $b_{12},b_{21}\ge0$, with $X_0=x_0\in(0,\infty)^2$.
	
	(\romannumeral1) Fix $x^\circ\in\mathcal D=(0,\infty)^2$. The coefficients are locally smooth, so Assumption~\ref{ass0} holds. For every $p_0>2$, using $b_{11},b_{22}>0$ and
	$b_{12},b_{21}\ge0$, we obtain
	\begin{flalign*}
		\langle x-x^\circ,f(x)\rangle+\frac{p_0-1}{2}\sum_{l=1}^2|g_l(x)|^2\le-b_{11}x_1^3-b_{22}x_2^3
		+C_{p_0}(1+|x|^2)\le-c_0|x|^3+C_{p_0}(1+|x|^2)
	\end{flalign*}
	for some $c_0>0$. Since $x^\circ$ is fixed, setting
	$r=|x-x^\circ|$, we have $|x|\ge \frac r2$ and $|x|\le Cr$ for all
	sufficiently large $r$. Hence
	\begin{flalign*}
		&-c_0|x|^3+C(1+|x|^2)
		\le
		K_1(1+r^2)-K_2(1+r^2)r
	\end{flalign*}
	for suitable $K_1,K_2>0$ and all sufficiently large $r$. Here $K_2>0$ can be chosen independently of $p_0$, while $K_1$ may depend on $p_0$. Therefore,
	\begin{flalign*}
		\langle x-x^\circ,f(x)\rangle
		+\frac{p_0-1}{2}\sum_{l=1}^2|g_l(x)|^2 \le
		K_1(1+|x-x^\circ|^2)
		-K_2(1+|x-x^\circ|^2)|x-x^\circ|.
	\end{flalign*}
	Thus Assumption~\ref{ass1} holds with $\gamma_1=1$. For the two coordinate faces, we take $\phi_i(x)=x_i$, $i=1,2$.
	For every $\widehat q_i>0$,
	\begin{align*}
		-x_i f_i(x) + \frac{\widehat q_i+1}{2}
		\sum_{l=1}^2 |e_i^\top g_l(x)|^2 = x_i^2
		\Big[ -\tilde{a}_i +b_{i1}x_1+b_{i2}x_2 +\frac{\widehat q_i+1}{2}\sigma_i^2 \Big]\le K_{3,i}x_i^2
		\left(1+|x-x^\circ|\right).
	\end{align*}
	Thus Assumption~\ref{ass2} holds with $K_{4,i}=0$ and $\gamma_{3,i}=1$. Since $\gamma_{3,i}=\gamma_1=1$, the inverse-moment exponent in Lemma~\ref{Lm4.1} is $q_i^\dagger = \frac{ p_0K_2\widehat q_i }{ p_0K_2+K_{3,i}\widehat q_i }$. Given any finite $Q_i>0$, first choose $\widehat q_i>Q_i$. Then $q_i^\dagger \rightarrow \widehat q_i$ as $p_0\to\infty$. Consequently, $p_0$ can be chosen sufficiently large so that $q_i^\dagger>Q_i$, $i=1,2$. Lemma~\ref{Lm4.1} then yields positivity of both components. 
	
	(\romannumeral2) The drift has quadratic growth, whereas the diffusion has linear growth: $|f(x)|\le C(1+|x|^2)$ and $\|g(x)\|_F\le C(1+|x|)$. Moreover, $\|Df(x)\|\le C(1+|x|)$ and $\|Dg_l(x)\|\le C$. Thus we have 
	\begin{flalign*}
		|L^0f(x)| + \Big(  \sum_l|L^lf(x)|^2 \Big)^{\frac 12}
		+ \Big( \sum_k|L^0g_k(x)|^2 \Big)^{\frac 12} \le C(1+|x|^3). 
	\end{flalign*}
	Furthermore,
	\begin{flalign*}
		H_{11}(x) =  \begin{pmatrix} \sigma_1^2x_1\\ 0
		\end{pmatrix}
		\quad \text{and }\quad
		H_{22}(x) = \begin{pmatrix} 0\\ \sigma_2^2x_2
		\end{pmatrix}, 
	\end{flalign*}
	while $H_{12}=H_{21}=0$. Hence $DH_{lk}$ is uniformly bounded. A convenient choice in Assumption~\ref{ass4} is therefore $\mu_I=2$, $\mu_R=3$, $\mu_S=1$, $\eta_I=\eta_R=\eta_S=0$. For Assumption~\ref{ass5}, since $Df$ has at most linear growth and $Dg_l$ is bounded,
	\begin{flalign*}
		\langle v,Df(x)v\rangle + \frac{\bar q-1}{2} \sum_{l=1}^2|Dg_l(x)v|^2 \le C\bigl(1+|x-x^\circ|\bigr)|v|^2. 
	\end{flalign*}
	Thus we may take $\rho(x) = 1+d_0|x-x^\circ| $ for $d_0>0$. Since $\kappa_0=\gamma_1=1$, Lemma~\ref{Lm4.3} gives $\mathscr X_\rho^\dagger = \frac{p_0K_2}{d_0}$. The diffusion vector fields commute. Indeed, $Dg_2(x)g_1(x)=0$ and $Dg_1(x)g_2(x)=0$, hence $H_{12}=H_{21}$. We next verify that the auxiliary parameters can be chosen consistently. Fix $p>1$, $0<\alpha<\frac 12$ and $0<\beta<\frac{1}{2}$. Since the moment and boundary inverse-moment orders can be taken sufficiently large, we may choose $q>s>r>\max\{p,2\}$ with $q\geq4$, and $\bar q>q$, together with $Q_b$ and $P_b$, sufficiently large so that
	\begin{flalign*}
		\Gamma \Big( 1-\frac{p}{P_b} \Big)>p\quad \text{and} \quad  \Gamma-\beta p>p,
	\end{flalign*}
	where $\Gamma= \min\{ \alpha Q_b,\,  \beta P_b,\, \bar q\left(\frac12-\alpha\right)-1,\, r(1-\alpha) \}$. Finally, $p_0$ can be chosen sufficiently large so that $p_0>\max\{m_I\mu_I,\,2q\mu_R,\,4q\mu_S,\,P_b\}$ and $\frac{C_{\rho,\varepsilon}\Lambda_0rs}{q(s-r)} <\mathscr X_\rho^\dagger  =\frac{p_0K_2}{d_0}$. Hence, all the conditions of Theorem~\ref{main-thm} can be fulfilled.
	
	(\romannumeral3) Hence, for every $p>1$ and all sufficiently small step sizes, the projected Milstein approximation $Y_n$ satisfies
	\begin{flalign*}
		\mathbb E\left[ \max_{0\le n\le N} \left|X(t_n)-Y_n\right|^p \right] \le Ch^p.
	\end{flalign*}
	The case $p=1$ follows from any higher-moment estimate by H\"older's inequality.

	\subsection{Two-island Wright--Fisher (WF) diffusion }
	
	Consider the following two-island WF diffusion \cite{SK1986,KE2008}
	\begin{equation*}
		\left\{
		\begin{aligned}
			\mathrm{d}X_t^1=\left[\tilde{u}_2-(\tilde{u}_1+\tilde{u}_2)X_t^1+\tilde{\vartheta}(X_t^2-X_t^1)\right]\mathrm{d}t
			+\sqrt{X_t^1(1-X_t^1)}\mathrm{d}B_t^1, \\
			\mathrm{d}X_t^2=\left[\tilde{u}_2-(\tilde{u}_1+\tilde{u}_2)X_t^2+\tilde{\vartheta}(X_t^1-X_t^2)\right]\mathrm{d}t
			+\sqrt{X_t^2(1-X_t^2)}\mathrm{d}B_t^2,
		\end{aligned}
		\right.
	\end{equation*}
	where $\tilde{\vartheta}\ge0$ and $B^1,B^2$ are independent Brownian motions. Assume $X_0=x_0\in(0,1)^2$ and take $\mathcal D=(0,1)^2$. 
	
	(\romannumeral1) For the left and right faces of coordinate $i$, define $\phi_{i,L}(x)=x_i$ and $\phi_{i,R}(x)=1-x_i$. For the left face,
	\begin{align*}
		-x_i f_i(x) + \frac{\widehat q_{i,L}+1}{2} x_i(1-x_i) =&
		- \Big( \tilde{u}_2-\frac{\widehat q_{i,L}+1}{2} \Big)x_i -\tilde{\vartheta} x_ix_j+\Big(\tilde{u}_1+\tilde{u}_2+\tilde{\vartheta}-\frac{\widehat q_{i,L}+1}{2}\Big)x_i^2\\
		\leq& - \Big( \tilde{u}_2-\frac{\widehat q_{i,L}+1}{2} \Big)x_i+\Big(\tilde{u}_1+\tilde{u}_2+\tilde{\vartheta}-\frac{\widehat q_{i,L}+1}{2}\Big)x_i^2
	\end{align*} 
	for $j\not= i$. Hence, if $0<\widehat q_{i,L}<2\tilde{u}_2-1$, Assumption~\ref{ass2} holds with $\gamma_{2,i,L}=1$ and $K_{4,i,L} = \tilde{u}_2-\frac{\widehat q_{i,L}+1}{2}>0$. Similarly, for the right face, if $0<\widehat q_{i,R}<2\tilde{u}_1-1$, then Assumption~\ref{ass2} holds with $\gamma_{2,i,R}=1$ and $ K_{4,i,R} = \tilde{u}_1-\frac{\widehat q_{i,R}+1}{2}>0$. Since $\mathcal D$ is bounded, $q_{i,L}^\dagger=\widehat q_{i,L}$ and  $q_{i,R}^\dagger=\widehat q_{i,R}$.
	In particular, if $\tilde{u}_1,\tilde{u}_2>\frac 12$, Lemma~\ref{Lm4.1} gives
	\begin{flalign*}
		\mathbb P\bigl(X_t\in(0,1)^2,\ 0\le t\le T\bigr)=1. 
	\end{flalign*}
	
	(\romannumeral2) The drift and diffusion are bounded on $\mathcal D$.
	Hence we may choose $\mu_I=0$ and $\eta_I=0$. For $g_i(x) = \sqrt{x_i(1-x_i)}\,e_i$, we have $Dg_i(x) = \frac{1-2x_i} {2\sqrt{x_i(1-x_i)}} e_i e_i^\top$. Thus $Dg_i$ has boundary singularities of order $\frac 12$. Moreover, $H_{ii}(x) = Dg_i(x)g_i(x) = \frac{1-2x_i}{2}e_i$ and $H_{ij}=0$, $i\ne j$. Hence the diffusion fields commute and $DH_{ij}$ is bounded. Since $f$ is affine, $L^0f$ and $L^lf$ are bounded on $\mathcal D$. Differentiation gives
	\begin{flalign*}
		|L^0g_i(x)| \le C\big(x_i^{-\frac 12}+(1-x_i)^{-\frac 12}\big),
		\quad i=1,2,\quad x\in\mathcal D.
	\end{flalign*}
	A convenient choice in Assumption~\ref{ass4} is therefore $\eta_I=0$ and, for $i=1,2$, $\eta_{R,i,L}= \eta_{R,i,R} = \frac12$ and $\eta_{S,i,L} = \eta_{S,i,R} = \frac12$. Lemma~\ref{Lm4.4} then requires $2q < \min_i \left\{ \widehat q_{i,L}, \widehat q_{i,R}
	\right\}$. Writing $u_*:=\min\{\tilde u_1,\tilde u_2\}$, the existence of some $q\ge4$ satisfying the preceding inequality requires $u_*>\frac92$. For Assumption~\ref{ass5}, observe that
	\begin{flalign*}
		|Dg_i(x)v|^2 = \frac{(1-2x_i)^2}{4x_i(1-x_i)} v_i^2
		\le \frac14 \Big( \frac1{x_i} + \frac1{1-x_i} \Big)v_i^2. 
	\end{flalign*}
	Since $Df$ is bounded, Assumption~\ref{ass5} holds with
	\begin{flalign*}
		\rho(x) = 1+ d_b\sum_{i=1}^2 \left( x_i^{-1}+(1-x_i)^{-1}
		\right)\quad \text{and}\quad
		d_b=\frac{\bar q-1}{8}. 
	\end{flalign*}
	For the selected boundary exponents, define $K_{4,i,L} = \tilde{u}_2-\frac{\widehat q_{i,L}+1}{2}$ and $K_{4,i,R} =\tilde{u}_1-\frac{\widehat q_{i,R}+1}{2}$.
	Then $\mathscr X_{i,L}^\dagger = \frac{ K_{4,i,L}\widehat q_{i,L}
	}{2} = \frac{ \widehat q_{i,L} (2\tilde{u}_2-1-\widehat q_{i,L})
	}{4}$ and $\mathscr X_{i,R}^\dagger =  \frac{ K_{4,i,R}\widehat q_{i,R} }{2}
	= \frac{ \widehat q_{i,R} (2\tilde{u}_1-1-\widehat q_{i,R}) }{4}$.
	Consequently,
	\begin{flalign*}
		\frac1{\mathscr X_\rho^\dagger} =
		d_b\sum_{i=1}^2 \Big( \frac1{\mathscr X_{i,L}^\dagger} + \frac1{\mathscr X_{i,R}^\dagger} \Big). 
	\end{flalign*}
	Furthermore, $\delta_S=\frac{\alpha}{2}<\frac12$. To apply Theorem~\ref{main-thm}, choose $\bar q>q>s>r>\max\{p,2\}$, $q\ge4$, $0<\alpha<\frac12$, together with boundary exponents satisfying $2q < \min_i \left\{ \widehat q_{i,L}, \widehat q_{i,R} \right\}$, $\widehat q_{i,L}<2\tilde{u}_2-1$ and $\widehat q_{i,R}<2\tilde{u}_1-1$. We further require $\frac{ C_{\rho,\varepsilon}\Lambda_0rs }{ q(s-r) } < \mathscr  X_\rho^\dagger$. Choose $0<Q_b< \min_i \left\{ \widehat q_{i,L}, \widehat q_{i,R} \right\}$. Since the domain is bounded, the globalization exponent is
	\begin{flalign*}
		\Gamma  = \min\Big\{ \alpha Q_b,\, \bar q\Big(\frac12-\alpha\Big)-1,\, r(1-\alpha) \Big\}. 
	\end{flalign*}
	Thus Theorem~\ref{main-thm} applies whenever $\Gamma\ge p$.

	(\romannumeral3) Let $u_*:=\min\{\tilde{u}_1,\tilde{u}_2\}>\frac92$. For every $1<p<u_*-\frac12$, suppose that the auxiliary parameters are chosen so that $\frac{C_{\rho,\varepsilon}\Lambda_0rs} {q(s-r)} <\mathscr X_\rho^\dagger$, $\Gamma\ge p$. Then, for each \(p\) satisfying the preceding moment and compatibility conditions, the projected Milstein approximation satisfies
	\begin{flalign*}
		\mathbb E\left[
		\max_{0\le n\le N}|X(t_n)-Y_n|^p
		\right]\le Ch^p 
	\end{flalign*}
	for sufficiently small \(h\). The case $p=1$ follows from a higher-moment estimate and H\"older's inequality.
	
	\subsection{Stochastic SIR model with population-dependent noise}
	
	We construct the following stochastic SIR model with
	population-dependent noise in It\^o form
	\begin{equation}\label{eq:sir_system}
		\left\{
		\begin{aligned}
			\mathrm{d}S_{t} &= S_{t} [ -a(\tilde{K}-N_{t}) + bI_{t} - \tilde{\delta}\frac{I_{t}}{N_{t}} ] \mathrm{d}t + \sigma S_{t}(\tilde{K}-N_{t}) \mathrm{d}B_{t}, \\
			\mathrm{d}I_{t} &= I_{t} [ -a(\tilde{K}-N_{t}) - bS_{t} - cR_{t} + \tilde{\delta}\frac{S_{t}}{N_{t}} ] \mathrm{d}t + \sigma I_{t}(\tilde{K}-N_{t}) \mathrm{d}B_{t}, \\
			\mathrm{d}R_{t} &= R_{t} [ cI_{t} - a(\tilde{K}-N_{t}) ] \mathrm{d}t + \sigma R_{t}(\tilde{K}-N_{t}) \mathrm{d}B_{t},
		\end{aligned}
		\right.
	\end{equation}
	where $N_{t} = S_{t} + I_{t} + R_{t}$, $\tilde{K} > 0$ and $a, b, c, \tilde{\delta},\sigma\ge 0$. Set $\mathcal{D}_{\tilde{K}} = \left\{ (S, I, R) : S, I, R > 0, S + I + R < \tilde{K} \right\}$. Assume that $X_0=(S_0,I_0,R_0)^\top\in\mathcal D_{\tilde{K}}$.

	(\romannumeral1) The coefficients are locally smooth on $\mathcal D_{\tilde{K}}$. For
	$x=(x_1,x_2,x_3)^\top=(S,I,R)^\top\in\mathcal D_{\tilde K}$, define $\mathcal N(x):=S+I+R$. Each drift component can be written as $f_i(x)=x_i\widetilde f_i(x)$, $i=1,2,3$, where $\widetilde f_i$ is uniformly bounded on $\mathcal D_{\tilde{K}}$, since $0<\mathcal N(x)<\tilde K$ and $0<\frac{x_i}{\mathcal N(x)}<1$. For the coordinate faces $\phi_i(x)=x_i$, we have $\nu_i=e_i$, $i=1,2,3$. Since $g(x)=\sigma(\tilde K-\mathcal N(x))x$, it follows that
	\[
	\nu_i^\top g(x)=\sigma x_i(\tilde K-\mathcal N(x)), \quad i=1,2,3.
	\]
	Hence, for every $q_i>0$,
	\begin{flalign*}
		-\phi_i(x)\nu_i^\top f(x)+\frac{q_i+1}{2}|\nu_i^\top g(x)|^2 =x_i^2\Big[-\widetilde f_i(x)+\frac{q_i+1}{2}\sigma^2(\tilde{K}-\mathcal N(x))^2\Big]\le C\phi_i(x)^2.
	\end{flalign*}
	Hence Assumption~\ref{ass2} holds on each coordinate face for every finite $q_i>0$. Define $\phi_4(x) = \frac{\tilde{K}-\mathcal N(x)}{\sqrt3}$, $
	\nu_4 = -\frac1{\sqrt3}(1,1,1)^\top$.
	Summing the three equations in \eqref{eq:sir_system} gives
	\begin{flalign*}
		\mathrm dN_t = -aN_t(\tilde{K}-N_t)\,\mathrm dt +
		\sigma N_t(\tilde{K}-N_t)\,\mathrm dB_t. 
	\end{flalign*}
	Hence $\sum_{i=1}^3f_i(x) = -a\mathcal{N}(x)(\tilde{K}-\mathcal N(x))$ and $(1,1,1)g(x)=\sigma \mathcal N(x)(\tilde K-\mathcal N(x))$. Therefore $\nu_4^\top f(x) = a\mathcal N(x)\phi_4(x)$ and $ \nu_4^\top g(x) = -\sigma \mathcal N(x)\phi_4(x)$. It follows that
	\begin{align*}
		-\phi_4(x)\nu_4^\top f(x) + \frac{q_4+1}{2} |\nu_4^\top g(x)|^2 = \Big[ -a\mathcal N(x)+ \frac{q_4+1}{2}\sigma^2\mathcal N(x)^2
		\Big]\phi_4(x)^2 \le C\phi_4(x)^2.
	\end{align*}
	Thus Assumption~\ref{ass2} holds for arbitrary finite boundary exponents. Since $\mathcal D_{\tilde{K}}$ is bounded, Lemma~\ref{Lm4.1} yields
	\begin{flalign*}
		\mathbb P
		\bigl(
		X_t\in\mathcal D_{\tilde{K}},\ 0\le t\le T
		\bigr)=1. 
	\end{flalign*}
	
	(\romannumeral2) The only term requiring special attention near $\mathcal{N}(x)=0$ is $\mathcal H(S,I,R):=\frac{SI}{\mathcal N(x)}$ since its higher-order derivatives may become singular as $\mathcal N(x)\downarrow0$.
	Its first derivatives are $\frac{\partial \mathcal H}{\partial S} = \frac{I(I+R)}{\mathcal N(x)^2}$, $\frac{\partial \mathcal H}{\partial I} = \frac{S(S+R)}{\mathcal N(x)^2}$, and $\frac{\partial \mathcal H}{\partial R} = -\frac{SI}{\mathcal N(x)^2}$. Hence $|D\mathcal H(x)|\le C$, $x\in\mathcal D_{\tilde{K}}$. Although $D^2\mathcal H$ need not be uniformly bounded as $\mathcal{N}(x)\downarrow0$, the function $\mathcal H$ is homogeneous of degree one: $\mathcal H(\lambda x)=\lambda \mathcal H(x)$, $\lambda>0$. Euler's identity gives $D\mathcal H(x)x=\mathcal H(x)$. Differentiating this identity yields $D^2\mathcal H(x)[x,x]=0$. Since the diffusion vector is $g(x)
	= \sigma(\tilde{K}-\mathcal{N}(x))x$, we obtain 
	\begin{flalign*}
		D^2\mathcal H(x)[g(x),g(x)] = \sigma^2(\tilde{K}-\mathcal{N}(x))^2 D^2\mathcal H(x)[x,x] =0. 
	\end{flalign*}
	Consequently, the potentially singular second derivative does not contribute to the It\^o term $L^0f$. All the quantities required in Assumption~\ref{ass4} are therefore uniformly bounded on $\mathcal D_{\tilde{K}}$, and we may take $\mu_I=\mu_R=\mu_S=0$,
	$\eta_I=\eta_R=\eta_S=0$ and $\delta_S=0$. Moreover, $Df$ and $Dg$ are uniformly bounded. Thus
	\begin{flalign*}
		\langle v,Df(x)v\rangle+\frac{\bar q-1}{2}|Dg(x)v|^2
		\le C|v|^2, 
	\end{flalign*}
	and Assumption~\ref{ass5} holds with $\rho\equiv1$ and $\mathscr X_\rho^\dagger=\infty$.
	Since the system is driven by a single Brownian motion, the
	commutativity condition is automatic. For any $p>1$, fix $0<\alpha<\frac 12$ and choose $r>\max\left\{p,2,\frac{p}{1-\alpha}\right\}$, $q>s>r$, $q\ge4$. Since arbitrary finite boundary inverse-moment orders are available, $Q_b$ and $\bar q$ can be chosen sufficiently large so that
	\begin{flalign*}
		\Gamma = \min\Big\{ \alpha Q_b,\, \bar q\Big(\frac12-\alpha\Big)-1,\, r(1-\alpha) \Big\}>p. 
	\end{flalign*}
	Hence all the assumptions of Theorem~\ref{main-thm} are satisfied.
	
	(\romannumeral3) For every $p>1$ and all sufficiently small step sizes, the projected Milstein approximation $Y_{n}$ satisfies
	\begin{equation*}
		\mathbb{E}\left[ \max_{0 \le n \le N} |X(t_{n}) - Y_{n}|^{p} \right] \le C h^{p}.
	\end{equation*}
	The case $p = 1$ follows from any higher-moment estimate by H\"older's inequality.

	\section{Numerical experiments}\label{sec:numerics}
	\label{sec5}
	In this section, we present numerical experiments to illustrate the strong convergence order of the proposed projected Milstein scheme. We  consider the representative models studied in Section~\ref{sec4} and compare the proposed scheme with the first-order domain-preserving scheme of \cite{br26} on a common problem.
	
	Throughout the numerical experiments, we take $T=1$ and
	$M_{\rm MC}=2000$. The tested step sizes are $h=2^{-11},\,2^{-12},\,2^{-13},\,2^{-14},\,2^{-15}$, and the reference solution is computed by the same projected Milstein scheme with the step size $h_{\rm ref}=2^{-18}$. For each Monte Carlo sample path, all discretization levels and the reference solution are driven by the same Brownian trajectory. For the $\ell$-th sample path, let $Y_n^{h,\ell}$ denote the projected Milstein approximation with step size $h$, and let $X_{\rm ref}^{\ell}(t_n)$ denote the corresponding reference value. Writing $h=Rh_{\rm ref}$, we have $X_{\rm ref}^{\ell}(t_n)=X_{nR}^{\rm ref,\ell}$. A reference line of slope one is included in each log--log plot for comparison with the convergence order. We measure the strong error by
	\begin{flalign*}
		{\rm Err}(h)
		:=
		\Bigg(
		\frac{1}{M_{\rm MC}}
		\sum_{\ell=1}^{M_{\rm MC}}
		\max_{0\le n\le N}
		\left|
		Y_n^{h,\ell}
		-
		X_{\rm ref}^{\ell}(t_n)
		\right|^2
		\Bigg)^{\frac 12}.
	\end{flalign*}
	This quantity estimates the discrete-time $L^2$ maximum error
	using a fine-grid numerical reference solution. For the interior projection, different values of $\alpha$ are used in the numerical experiments, all satisfying $0<\alpha<\frac 12$. In the unbounded case, we additionally take $\beta=0.25$.
	
	The projections for the SIS, WF and hypercube comparison examples are computed by componentwise clipping.
	
	For the SIRS and SIR examples, let $z$ denote the Milstein predictor and set $\ell_h = h^\alpha$. We shift the predictor by defining $w_i = z_i - \ell_h$ for $i = 1, 2, 3$. In the shifted variables $v_i = y_i - \ell_h$, the feasible set becomes $\{v \in \mathbb{R}^3 : v_i \ge 0,\ i = 1, 2, 3, \quad \sum_{i=1}^3 v_i \le b_h\}$, where $b_h = M - (3 + \sqrt{3})\ell_h$ for SIRS and $b_h = \tilde K - (3 + \sqrt{3})\ell_h$ for SIR. We assume $h$ is sufficiently small so that $b_h > 0$. Let $(a)_+ = \max\{a, 0\}$. If $\sum_{i=1}^3 (w_i)_+ \le b_h$, we set $\tau = 0$. Otherwise, we sort the components of $w$ so that $w_{(1)} \ge w_{(2)} \ge w_{(3)}$ and compute
	\begin{flalign*}
		\tau = \max_{1 \le k \le 3} \frac{\sum_{j=1}^k w_{(j)} - b_h}{k}.
	\end{flalign*}
	In this case, $\tau > 0$ and $\sum_{i=1}^3 (w_i - \tau)_+ = b_h$. The projected point in the original coordinates is then given by $y_i = \ell_h+(w_i-\tau)_+, \quad i = 1, 2, 3$.
	
	For the LV example, set $\ell_h=h^\alpha$ and $R_h=h^{-\beta}$, and let $z$ be the Milstein predictor. Define  $y_i(\lambda)=\max\{\ell_h, x_i^\circ+\frac{z_i-x_i^\circ}{1+\lambda}\}$, $i=1,2$. If \(|y(0)-x^\circ|\le R_h\), we return \(y(0)\). Otherwise, we solve \(|y(\lambda)-x^\circ|=R_h\) by bisection, first doubling the upper bound for \(\lambda\) until the corresponding point lies inside the ball. Let \(y^{\mathrm{lo}}\) and \(y^{\mathrm{hi}}\) denote the points at the current lower and upper multiplier bounds. We stop when $|y^{\mathrm{hi}}-y^{\mathrm{lo}}|
	\le10^{-13}(1+|y^{\mathrm{hi}}|)$ and return the feasible point \(y^{\mathrm{hi}}\).
	
	\begin{example}{\textbf{Stochastic SIS model}}
		
		We first consider the stochastic SIS model~\eqref{eq:SIS}. The parameters are chosen as $N_{\rm p}=1$,
		$\beta_e=2$, $\mu_e=0.4$, $\gamma_e=0.3$,
		$\delta_e=0.35$ and $I_0=0.4$. The interval is $\mathcal D=(0,N_{\rm p})$, and the projection interval is 
		\begin{flalign*}
			&\mathcal D_h=[h^\alpha,N_{\rm p}-h^\alpha],
		\end{flalign*}
		with $\alpha=0.4$. Figure~\ref{SIS} shows that the strong convergence order of the projected Milstein is in agreement with the expected order. 
		
		\begin{figure}[htbp]
			\centering
			\includegraphics[width=0.8\textwidth]{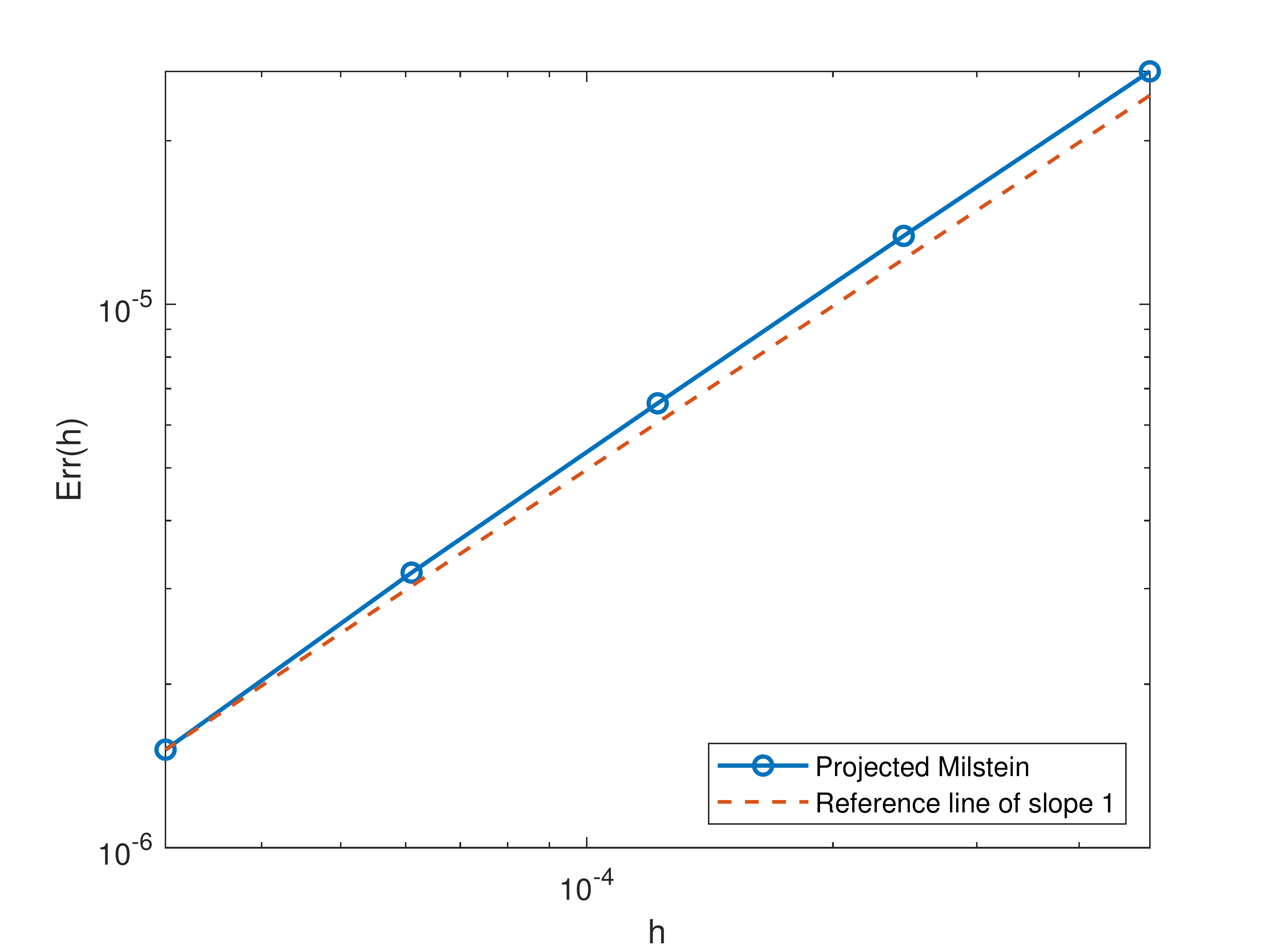}
			\caption{Stochastic SIS model}
			\label{SIS}
		\end{figure}
		
	\end{example}
	
	\begin{example}{\textbf{Stochastic SIRS epidemic model with vaccination}}
		
		We next consider the stochastic SIRS model~\eqref{eq:SIRS}. 
		We take $\alpha=0.4$, $\Lambda=0.2$, $\mu_s=0.5$,
		$\beta_s=1.2$, $\gamma_s=0.3$, $\nu=0.15$,
		$\omega=0.1$ and $\sigma_s=0.25$,
		with $X_0=(0.60,0.25,0.15)^\top$ and $M=2$.
		The projected set is
		\begin{flalign*}
			\mathcal D_h = \left\{ (S,I,R)^\top:
			S,I,R\ge h^\alpha,\quad
			S+I+R\le M-\sqrt{3}\,h^\alpha
			\right\}.
		\end{flalign*}
		Figure~\ref{SIR} shows that the strong convergence order of the projected Milstein is close to $1$. 
		
		\begin{figure}[htbp]
			\centering
			\includegraphics[width=0.8\textwidth]{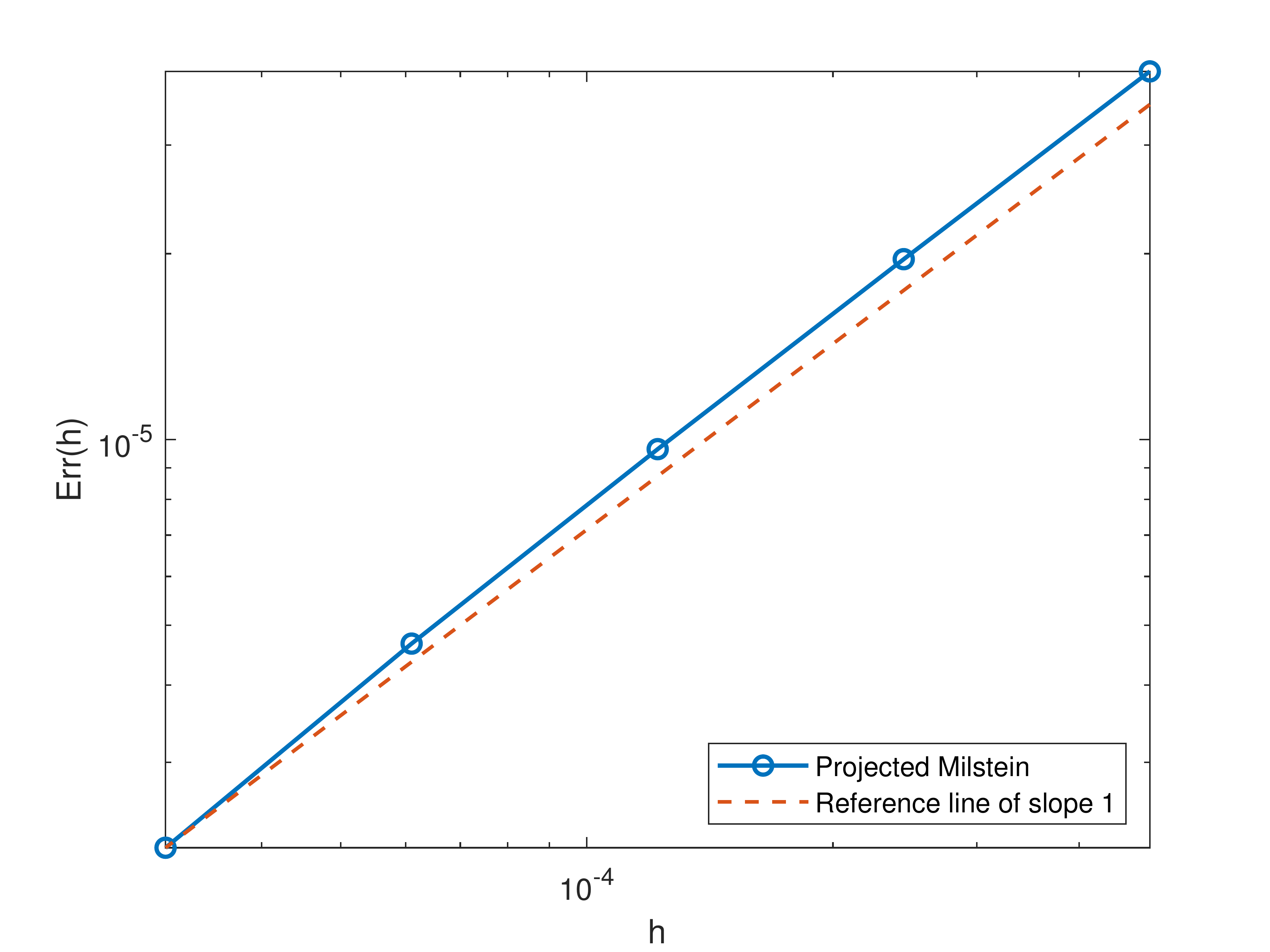}
			\caption{Stochastic SIRS epidemic model with vaccination}
			\label{SIR}
		\end{figure}
	\end{example}

	\begin{example}{\textbf{Competitive stochastic LV system}}
		
		We choose $\tilde a_1=2$, $\tilde a_2=1.6$, $b_{11}=b_{22}=1$, $b_{12}=0.3$, $b_{21}=0.25$, $\sigma_1=0.35$, $\sigma_2=0.30$, $X_0=(0.8,1.0)^\top$,
		and set $x^\circ=X_0$. The projection domain is
		\begin{flalign*}
			\mathcal D_h
			=
			\left\{
			x\in\mathbb R^2:
			x_1,x_2\ge h^\alpha,\quad
			|x-x^\circ|\le h^{-\beta}
			\right\},
		\end{flalign*}
		with $\alpha=0.25$ and $\beta=0.25$. Figure~\ref{LV} shows that the strong convergence order of the projected Milstein is consistent with the expected order. 
		
		\begin{figure}[htbp]
			\centering
			\includegraphics[width=0.8\textwidth]{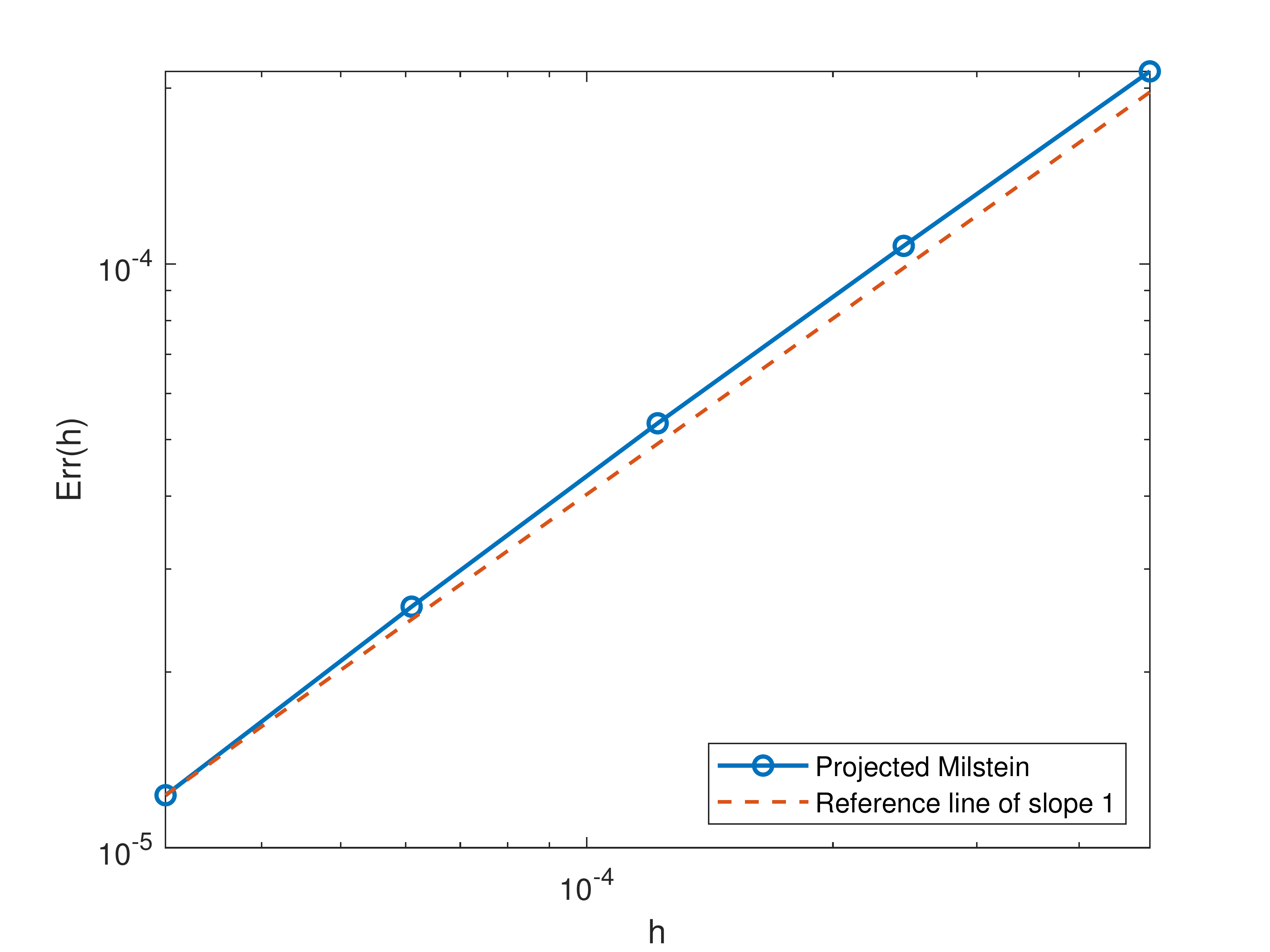}
			\caption{Competitive stochastic LV system}
			\label{LV}
		\end{figure}
	\end{example}
	
	\begin{example}{\textbf{Two-island WF diffusion}}
		
		We consider the two-island Wright--Fisher diffusion with $\tilde{u}_1=\tilde{u}_2=6$, $\tilde{\vartheta}=0.4$, $X_0=(0.3,0.7)^\top$. The projection set is 
		\begin{flalign*}
			\mathcal D_h = [h^\alpha,1-h^\alpha]^2,
		\end{flalign*}
		with $\alpha=0.3$. The strong convergence order of the projected Milstein is close to one, as shown in Figure \ref{TWF}. 
		
		\begin{figure}[htbp]
			\centering
			\includegraphics[width=0.8\textwidth]{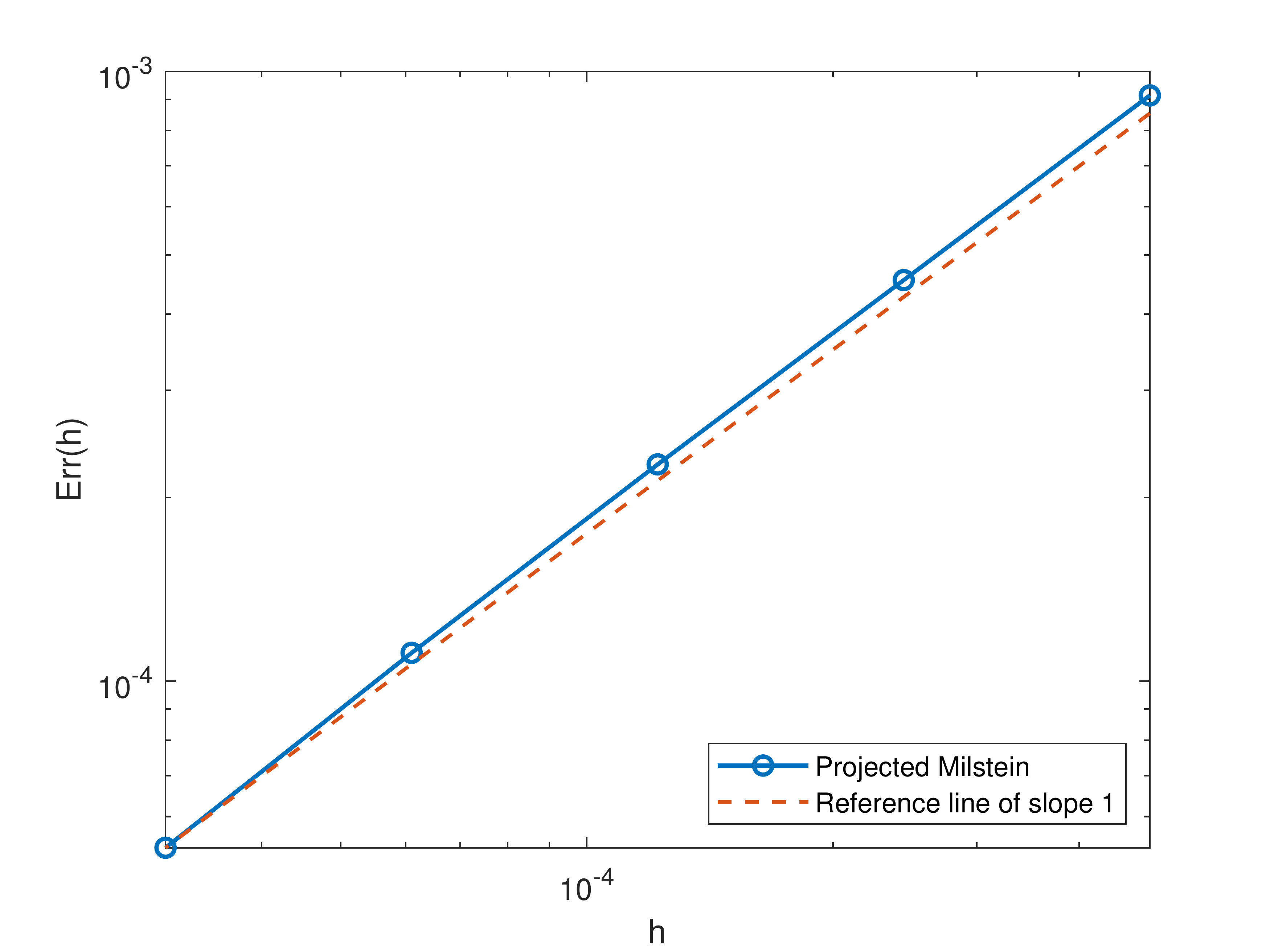}
			\caption{Two-island WF diffusion}
			\label{TWF}
		\end{figure}
		
	\end{example}

	\begin{example}{\textbf{Stochastic SIR model with population-dependent noise}}
		
		We take $\tilde K=2$, $a=0.5$, $b=0.3$,
		$c=0.4$, $\tilde\delta=0.7$ and $\sigma=0.2$,
		with $X_0=(0.5,0.2,0.2)^\top$. The projection domain is
		\begin{flalign*}
			\mathcal D_h =
			\left\{ (S,I,R)^\top: S,I,R\ge h^\alpha,\quad S+I+R\le \tilde{K}-\sqrt{3}\,h^\alpha\right\},
		\end{flalign*}
		with $\alpha=0.45$. Figure~\ref{NSIR} shows that the strong convergence order of the projected Milstein is close to $1$. 
		
		\begin{figure}[htbp]
			\centering
			\includegraphics[width=0.8\textwidth]{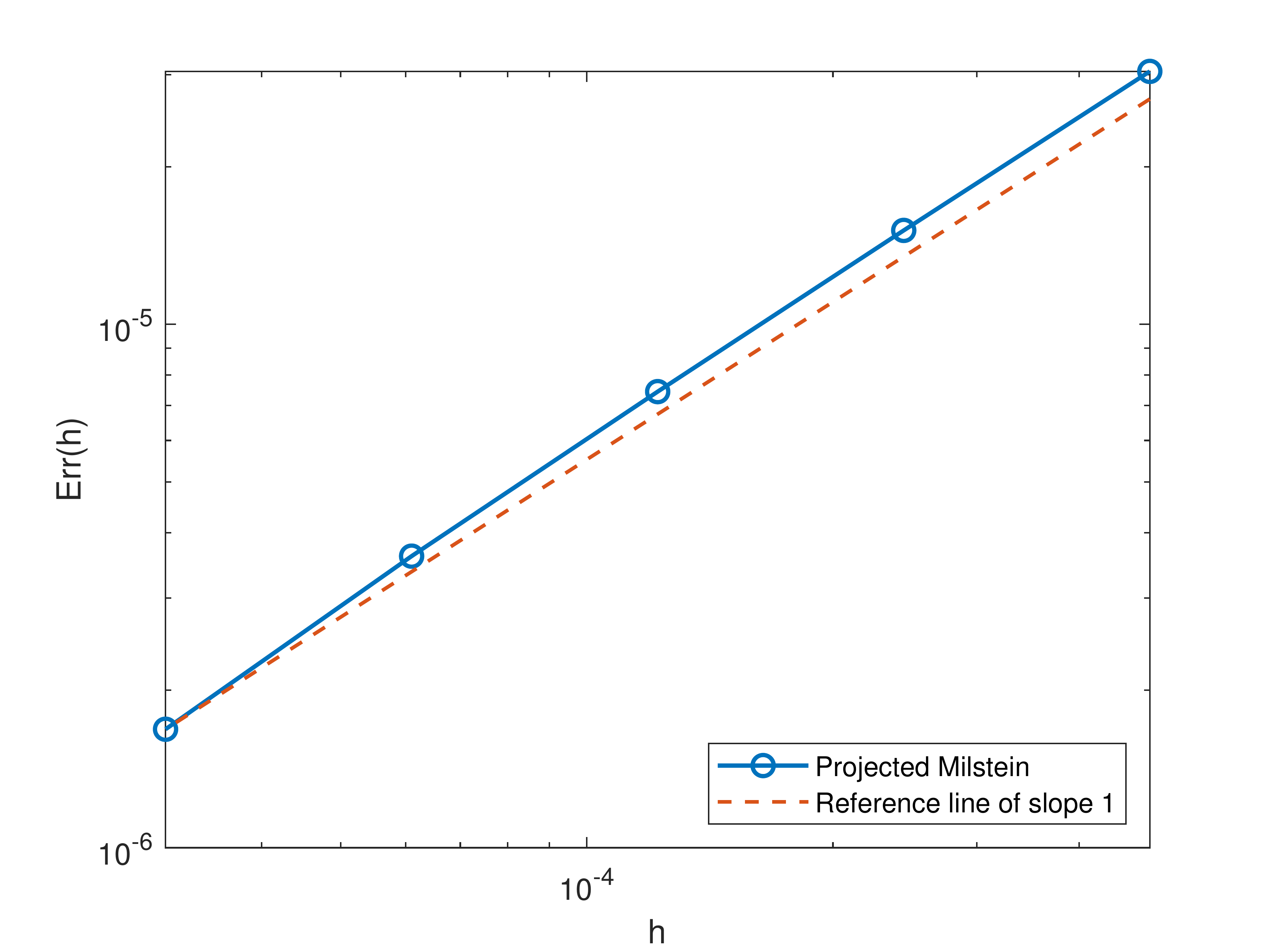}
			\caption{Stochastic SIR model with population-dependent noise}
			\label{NSIR}
		\end{figure}
	\end{example}
	
	The strong errors for all five examples are summarized in Table~\ref{tab:strong_errors_all}. For each model, the reported convergence order corresponds to the slope of a least-squares fit of $\log(\mathrm{Err}(h))$ versus $\log h$ across the five tested step sizes. All fitted orders are close to one, in agreement with the strong order-one convergence predicted by Theorem~\ref{main-thm}.
	
	\begin{table}[htbp]
		\centering
		\caption{Strong errors and convergence orders for the five test models.}
		\label{tab:strong_errors_all}
		\begin{tabular}{c|ccccc}
			\hline
			
			$h$
			& SIS
			& SIRS
			& Lotka--Volterra
			& Wright--Fisher
			& SIR \\
			\hline
			
			$2^{-11}$
			& $2.6846\times10^{-5}$
			& $3.9481\times10^{-5}$
			& $2.1364\times10^{-4}$
			& $9.1329\times10^{-4}$
			& $3.0432\times10^{-5}$ \\
			
			$2^{-12}$
			& $1.3371\times10^{-5}$
			& $1.9596\times10^{-5}$
			& $1.0736\times10^{-4}$
			& $4.5479\times10^{-4}$
			& $1.5119\times10^{-5}$ \\
			
			$2^{-13}$
			& $6.5763\times10^{-6}$
			& $9.6467\times10^{-6}$
			& $5.3332\times10^{-5}$
			& $2.2678\times10^{-4}$
			& $7.4384\times10^{-6}$ \\
			
			$2^{-14}$
			& $3.2088\times10^{-6}$
			& $4.6700\times10^{-6}$
			& $2.5875\times10^{-5}$
			& $1.1133\times10^{-4}$
			& $3.6069\times10^{-6}$ \\
			
			$2^{-15}$
			& $1.5152\times10^{-6}$
			& $2.1808\times10^{-6}$
			& $1.2301\times10^{-5}$
			& $5.3338\times10^{-5}$
			& $1.6842\times10^{-6}$ \\
			\hline
			Order
			& $1.0353$
			& $1.0426$
			& $1.0289$
			& $1.0226$
			& $1.0543$ \\
			\hline
		\end{tabular}
	\end{table}
	
	Finally, we compare the proposed projected Milstein scheme with the first-order explicit domain-preserving scheme of ~\cite{br26}. Since the strong order-one scheme in~\cite{br26} is developed for a scalar Brownian driver and a hypercube-type invariant domain, we perform the comparison on a common test problem satisfying
	the assumptions of both schemes.
	\begin{example}{\textbf{Comparison on a two-dimensional test problem}}
		
		Consider the two-dimensional SDE on $\mathcal D=(-1,1)^2$,
		\begin{equation*}
			\left\{
			\begin{aligned}
				\mathrm dX_{1,t} &= (X_{1,t}^2-1)X_{1,t}X_{2,t}\,\mathrm dt +(X_{1,t}^2-1)X_{1,t}^2X_{2,t}\,\mathrm dB_t,
				\nonumber\\ \mathrm dX_{2,t} &=
				(X_{2,t}^2-1)\cos(X_{1,t})\sin(X_{2,t})\,\mathrm dt
				+(X_{2,t}^2-1)X_{1,t}X_{2,t}^2\,\mathrm dB_t,
			\end{aligned}
			\right.
		\end{equation*}
		with $X_0=(0.25,-0.25)^\top$. Both the drift and diffusion components vanish on the corresponding faces $x_i=\pm1$, so that the test problem is compatible with the hypercube framework in~\cite{br26}. For the proposed scheme, we use
		\begin{flalign*}
			\mathcal D_h  = [-1+h^\alpha,1-h^\alpha]^2,
		\end{flalign*}
		with $\alpha=0.4$. We compare the proposed projected Milstein scheme with the DP$(1.0;2)$ scheme of~\cite{br26}. For each scheme, the order reported at step size $h$ in Table~\ref{tab:BC-comparison} is computed as $\log_2\!\bigl(\mathrm{Err}(2h)/\mathrm{Err}(h)\bigr)$.
		As shown in Figure~\ref{Compare}, both schemes exhibit approximately first-order strong convergence in the discrete-time $L^2$ maximum norm. However, the convergence theorem in~\cite{br26} establishes an error bound for the maximum of the pointwise $L^2$ errors.
		
		\begin{figure}[htbp]
			\centering
			\includegraphics[width=0.8\textwidth]{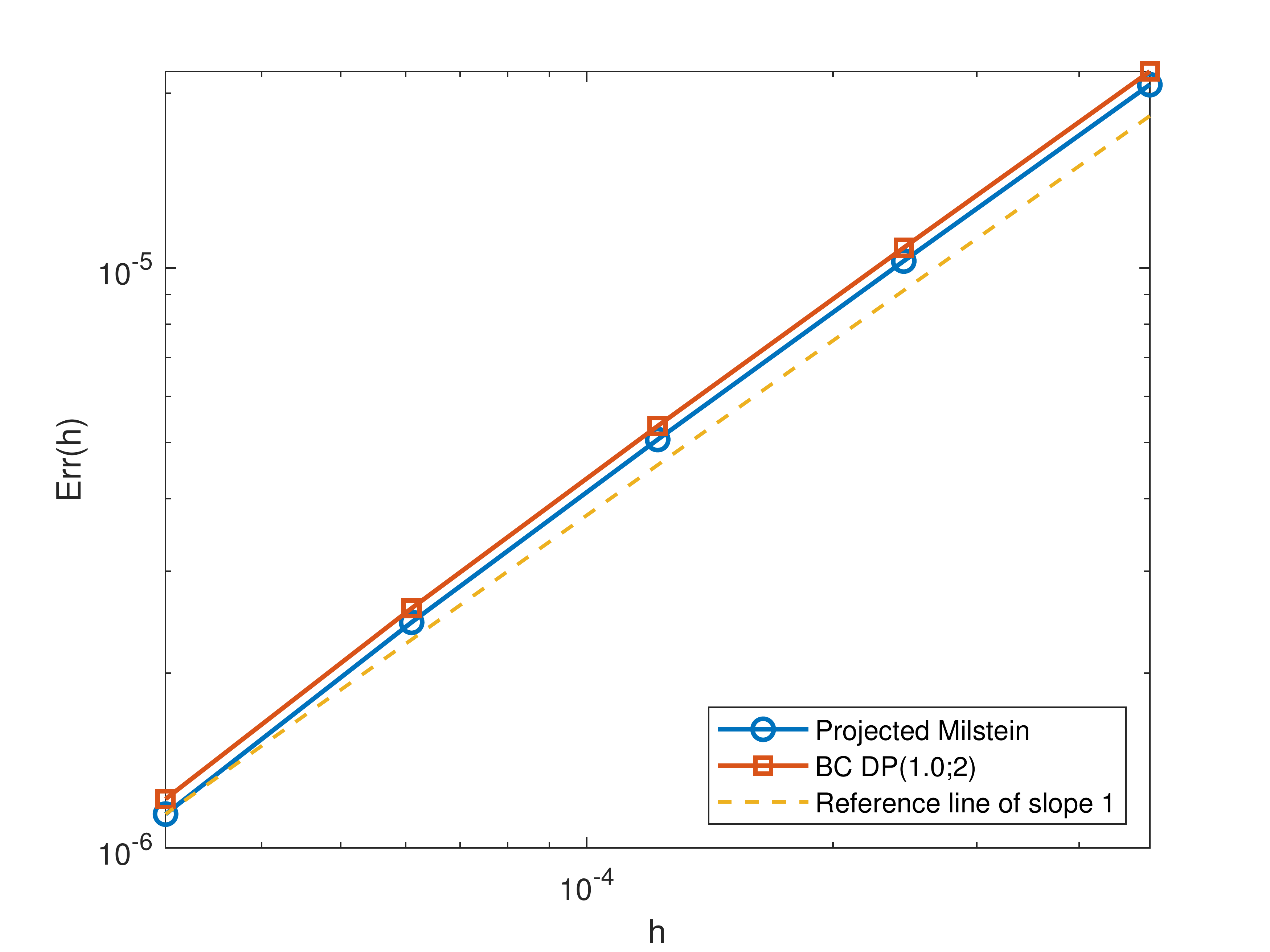}
			\caption{Comparison on a two-dimensional test problem}
			\label{Compare}
		\end{figure}
		
		\begin{table}[htbp]
			\centering
			\caption{Strong convergence of the projected Milstein scheme and the DP$(1.0;2)$ scheme of Br\'ehier and Cohen.}
			\label{tab:BC-comparison}
			\begin{tabular}{ccccc}
				\hline
				& \multicolumn{2}{c}{Projected Milstein}
				& \multicolumn{2}{c}{DP$(1.0;2)$}\\
				\cline{2-3}\cline{4-5}
				$h$
				& Error & Order
				& Error & Order\\
				\hline
				
				$2^{-11}$
				& $2.0700\times10^{-5}$ & --
				& $2.1821\times10^{-5}$ & --\\
				
				$2^{-12}$
				& $1.0276\times10^{-5}$ & $1.0104$
				& $1.0836\times10^{-5}$ & $1.0099$\\
				
				$2^{-13}$
				& $5.0575\times10^{-6}$ & $1.0227$
				& $5.3374\times10^{-6}$ & $1.0216$\\
				
				$2^{-14}$
				& $2.4479\times10^{-6}$ & $1.0469$
				& $2.5883\times10^{-6}$ & $1.0441$\\
				
				$2^{-15}$
				& $1.1429\times10^{-6}$ & $1.0989$
				& $1.2134\times10^{-6}$ & $1.0930$\\
				
				\hline
			\end{tabular}
		\end{table}
	\end{example}
	
	The comparison has two purposes. First, it checks that the proposed projection does not destroy the first-order Milstein accuracy on a bounded product domain. Second, it provides a direct computational comparison with an existing explicit domain-preserving first-order scheme on a problem for which both approaches are theoretically applicable. For the non-product and unbounded domains considered in the preceding subsections, such a direct comparison is not appropriate, because those geometries fall outside the basic hypercube setting of~\cite{br26}.

	\FloatBarrier
	\section{Conclusion}\label{sec:conclusion}
	\label{sec6}
	Under the stated assumptions and parameter conditions,
	the domain-preserving projected Milstein scheme achieves
	strong order-one convergence for the maximum grid error
	in $L^p$, for suitable $p>1$ and sufficiently small step sizes.
	The proof combines integrability estimates for the exact
	solution with local Milstein consistency, localization
	and weighted stability. Numerical experiments show
	approximately first-order convergence over the tested
	step sizes.
	
	Future work will focus on the weak convergence analysis of the proposed scheme and the development of higher-order domain-preserving schemes.

	
	
	\section*{Acknowledgements}
	The authors use ChatGPT (OpenAI, GPT-6.1sol Ultra) to assist with English-language editing, enhance the clarity, organization of the manuscript and standardize the formatting of the references. The authors review and approve all AI-assisted edits incorporated into the final manuscript. A. Xiao is supported by the National Natural Science Foundation of China (No. 12471391). X. Hu is supported by the Postgraduate Scientific Research Innovation Project of Hunan Province ( No. CX20250930) and the China Scholarship Council (No. 202508430056). Y. Hu is supported by an NSERC Discovery Grant (No. RGPIN 2024-05941) and a centennial fund from the University of Alberta.
	

\end{document}